\documentclass[10pt]{article}
\usepackage{graphics}                     
\usepackage{graphicx}                     
\date{}
\newtheorem{lemma}{Lemma}[section]
\newtheorem{theorem}{Theorem}[section]
\newtheorem{proposition}{Proposition}[section]
\newcommand{\eproof}{\makebox[1cm]{}\hfill{\framebox[2.5mm]{}}}
\begin{document}
\title{\bf Snark Designs\\}
\author{A.~D.~Forbes\\
        Department of Mathematics and Statistics\\
        The Open University\\
        Walton Hall\\
        Milton Keynes MK7 6AA\\
        UNITED KINGDOM}
\maketitle
\begin{abstract}
The main aim of this paper is to solve the design spectrum problem for
Tietze's graph,
the two 18-vertex Blanu\v{s}a snarks,
the six snarks on 20 vertices (including the flower snark J5),
the twenty non-trivial snarks on 22 vertices (including the two Loupekine snarks) and
Goldberg's snark \#3.
Together with the Petersen graph (for which the spectrum has already been computed)
this list includes all non-trivial snarks of up to 22 vertices.
We also give partial results for a selection of larger graphs:
the two Celmins--Swart snarks,
the 26- and 34-vertex Blanu\v{s}a snarks,
the flower snark J7,
the double star snark,
Zamfirescu's graph,
Goldberg's snark \#5,
the Szekeres snark and the Watkins snark.
\end{abstract}

\newcommand{\adfsplit}{\\ \makebox{~~~~~~~~}}
\newcommand{\adfLgap}{\vskip 2mm}

\newcommand{\adfVfy}[1]{}

\section{Introduction}
\label{sec:introduction}
Throughout this paper all graphs are simple.
Let $G$ be a graph. If the edge set of a graph $K$ can be partitioned into edge sets
of graphs each isomorphic to $G$, we say that there exists a {\em decomposition} of $K$ into $G$.
In the case where $K$ is the complete graph $K_n$ we refer to the decomposition as a
$G$ {\em design} of order $n$. The {\em spectrum} of $G$ is the set of positive integers $n$ for
which there exists a $G$ design of order $n$.
For completeness we remark here that the empty set is a $G$ design of order 1,
and that this trivial case will be omitted from discussion henceforth.

A complete solution of the spectrum problem often seems to be very difficult.
However it has been achieved in many cases, especially amongst the smaller graphs.
We refer the reader to the survey article of Adams, Bryant and Buchanan, \cite{ABB} and,
for more up to date results, the Web site maintained by Bryant and McCourt, \cite{BM}.
In particular, as a result of recent work on the graphs of the Platonic solids,
\cite{ABFG-d}, \cite{FG-i}, the spectrum of the dodecahedron is now known
and for the icosahedron the spectrum has been determined up to a small number of unresolved cases.
The Archimedean solids have also received recent attention, \cite{FGH-ttoc}, \cite{FG-arch},
and we now have known spectra for five of the graphs and partial results for a further seven.
As far as the author is aware, the current state of knowledge is as follows.
\begin{enumerate}
\item Tetrahedron designs of order $n$ exist if and only if
$n \equiv$ 1 or 4 (mod 12), \cite{Han}.
\item Octahedron designs of order $n$ exist if and only if
$n \equiv$ 1 or 9 (mod 24), $n \neq 9$, \cite{GRR}, \cite{ABR}.
\item Cube designs of order $n$ exist if and only if
$n \equiv$ 1 or 16 (mod 24), \cite{Mah}, \cite{Kot}, \cite{BZG}.
\item Dodecahedron designs of order $n$ exist if and only if
$n \equiv$ 1, 16, 25 or 40 (mod 60) and $n \neq 16$, \cite{AB-tocdi}, \cite{ABB}, \cite{ABFG-d}.
\item Icosahedron designs of order $n$ exist if
$n \equiv$ 1, 16, 21 or 36 (mod 60) except possibly $n = 21$, $141$, $156$, $201$, $261$, $276$;
\cite{AB-tocdi}, \cite{ABB}, \cite{FG-i}.
\item Cuboctahedron designs of order $n$ exist if and only if
$n \equiv$ 1 or 33 (mod 48), \cite{GGH-co}.
\item Rhombicuboctahedron designs of order $n$ exist if and only if
$n \equiv$ 1 or 33 (mod 96), \cite{FGH-rco}.
\item Truncated tetrahedron design of order $n$ exist if and only if
$n \equiv 1$ or 28 (mod 36), \cite{FGH-ttoc}.
\item Truncated octahedron designs of order $n$ exist if and only if
$n \equiv 1$ or 64 (mod 72), \cite{FGH-ttoc}.
\item Truncated cube designs of order $n$ exist if and only if
$n \equiv 1$ or 64 (mod 72), \cite{FGH-ttoc}.
\item Truncated cuboctahedron, icosidodecahedron, rhombicosidodecahedron,
truncated icosahedron, truncated dodecahedron, truncated icosidodecahedron and
snub cube designs of order $n$ exist if $n \equiv 1$ \textup{(mod $2e$)},
where $e$ is the number of edges in the Archimedean graph \cite{FG-arch}.
\item There exists an icosidodecahedron design of order $81$
and a truncated dodecahedron design of order $145$, \cite{FG-arch}.
\end{enumerate}

A {\em snark} is a connected, bridgeless 3-regular graph with chromatic index 4.
However, a snark is usually regarded as trivial (or reducible) if it has girth less than 5 or if it has
three edges the deletion of which results in a disconnected graph each of whose components
is non-trivial (as a graph).
The term is due to Martin Gardner who, in a popular account of the subject, \cite{Gard}, likened
the elusive nature of these graphs to that of the quarry in Lewis Carroll's poem, \cite{Carr}.
It appears that the only non-trivial snark where the spectrum has been found is the
smallest, namely the Petersen graph, \cite{AB-pet}.
The main purpose of this paper is to determine the spectrum for each of the remaining non-trivial snarks
with up to 22 vertices: the two Blanu\v{s}a snarks on 18 vertices,
the six snarks on 20 vertices and the twenty non-trivial snarks on 22 vertices
(Theorems~\ref{thm:Blanusa18 18}, \ref{thm:snark20 20} and \ref{thm:snark22 22}).
There are no non-trivial snarks on 12, 14 or 16 vertices.
In addition, we compute the spectra of two trivial snarks,
the 12-vertex Tietze's graph (Theorem~\ref{thm:Tietze 12}) and
the 24-vertex Goldberg's snark \#3 (Theorem~\ref{thm:Goldberg3 24}),
and we report partial results for some of the larger snarks (Theorem~\ref{thm:Snark partial}).

Let $v$ denote the number of vertices of $G$.
It is clear that for a 3-regular graph $G$, a $G$ design of order $n$ can exist only if
(i) $n = 1$ or $n \ge v$,
(ii) $n(n-1) \equiv 0 ~(\mathrm{mod~} 3v)$ and
(iii) $n \equiv 1 ~(\mathrm{mod~} 3)$.
These conditions are determined by elementary counting and given explicitly in
Table~\ref{tab:necessary} for some values of $v$.
\begin{table}
\caption{Design existence conditions for 3-regular graphs}
\label{tab:necessary}
\begin{center}
\begin{tabular}{@{}c|l|l@{}}
 $v$     & typical graphs               & conditions\\
\hline
  10     & Petersen graph               & $n \equiv 1$ or 10 (mod 15) \rule{0mm}{5mm}\\
  12     & Tietze graph                 & $n \equiv 1$ or 28 (mod 36) \\
  18     & Blanu\v{s}a snarks           & $n \equiv 1$ (mod 27) \\
  20     & Flower snark J5              & $n \equiv 1,$ 16, 25, 40 (mod 60), $n \neq 16$\\
  22     & Loupekine snarks             & $n \equiv 1$ or 22 (mod 33) \\
  24     & Goldberg's snark \#3         & $n \equiv 1$ or 64 (mod 72) \\
  26     & $\begin{array}{l}\textrm{Celmins--Swart~snarks}\\
                            \textrm{Blanu\v{s}a~snarks}
            \end{array}$                & $n \equiv 1$ or 13 (mod 39), $n \neq 13$\\
  28     & Flower snark J7              & $n \equiv 1,$ 28, 49, 64 (mod 84)\\
  30     & Double star snark            & $n \equiv 1$ or 10 (mod 45), $n \neq 10$\\
  34     & Blanu\v{s}a snarks           & $n \equiv 1$ or 34 (mod 51) \\
  36     & Zamfirescu's graph           & $n \equiv 1$ or 28 (mod 108), $n \neq 28$\\
  40     & Goldberg's snark \#5         & $n \equiv 1,$ 16, 25, 40 (mod 120), $n \neq 16$, 25\\
  50     & Szekeres \& Watkins snarks   & $n \equiv 1$ or 25 (mod 75), $n \neq 25$
\end{tabular}
\end{center}
\end{table}
\noindent We now state our results.

{\theorem \label{thm:Tietze 12}
Designs of order $n$ exist for the Tietze graph if and only if
$n \equiv 1$ or $28$ \textup{(mod 36)}.}

{\theorem \label{thm:Blanusa18 18}
Designs of order $n$ exist for each of the two 18-vertex Blanu\v{s}a snarks if and only if
$n \equiv 1$ \textup{(mod 27)}.}

{\theorem \label{thm:snark20 20}
Designs of order $n$ exist for each of the six snarks on $20$ vertices if and only if
$n \equiv 1,$ $16$, $25$ or $40$ \textup{(mod 60)} and $n \neq 16$.}

{\theorem \label{thm:snark22 22}
Designs of order $n$ exist for each of the twenty non-trivial snarks on $22$ vertices if and only if
$n \equiv 1$ or $22$ \textup{(mod 33)}.}

{\theorem \label{thm:Goldberg3 24}
Designs of order $n$ exist for Goldberg's snark \#3 if and only if $n \equiv 1$ or $64$ \textup{(mod 72)}.}\\

\noindent Our last theorem is a collection of partial results, with half of the possible
residue classes being resolved in each case.

{\theorem \label{thm:Snark partial}
(i) Designs of order $n$ exist for the two Celmins--Swart snarks
     and the two 26-vertex Blanu\v{s}a snarks if $n \equiv 1$ \textup{(mod 39)}.\\
(ii) Designs of order $n$ exist for the flower snark J7 if $n \equiv 1$ or $28$ \textup{(mod 84)}.\\
(iii) Designs of order $n$ exist for the double star snark if $n \equiv 1$ \textup{(mod 45)}.\\
(iv) Designs of order $n$ exist for the two 34-vertex Blanu\v{s}a snarks if $n \equiv 1$ \textup{(mod 51)}.\\
(v) Designs of order $n$ exist for Zamfirescu's graph if $n \equiv 1$ \textup{(mod 108)}.\\
(vi) Designs of order $n$ exist for Goldberg's snark \#5 if $n \equiv 1$ or $40$ \textup{(mod 120)}.\\
(vii) Designs of order $n$ exist for both the Szekeres snark and the Watkins snark if $n \equiv 1$ \textup{(mod 75)}.}\\

Our method of proof uses Wilson's fundamental construction.
For this we need the concept of a \textit{group divisible design} (GDD).
Recall that a $k$-GDD of type $g^t$ is an ordered triple ($V,\mathcal{G},\mathcal{B}$)
where $V$ is a base set of cardinality $v = tg$, $\mathcal{G}$ is a partition of $V$ into $t$
subsets of cardinality $g$ called \textit{groups} and $\mathcal{B}$ is a collection of subsets of
cardinality $k$ called \textit{blocks} which collectively have the property that each
pair of elements from different groups occurs in precisely one block but no pair of
elements from the same group occurs in any block.
We will also need $k$-GDDs of type $g^{t} h^{1}$ where $h \neq g$.
These are defined analogously, with the base set $V$ having cardinality $t g + h$
and partitioned into $t$ subsets of cardinality $g$ and one subset of cardinality $h$.
In the propositions of this section and elsewhere we give primary references for the existence
of the 3- and 4-GDDs that we use in our constructions.
The same information is conveniently presented in Ge's article in the
{\em Handbook of Combinatorial Designs}, \cite{Ge}.
As is well known, a $k$-GDD of type $g^k$ exists whenever there exist $k-2$
mutually orthogonal Latin squares (MOLS) of side $g$.

The remaining sections deal with each of the snarks stated in the theorems. We give actual
designs for a few chosen numbers as well as decompositions of certain multipartite graphs.
These are used these in combination with group divisible designs
to construct the complete graphs required to prove the theorems.
The decompositions were created and checked by a computer program written in the C language.
They were checked again by a simple {\sc Mathematica} program. A further check was
made by a more sophisticated {\sc Mathematica} program which reads the text of the paper and
extracts for appropriate analysis anything that looks like a graph edge set or a graph decomposition.
It is hoped that this threefold verification, perhaps reinforced by the Bellman's assertion,
``What I tell you three times is true'', \cite{Carr}, will help persuade
the reader that the results presented in Sections~\ref{sec:first}--\ref{sec:last} are correct.

For our purposes a graph $G$ has vertices $V(G) = \{1, 2, \dots, v\}$ for some positive integer $v$
and, since there are no isolated vertices, it is defined by its edge set
$E(G) = \{\{i,j\}: i,j \in V(G), i \sim j\}$.
In the main body of the paper we specify the edge sets that we actually used in our computations.
To save space we do not provide illustrations for the graphs of 22 vertices;
the interested reader will have no trouble in confirming that our edge sets
account for the twenty isomorphically distinct non-trivial snarks of this order.
Pictures can be found in \cite{RW}.
If $G$ has $e$ edges, the numbers of occurrences of $G$ in a decomposition into $G$ of the complete graph $K_n$,
the complete $r$-partite graph $K_{n^r}$ and the complete $(r+1)$-partite graph $K_{n^r m^1}$ are
respectively
$$ \frac{n(n-1)}{2e},~~~ \frac{n^2 r (r-1)}{2e}~~~ \textrm{and}~~~ \frac{nr(n(r - 1) + 2m)}{2e}.$$

{\proposition \label{prop:2e + 1, 4e + 1, e^3 => 2et + 1}
Let $G$ be a graph with $e > 0$ edges.
Suppose there exist $G$ designs of order $2e+1$ and $4e+1$. Suppose also there exists
a decomposition of the complete tripartite graph $K_{e,e,e}$ into $G$.
Then there exist $G$ designs of order $n$ for $n \equiv 1 ~(\textup{mod}~ 2e)$.}\\

\noindent\textbf{Proof.}
There exist 3-GDDs of types $2^{3t}$ and $2^{3t+1}$ for $t \ge 1$, \cite{Han} (see also \cite{Ge}).
In each case we perform the following construction.
Replace each point of the 3-GDD by $e$ elements; that is, inflate by a factor of $e$.
Add a further point, which we denote by the symbol $\infty$.
Lay a complete graph $K_{2e + 1}$ on each of the inflated groups together with $\infty$ and
replace each block of the 3-GDD by a complete tripartite graph $K_{e,e,e}$.
This yields designs of order $6et + 1$  for $t \ge 1$ in the first case and
$6et + 2e + 1$ for $t \ge 1$ in the second case.

There exists a 3-GDD of type $6^t 4^1$ for $t \ge 3$, \cite{CHR} (see also \cite{Ge}).
Inflate by a factor of $e$, add $\infty$,
lay a complete graph $K_{6e + 1}$ (from the previous construction) or $K_{4e + 1}$
on each of the inflated groups together with $\infty$ and
replace each block of the 3-GDD by a complete tripartite graph $K_{e,e,e}$
to give designs of order $6et + 4e + 1$ for $t \ge 3$.

Finally, there exist 3-GDDs of types $2^3 4^1$ and $4^4$ (see \cite{Ge}).
Inflate by a factor of $e$, add $\infty$,
lay a complete graph $K_{2e + 1}$ or $K_{4e + 1}$
on each of the inflated groups together with $\infty$ and
replace each block of the 3-GDD by a complete tripartite graph $K_{e,e,e}$
to give designs of orders $10e + 1$ and $16e + 1$, the two orders not covered
by the previous constructions. \eproof\\

{\proposition \label{prop:2e + 1, 4e + 1, (e/3)^3 => 2et + 1}
Let $G$ be a graph with $e > 0$ edges, $3|e$.
Suppose there exist $G$ designs of order $2e+1$ and $4e+1$. Suppose also there exists
a decomposition of the complete tripartite graph $K_{e/3,e/3,e/3}$ into $G$.
Then there exist $G$ designs of order $n$ for $n \equiv 1 ~(\textup{mod}~ 2e)$.}\\

\noindent\textbf{Proof.} There exists a 3-GDD of type $3^3$ (Latin square of side 3).
Inflate by a factor of $e/3$, and replace each block by a complete tripartite graph $K_{e/3,e/3,e/3}$
to yield a complete tripartite graph $K_{e,e,e}$.
Now use Proposition \ref{prop:2e + 1, 4e + 1, e^3 => 2et + 1}.\eproof\\

{\proposition \label{prop:e + 1, 2e + 1, 4e + 1, (e/3)^3 => et + 1}
Let $G$ be a graph with $e > 0$ edges, $3|e$.
Suppose there exist $G$ designs of order $e+1$, $2e+1$, and $4e+1$. Suppose also there exists
a decomposition of the complete tripartite graph $K_{e/3,e/3,e/3}$ into $G$.
Then there exist $G$ designs of order $n$ for $n \equiv 1 ~(\textup{mod}~ e)$.}\\

\noindent\textbf{Proof.} Take a 3-GDD of type $6^t$, $t\ge 3$, \cite{Han} (see also \cite{Ge}).
Inflate by a factor of $e/3$, add a further point, $\infty$,
lay a complete graph $K_{2e+1}$ on each of the inflated groups together with $\infty$ and
replace each block of the 3-GDD by a complete tripartite graph $K_{e/3,e/3,e/3}$.
This gives a $G$ design of order $2et + 1$ for $t \ge 3$.

There exists a 3-GDD of type $3^{2t + 1}$ for $t \ge 1$, \cite{Han} (see also \cite{Ge}).
Inflate by a factor of $e/3$, add $\infty$,
lay a complete graph $K_{e+1}$ on each of the inflated groups together with $\infty$ and
replace each block of the 3-GDD by a complete tripartite graph $K_{e/3,e/3,e/3}$.
This gives a $G$ design of order $2et + e + 1$ for $t \ge 1$.\eproof\\

{\proposition \label{prop:e + 1, 2e + 1, e^3, (e/3)^4 => et + 1}
Let $G$ be a graph with $e > 0$ edges, $3|e$.
Suppose there exist $G$ designs of order
$e+1$ and $2e+1$. Suppose also there exist decompositions of the
complete tripartite graph $K_{e,e,e}$ and the
complete 4-partite graph $K_{e/3,e/3,e/3,e/3}$ into $G$.
Then there exist $G$ designs of order $n$ for $n \equiv 1 ~(\textup{mod}~ e)$.}\\

\noindent\textbf{Proof.} There exists a 4-GDD of type $6^t$ for $t\ge 5$, \cite{BSH} (see also \cite{Ge}).
Inflate by a factor of $e/3$, add $\infty$,
lay a complete graph $K_{2e+1}$ on each of the inflated groups together with $\infty$ and
replace each block by a complete 4-partite graph $K_{e/3,e/3,e/3,e/3}$.
This gives a $G$ design of order $2et + 1$ for $t \ge 5$.

There exists a 4-GDD of type $6^t 3^1$ for $t \ge 4$, \cite{GL} (see also \cite{Ge}).
Inflate by a factor of $e/3$, add $\infty$,
lay a complete graph $K_{e+1}$ or $K_{2e+1}$ on each of the inflated groups together with $\infty$ and
replace each block by a complete 4-partite graph $K_{e/3,e/3,e/3,e/3}$.
This gives a $G$ design of order $2et + e + 1$ for $t \ge 4$.

Using 4-GDDs of types $3^4$, $3^5$, $3^5 6^1$ and $3^8$ (see \cite{Ge})
we proceed as above to yield designs of order $4e+1$, $5e+1$, $7e+1$ and $8e+1$ respectively.

To deal with the remaining cases we use a 3-GDD of type $1^3$ or $2^3$ (Latin square of side 2).
Inflate by a factor of $e$, add  $\infty$,
lay a complete graph $K_{e+1}$ or $K_{2e+1}$ on each of the inflated groups together with $\infty$ and
replace each block by a complete tripartite graph $K_{e,e,e}$ to produce a design of order
$3e + 1$ or $6e + 1$.\eproof\\

{\proposition \label{prop:d=3, v=12}
Let $G$ be a 3-regular graph on $12$ vertices.
Suppose there exist $G$ designs of order $28$, $37$, $64$, $73$ and $100$.
Suppose also that there exist decompositions of the complete multipartite graphs
$K_{6,6,6}$ and $K_{3,3,3,3}$ into $G$.
Then there exists a $G$ design of order $n$ if and only if $n \equiv 1$ or $28~(\textup{mod}~ 36)$.}\\

\noindent\textbf{Proof.} The truncated tetrahedron is a 3-regular graph on 12 vertices, which is known to have
the spectrum stated in the proposition, \cite{FGH-ttoc}.
Since the constructions in \cite{FGH-ttoc} for the truncated tetrahedron use precisely the decompositions
stated in this proposition, any graph $G$ satisfying all the conditions of the proposition
has the same spectrum.\eproof\\

{\proposition \label{prop:d=3, v=20}
Let $G$ be a 3-regular graph on $20$ vertices.
Suppose there exist $G$ designs of order $25$, $40$, $61$, $76$, $85$, $121$, $160$, $181$ and $220$.
Suppose also that there exist decompositions of the complete multipartite graphs
$K_{5,5,5,5}$, $K_{10,10,10,10}$, $K_{6,6,6,6,6}$,
$K_{60,60,75}$, $K_{60,60,60,75}$, $K_{15,15,15,21}$, $K_{24,24,24,24,39}$ and $K_{24,24,24,24,24,24,60}$ into $G$.
Then there exists a $G$ design of order $n$ if and only if
$n \equiv 1,~16,~25$ or $40~(\textup{mod}~ 60)$ and $n \neq 16$.}\\

\noindent\textbf{Proof.} The dodecahedron is a 3-regular graph on 20 vertices, which is known to have
the spectrum stated in the proposition, \cite{ABFG-d}.
We therefore follow the proof given in \cite{ABFG-d}, except that the sporadic value 340
of \cite[Section 4]{ABFG-d} is handled in a different way.
For this purpose an additional decomposition into $G$, namely that of $K_{10,10,10,10}$, is required.
The other decompositions stated in the proposition correspond to those used by the constructions
described in \cite{ABFG-d}.

There exists a 4-GDD of type $4^6 10^1$, \cite{GL} (see also \cite{Ge}).
Inflate by a factor of $10$, lay a complete graph $K_{40}$ or $K_{100}$
(for the construction of the latter, see \cite[Proposition 3.3]{ABFG-d}) on each of the inflated groups and
replace each block of the 4-GDD by a complete 4-partite graph $K_{10,10,10,10}$.
This construction yields a $G$ design of order 340.\eproof\\

{\proposition \label{prop:d=3, v=22}
Let $G$ be a 3-regular graph on $22$ vertices.
Suppose there exist $G$ designs of order $22$, $34$, $55$, $67$ and $88$.
Suppose also that there exist decompositions of the complete multipartite graphs
$K_{33,33,33}$, $K_{11,11,11,11}$, $K_{22,22,22,55}$ and $K_{22,22,22,22,22,22,55}$ into $G$.
Then there exists a $G$ design of order $n$ if and only if
$n \equiv 1$ or $22~(\textup{mod}~ 33)$.}\\

\noindent\textbf{Proof.}
The `only if' part follows from the standard divisibility conditions---see Table~\ref{tab:necessary}.
Decompositions of $K_{34}$, $K_{67}$, $K_{33,33,33}$ and $K_{11,11,11,11}$ together with
Proposition~\ref{prop:e + 1, 2e + 1, e^3, (e/3)^4 => et + 1}
yield designs of order $1$ (mod $33$).

There exist 4-GDDs of type $2^{3t+1}$ for $t \ge 2$ and type $2^{3t} 5^1$ for $t \ge 3$,
\cite{BSH}, \cite{GL} (see also \cite{Ge}).
In each case, inflate by a factor of 11,
lay a complete graph $K_{22}$ or $K_{55}$ on each of the inflated groups and replace
each block by a complete 4-partite graph $K_{11,11,11,11}$.
Since decompositions exist for $K_{22}$, $K_{55}$ and $K_{11,11,11,11}$,
this construction yields designs of order 22 (mod 33) except for 88, 121 and 187.

For 121 and 187 we use the decompositions of
$K_{22}$, $K_{55}$, $K_{22,22,22,55}$ and $K_{22,22,22,22,22,22,55}$ to
construct decompositions of the complete graphs $K_{121}$ and $K_{187}$.\eproof\\

{\proposition \label{prop:d=3, v=24}
Let $G$ be a 3-regular graph on $24$ vertices.
Suppose there exist $G$ designs of order $64$, $73$, $136$ and $145$.
Suppose also that there exist decompositions of the complete multipartite graphs
$K_{12,12,12}$, $K_{24,24,15}$, $K_{72,72,63}$, $K_{24,24,24,24}$ and $K_{24,24,24,21}$ into $G$.
Then there exists a $G$ design of order $n$ if and only if
$n \equiv 1$ or $64~(\textup{mod}~ 72)$.}\\

\noindent\textbf{Proof.}
The `only if' part follows from the standard divisibility conditions---see Table~\ref{tab:necessary}.
Decompositions of $K_{73}$, $K_{145}$ and $K_{12,12,12}$ together with
Proposition~\ref{prop:2e + 1, 4e + 1, (e/3)^3 => 2et + 1}
yield designs of order $1$ (mod $72$).

For the other residue class, we first construct a decomposition of the complete tripartite
graph $K_{24,24,24}$ into $G$ using the decomposition of $K_{12,12,12}$ and a
3-GDD of type $2^3$ (Latin square of side 2), as in the proof of
Proposition~\ref{prop:2e + 1, 4e + 1, (e/3)^3 => 2et + 1}.

There exists a 4-GDD of type $6^{t} 3^1$ for $t \ge 4$, \cite{GL} (see also \cite{Ge}).
Inflate the 3-element group by a factor of 21,
inflate all other groups by a factor of 24, and add a new point, $\infty$.
Lay a complete graph $K_{64}$ over the inflated 3-element group together with $\infty$ and
lay a complete graph $K_{145}$ over each of the other inflated groups together with $\infty$.
Replace each block of the 4-GDD by one of the complete 4-partite graphs
$K_{24,24,24,24}$ or $K_{24,24,24,21}$, as appropriate.
Since decompositions exist for $K_{64}$, $K_{145}$, $K_{24,24,24,24}$ and $K_{24,24,24,21}$,
this construction yields designs of order 64 (mod 144) except for 208, 352 and 496.

There exists a 3-GDD of type $3^{2t} 9^1$ for $t \ge 2$, \cite{CHR} (see also \cite{Ge}).
Inflate the 9-element group by a factor of 15,
inflate all other groups by a factor of 24, add $\infty$,
lay a complete graph $K_{136}$ over the inflated 9-element group together with $\infty$ and
lay a complete graph $K_{73}$ over each of the other groups together with $\infty$.
Replace each block of the 3-GDD by one of the complete tripartite graphs
$K_{24,24,24}$ or $K_{24,24,15}$, as appropriate.
Since decompositions exist for $K_{73}$, $K_{136}$, $K_{24,24,24}$ and $K_{24,24,15}$,
this construction yields designs of order 136 (mod 144) except for 280.

It only remains to deal with the four cases not covered by the previous constructions.
For 208, we construct a complete graph $K_{208}$ by augmenting the complete tripartite
graph $K_{72,72,63}$ with a common point, $\infty$, and laying complete graphs
$K_{73}$ and $K_{64}$ over the augmented partitions.

There exist 4-GDDs of types $3^4$, $3^5$ and $3^5 6^1$ (see \cite{Ge}).
In each case, inflate one of the 3-element groups by a factor of 21,
inflate all other groups by a factor of 24, add $\infty$,
lay a complete graph $K_{64}$ over the 21-inflated 3-element group together with $\infty$,
lay a complete graph $K_{73}$ over each of the other inflated 3-element groups together with $\infty$
and, if present, lay a complete graph $K_{145}$ over the inflated 6-element group together with $\infty$.
Replace each block of the 3-GDD by $K_{24,24,24,24}$ or $K_{24,24,24,21}$, as appropriate.
These constructions give $G$ designs of orders 280, 352 and 496 respectively.\eproof\\


\newcommand{\adfGsmall}{60mm}
\newcommand{\adfGlarge}{70mm}

\section{The Tietze graph}
\label{sec:Tietze}
\label{sec:first}
\begin{figure}[h]
\begin{center}
\includegraphics[width=\adfGsmall]{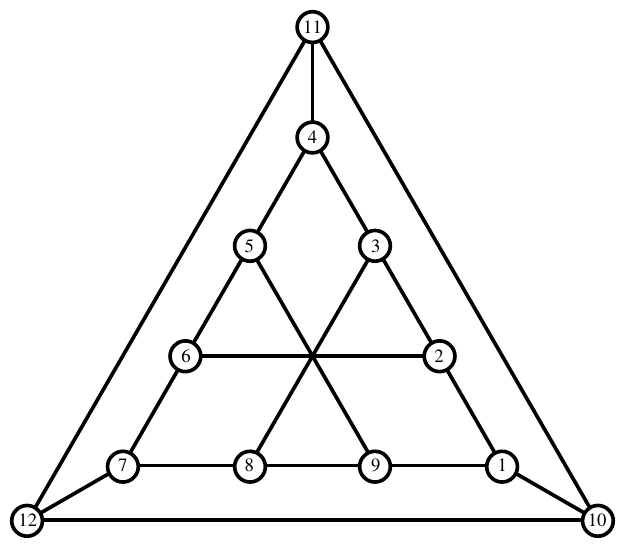}
\caption{The Tietze graph}
\end{center}
\end{figure}
\newcommand{\adfT}{\mathrm{T}}
The Tietze graph is represented by the ordered 12-tuple of its vertices (1, 2, \dots, 12)${}_{\adfT}$.
The edge set that we use, as supplied with the {\sc Mathematica} system, \cite{Math-Snarks}, is
\{$\{1,2\}$, $\{2,3\}$, $\{3,4\}$, $\{4,5\}$, $\{5,6\}$, $\{6,7\}$, $\{7,8\}$, $\{8,9\}$,
$\{1,9\}$, $\{1,10\}$, $\{4,11\}$, $\{7,12\}$, $\{3,8\}$, $\{2,6\}$, $\{5,9\}$,
$\{10,11\}$, $\{11,12\}$, $\{10,12\}$\}.

{\lemma \label{lem:Tietze designs}
There exist Tietze graph designs of order $28$, $37$, $64$, $73$ and $100$.}\\


\noindent\textbf{Proof.}
Let the vertex set of $K_{28}$ be $Z_{28}$. The decomposition consists of the graphs

$(23,26,13,6,0,1,2,3,4,7,11,18)_{\adfT}$,

$(0,2,4,1,3,8,11,13,17,10,6,19)_{\adfT}$,

$(0,5,11,1,13,20,12,24,7,18,9,26)_{\adfT}$

\noindent under the action of the mapping $x \mapsto x + 4$ (mod 28).

\adfVfy{28, \{\{28,7,4\}\}, -1, -1, -1} 

Let the vertex set of $K_{37}$ be $Z_{37}$. The decomposition consists of the graph

$(0,1,3,6,2,7,14,23,15,10,22,33)_{\adfT}$

\noindent under the action of the mapping $x \mapsto x + 1$ (mod 37).

\adfVfy{37, \{\{37,37,1\}\}, -1, -1, -1} 

Let the vertex set of $K_{64}$ be $Z_{63} \cup \{\infty\}$. The decomposition consists of the graphs

$(\infty,6,30,11,56,16,62,39,41,19,22,24)_{\adfT}$,

$(30,26,7,33,58,32,29,3,15,14,55,42)_{\adfT}$,

$(20,52,50,25,58,13,23,37,29,51,5,56)_{\adfT}$,

$(44,3,43,49,7,2,36,54,24,17,56,52)_{\adfT}$,

$(10,61,52,25,39,33,12,4,3,54,57,2)_{\adfT}$,

\vskip 1mm
$(1,2,9,43,44,51,22,23,30,17,59,38)_{\adfT}$

\noindent under the action of the mapping $\infty \mapsto \infty$, $x \mapsto x + 3$ (mod 63)
for the first five, and $x \mapsto x + 9$ (mod 63) for the last one.

\adfVfy{64, \{\{63,7,9\},\{1,1,1\}\}, 63, -1, \{5,\{\{63,3,3\},\{1,1,1\}\}\}} 

Let the vertex set of $K_{73}$ be $Z_{73}$. The decomposition consists of the graphs

$(0,1,3,6,2,7,14,22,13,10,24,36)_{\adfT}$,

$(0,15,31,1,18,51,26,70,38,23,47,68)_{\adfT}$

\noindent under the action of the mapping $x \mapsto x + 1$ (mod 73).

\adfVfy{73, \{\{73,73,1\}\}, -1, -1, -1} 

Let the vertex set of $K_{100}$ be $Z_{100}$. The decomposition consists of the graphs

$(10,43,83,87,79,15,80,92,63,46,73,35)_{\adfT}$,

$(81,7,5,21,12,92,11,87,49,84,71,25)_{\adfT}$,

$(71,14,84,79,72,82,19,22,28,61,37,41)_{\adfT}$,

$(34,38,24,32,74,50,6,77,47,59,91,80)_{\adfT}$,

$(6,91,47,96,21,86,12,2,85,54,17,69)_{\adfT}$,

$(87,20,33,64,16,44,98,32,65,60,22,91)_{\adfT}$,

$(89,30,50,78,81,39,18,41,42,12,94,16)_{\adfT}$,

$(36,34,74,68,4,21,29,69,97,35,52,86)_{\adfT}$,

$(58,40,79,38,15,63,2,13,25,50,69,80)_{\adfT}$,

$(5,20,65,10,60,6,87,58,23,59,85,12)_{\adfT}$,

$(64,45,43,55,49,79,59,42,93,96,61,56)_{\adfT}$

\noindent under the action of the mapping $x \mapsto x + 4$ (mod 100).\eproof

\adfVfy{100, \{\{100,25,4\}\}, -1, -1, -1} 

{\lemma \label{lem:Tietze multipartite}
There exist decompositions of the complete multipartite graphs
$K_{6,6,6}$ and $K_{3,3,3,3}$ into the Tietze graph.}\\


\noindent\textbf{Proof.}
Let the vertex set of $K_{6,6,6}$ be $Z_{18}$ partitioned according to residue classes modulo 3.
The decomposition consists of the graph

$(0,1,2,3,7,11,4,9,5,8,13,6)_{\adfT}$

\noindent under the action of the mapping $x \mapsto x + 3$ (mod 18).

\adfVfy{18, \{\{18,6,3\}\}, -1, \{\{6, \{0,1,2\}\}\}, -1} 

Let the vertex set of $K_{3,3,3,3}$ be $Z_{12}$ partitioned according to residue classes modulo 4.
The decomposition consists of the graph

$(0,1,2,3,4,6,9,7,10,5,8,11)_{\adfT}$

\noindent under the action of the mapping $x \mapsto x + 4$ (mod 12).\eproof

\adfVfy{12, \{\{12,3,4\}\}, -1, \{\{3, \{0,1,2,3\}\}\}, -1} 

\vskip 2mm
Theorem~\ref{thm:Tietze 12} follows from Lemmas~\ref{lem:Tietze designs} and \ref{lem:Tietze multipartite},
and Proposition~\ref{prop:d=3, v=12}.


\section{The two 18-vertex Blanu\v{s}a snarks}
\label{sec:Blanusa18 18}
\begin{figure}[h]
\begin{center}
\includegraphics[width=\adfGsmall]{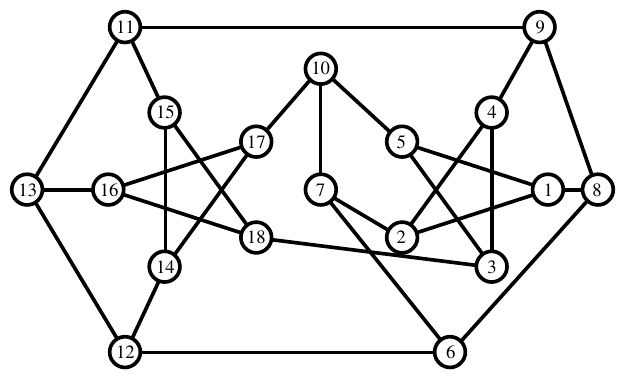}
\includegraphics[width=\adfGsmall]{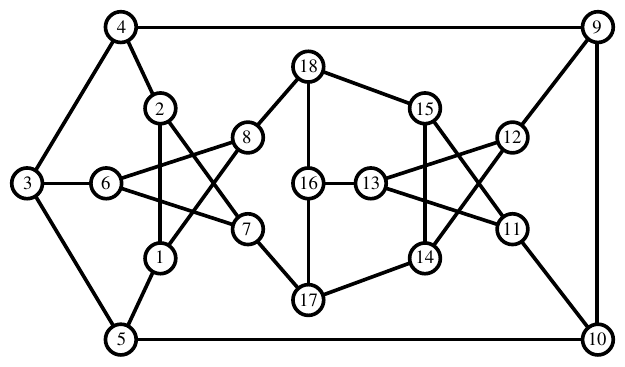}
\caption{The 18-vertex Blanu\v{s}a snarks}
\end{center}
\end{figure}
\newcommand{\adfBaa}{\mathrm{B11}}
\newcommand{\adfBab}{\mathrm{B12}}
The 18-vertex Blanu\v{s}a snarks are represented by the ordered 18-tuples
(1, 2, \dots, 18)${}_{\adfBaa}$ for the $(1,1)$-Blanu\v{s}a snark and
(1, 2, \dots, 18)${}_{\adfBab}$ for the $(1,2)$-Blanu\v{s}a snark.
The edge sets, as supplied with the {\sc Mathematica} system, \cite{Math-Snarks}, are respectively

$\adfBaa$:
\{$\{2,4\}$, $\{3,4\}$, $\{1,2\}$, $\{1,5\}$, $\{3,5\}$, $\{6,8\}$, $\{16,18\}$, $\{16,17\}$, $\{6,7\}$, $\{2,7\}$, $\{1,8\}$, $\{4,9\}$, $\{5,10\}$, $\{11,15\}$, $\{14,15\}$, $\{12,14\}$, $\{12,13\}$, $\{11,13\}$, $\{14,17\}$, $\{13,16\}$, $\{15,18\}$, $\{8,9\}$, $\{7,10\}$, $\{9,11\}$, $\{6,12\}$, $\{3,18\}$, $\{10,17\}$\}
and

$\adfBab$:
\{$\{2,4\}$, $\{3,4\}$, $\{1,2\}$, $\{1,5\}$, $\{3,5\}$, $\{3,6\}$, $\{6,8\}$, $\{8,18\}$, $\{16,18\}$, $\{16,17\}$, $\{7,17\}$, $\{6,7\}$, $\{2,7\}$, $\{1,8\}$, $\{4,9\}$, $\{9,12\}$, $\{5,10\}$, $\{10,11\}$, $\{11,15\}$, $\{14,15\}$, $\{12,14\}$, $\{12,13\}$, $\{11,13\}$, $\{14,17\}$, $\{13,16\}$, $\{15,18\}$, $\{9,10\}$\}.


{\lemma \label{lem:Blanusa18 designs}
There exist designs of order $28$, $55$ and $109$ for each of the two 18-vertex Blanu\v{s}a snarks.}\\


\noindent\textbf{Proof.}
Let the vertex set of $K_{28}$ be $Z_{28}$. The decompositions consist of

$(0,1,2,3,4,5,6,8,7,9,10,11,16,13,17,18,23,26)_{\adfBaa}$,

$(0,3,1,9,12,2,14,11,27,25,19,16,26,23,4,20,10,13)_{\adfBaa}$

\noindent and

$(0,1,2,3,4,5,6,8,7,9,13,10,19,14,20,11,21,22)_{\adfBab}$,

$(0,3,1,9,10,11,20,16,24,23,7,2,14,12,21,22,27,5)_{\adfBab}$

\noindent under the action of the mapping $x \mapsto x + 4$ (mod 28).

\adfVfy{28, \{\{28,7,4\}\}, -1, -1, -1} 

Let the vertex set of $K_{55}$ be $Z_{55}$. The decompositions consist of

$(0,1,2,4,6,3,8,11,20,18,5,21,38,53,25,13,39,44)_{\adfBaa}$

\noindent and

$(0,1,2,4,6,7,14,15,13,23,3,24,36,8,22,12,37,49)_{\adfBab}$

\noindent under the action of the mapping $x \mapsto x + 1$ (mod 55).

\adfVfy{55, \{\{55,55,1\}\}, -1, -1, -1} 

Let the vertex set of $K_{109}$ be $Z_{109}$. The decompositions consist of

$(87,66,38,101,5,88,85,75,3,86,61,71,0,1,6,2,8,16)_{\adfBaa}$,

$(0,4,1,19,9,2,27,26,55,59,3,81,47,28,92,7,75,50)_{\adfBaa}$

\noindent and

$(36,48,25,43,35,78,76,2,64,58,89,93,0,1,4,7,15,29)_{\adfBab}$,

$(0,4,1,36,37,10,23,40,21,63,3,62,73,8,55,28,74,99)_{\adfBab}$

\noindent under the action of the mapping $x \mapsto x + 1$ (mod 109).\eproof

\adfVfy{109, \{\{109,109,1\}\}, -1, -1, -1} 

{\lemma \label{lem:Blanusa18 multipartite}
There exist decompositions of $K_{9,9,9}$ into each of the two 18-vertex Blanu\v{s}a snarks.}\\


\noindent\textbf{Proof.}
Let the vertex set of $K_{9,9,9}$ be $Z_{27}$ partitioned according to residue classes modulo 3.
The decompositions consist of

$(0,1,2,3,4,5,12,19,20,23,24,10,26,14,7,16,9,21)_{\adfBaa}$

\noindent and

$(0,1,2,3,4,6,8,11,10,18,5,21,16,14,24,17,22,7)_{\adfBab}$

\noindent under the action of the mapping $x \mapsto x + 3$ (mod 27).\eproof

\adfVfy{27, \{\{27,9,3\}\}, -1, \{\{9, \{0,1,2\}\}\}, -1} 

\vskip 2mm
Theorem~\ref{thm:Blanusa18 18} follows from
Lemmas~\ref{lem:Blanusa18 designs} and \ref{lem:Blanusa18 multipartite}, and
Proposition~\ref{prop:e + 1, 2e + 1, 4e + 1, (e/3)^3 => et + 1}.


\section{The six snarks on 20 vertices}
\label{sec:snark20}
\begin{figure}[h]
\begin{center}
\includegraphics[width=\adfGsmall]{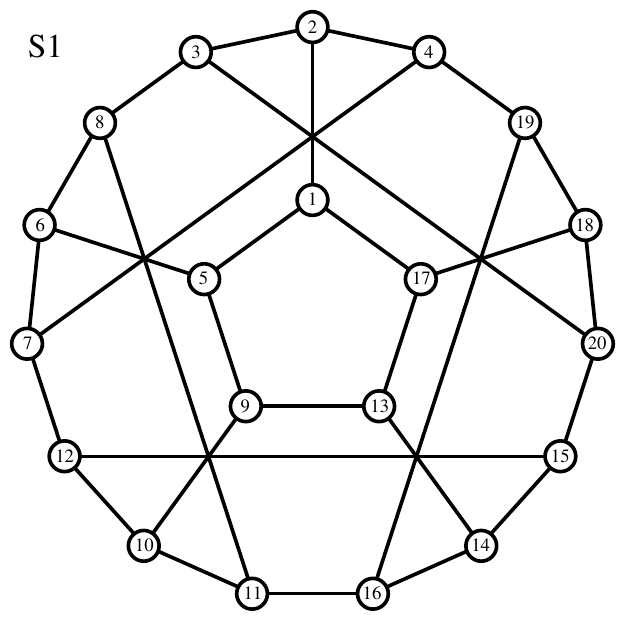}
\includegraphics[width=\adfGsmall]{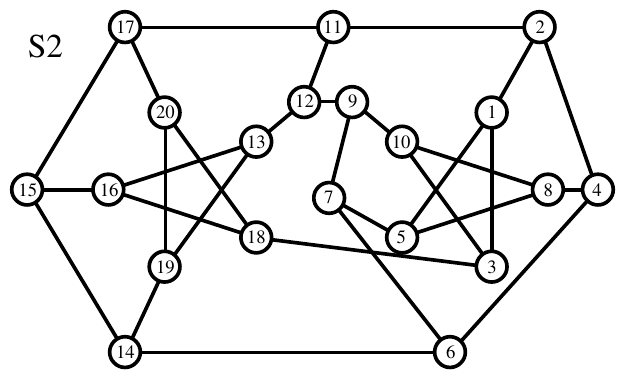}
\caption{20-vertex non-trivial snarks}
\end{center}
\end{figure}
\begin{figure}[h]
\begin{center}
\includegraphics[width=\adfGsmall]{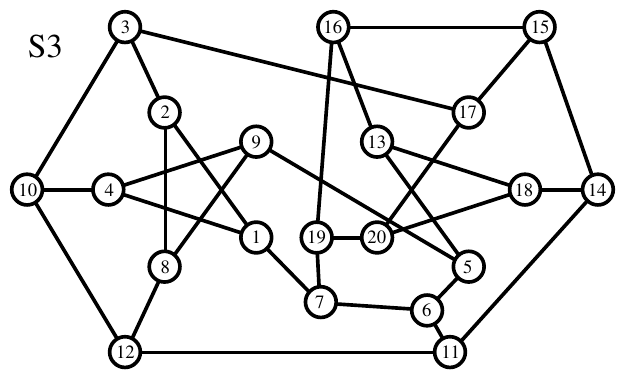}
\includegraphics[width=\adfGsmall]{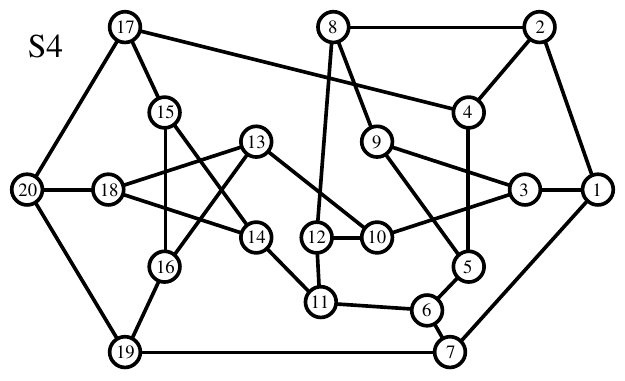}
\caption{20-vertex non-trivial snarks}
\end{center}
\end{figure}
\begin{figure}[h]
\begin{center}
\includegraphics[width=\adfGsmall]{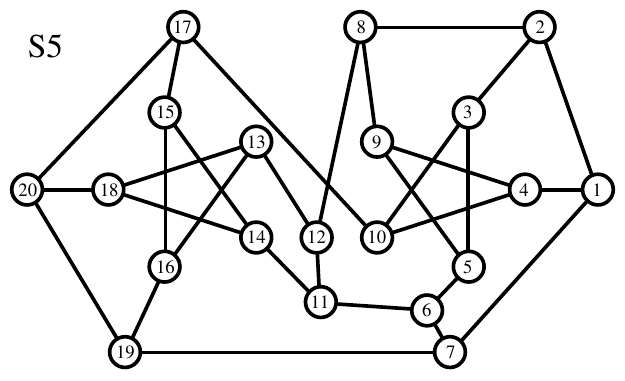}
\includegraphics[width=\adfGsmall]{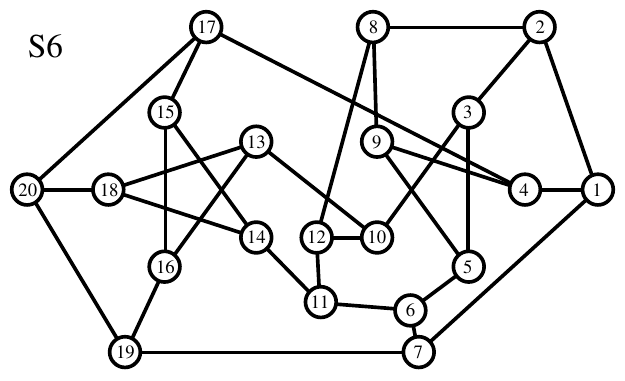}
\caption{20-vertex non-trivial snarks}
\end{center}
\end{figure}
\newcommand{\adfS}[1]{{\mathrm{S}#1}}
The six snarks on 20 vertices are represented by the ordered 20-tuples
(1, 2, \dots, 20)${}_{\adfS{1}}$, (1, 2, \dots, 20)${}_{\adfS{2}}$, \dots, (1, 2, \dots, 20)${}_{\adfS{6}}$, with edge sets

$\adfS{1}$:
   \{$\{1,2\}$, $\{2,3\}$, $\{2,4\}$, $\{5,6\}$, $\{6,7\}$, $\{6,8\}$, $\{9,10\}$, $\{10,11\}$, $\{10,12\}$, $\{13,14\}$, $\{14,15\}$, $\{14,16\}$,
   $\{17,18\}$, $\{18,20\}$, $\{18,19\}$, $\{1,17\}$, $\{13,17\}$, $\{5,9\}$, $\{9,13\}$, $\{1,5\}$, $\{3,20\}$, $\{16,19\}$, $\{4,7\}$,
   $\{8,11\}$, $\{12,15\}$, $\{3,8\}$, $\{7,12\}$, $\{11,16\}$, $\{15,20\}$, $\{4,19\}$\},

$\adfS{2}$:
  \{$\{1,2\}$, $\{1,3\}$, $\{1,5\}$, $\{2,4\}$, $\{2,11\}$, $\{3,10\}$, $\{3,18\}$, $\{4,6\}$, $\{4,8\}$, $\{5,7\}$, $\{5,8\}$, $\{6,7\}$,
   $\{6,14\}$, $\{7,9\}$, $\{8,10\}$, $\{9,10\}$, $\{9,12\}$, $\{11,12\}$, $\{11,17\}$, $\{12,13\}$, $\{13,16\}$, $\{13,19\}$, $\{14,15\}$,
   $\{14,19\}$, $\{15,16\}$, $\{15,17\}$, $\{16,18\}$, $\{17,20\}$, $\{18,20\}$, $\{19,20\}$\},

$\adfS{3}$:
  \{$\{1,2\}$, $\{1,4\}$, $\{1,7\}$, $\{2,3\}$, $\{2,8\}$, $\{3,10\}$, $\{3,17\}$, $\{4,9\}$, $\{4,10\}$, $\{5,6\}$, $\{5,9\}$, $\{5,13\}$,
   $\{6,7\}$, $\{6,11\}$, $\{7,19\}$, $\{8,9\}$, $\{8,12\}$, $\{10,12\}$, $\{11,12\}$, $\{11,14\}$, $\{13,16\}$, $\{13,18\}$, $\{14,15\}$,
   $\{14,18\}$, $\{15,16\}$, $\{15,17\}$, $\{16,19\}$, $\{17,20\}$, $\{18,20\}$, $\{19,20\}$\},

$\adfS{4}$:
  \{$\{1,2\}$, $\{1,3\}$, $\{1,7\}$, $\{2,4\}$, $\{2,8\}$, $\{3,9\}$, $\{3,10\}$, $\{4,5\}$, $\{4,17\}$, $\{5,6\}$, $\{5,9\}$, $\{6,7\}$,
   $\{6,11\}$, $\{7,19\}$, $\{8,9\}$, $\{8,12\}$, $\{10,12\}$, $\{10,13\}$, $\{11,12\}$, $\{11,14\}$, $\{13,16\}$, $\{13,18\}$,
   $\{14,15\}$, $\{14,18\}$, $\{15,16\}$, $\{15,17\}$, $\{16,19\}$, $\{17,20\}$, $\{18,20\}$, $\{19,20\}$\},

$\adfS{5}$:
  \{$\{1,2\}$, $\{1,4\}$, $\{1,7\}$, $\{2,3\}$, $\{2,8\}$, $\{3,5\}$, $\{3,10\}$, $\{4,9\}$, $\{4,10\}$, $\{5,6\}$, $\{5,9\}$, $\{6,7\}$,
   $\{6,11\}$, $\{7,19\}$, $\{8,9\}$, $\{8,12\}$, $\{10,17\}$, $\{11,12\}$, $\{11,14\}$, $\{12,13\}$, $\{13,16\}$, $\{13,18\}$,
   $\{14,15\}$, $\{14,18\}$, $\{15,16\}$, $\{15,17\}$, $\{16,19\}$, $\{17,20\}$, $\{18,20\}$, $\{19,20\}$\}
and

$\adfS{6}$:
  \{$\{1,2\}$, $\{1,4\}$, $\{1,7\}$, $\{2,3\}$, $\{2,8\}$, $\{3,5\}$, $\{3,10\}$, $\{4,9\}$, $\{4,17\}$, $\{5,6\}$, $\{5,9\}$, $\{6,7\}$,
   $\{6,11\}$, $\{7,19\}$, $\{8,9\}$, $\{8,12\}$, $\{10,12\}$, $\{10,13\}$, $\{11,12\}$, $\{11,14\}$, $\{13,16\}$, $\{13,18\}$,
   $\{14,15\}$, $\{14,18\}$, $\{15,16\}$, $\{15,17\}$, $\{16,19\}$, $\{17,20\}$, $\{18,20\}$, $\{19,20\}$\}

\noindent respectively.
The first one is the flower snark J5 as supplied with the {\sc Mathematica} system, \cite{Math-Snarks}.
Graphs S2, S3, \dots, S6 are as appear in \cite{OPRS} and \cite{RW},
where the graphs are called Sn5, Sn6, \dots, Sn9 respectively.
(We have not used the 20-vertex graph labelled Sn4 in \cite{OPRS} and \cite{RW}.)
Our graphs are pairwise isomorphic to the six 20-vertex snarks in Royale's list, \cite{Roy},
in order $(6, 1, 5, 3, 2, 4)$.


{\lemma \label{lem:snark20 designs}
For each of the six snarks on $20$ vertices, there exist designs of order
$25$, $40$, $61$, $76$, $85$, $121$, $160$, $181$ and $220$.}\\


\noindent\textbf{Proof.}
Let the vertex set of $K_{25}$ be $Z_{25}$. The decompositions consist of

$(16,22,10,20,18,12,3,17,11,14,6,7,0,1,19,2,4,8,5,9)_{\adfS{1}}$,

$(0,5,1,3,9,20,13,11,23,10,16,2,14,22,24,7,19,21,8,18)_{\adfS{1}}$,

\adfLgap
$(22,16,8,15,2,7,3,9,0,6,21,12,1,4,5,10,13,14,11,18)_{\adfS{2}}$,

$(0,2,7,4,9,13,23,21,1,22,5,15,3,20,18,12,19,24,16,14)_{\adfS{2}}$,

\adfLgap
$(1,23,2,12,15,7,17,21,16,9,19,6,0,3,4,5,8,14,11,20)_{\adfS{3}}$,

$(0,2,1,3,9,24,4,5,16,18,15,23,17,13,22,10,8,12,6,14)_{\adfS{3}}$,

\adfLgap
$(11,16,19,5,14,9,3,13,15,12,23,2,0,7,1,8,10,6,4,17)_{\adfS{4}}$,

$(0,1,2,3,10,20,17,5,24,8,23,18,7,11,4,12,19,9,21,6)_{\adfS{4}}$,

\adfLgap
$(17,0,10,13,21,12,23,8,3,15,9,4,1,2,5,6,11,7,14,20)_{\adfS{5}}$,

$(0,2,1,3,23,24,9,14,10,16,22,15,4,13,19,6,8,11,17,7)_{\adfS{5}}$

\noindent and

$(22,21,11,2,5,20,9,13,4,14,23,17,1,0,7,8,10,12,3,6)_{\adfS{6}}$,

$(0,1,3,2,20,14,16,22,13,18,24,6,5,23,19,10,12,9,21,4)_{\adfS{6}}$

\noindent under the action of the mapping $x \mapsto x + 5$ (mod 25).

\adfVfy{25, \{\{25,5,5\}\}, -1, -1, -1} 

Let the vertex set of $K_{40}$ be $Z_{39} \cup \{\infty\}$. The decompositions consist of

$(\infty,26,24,38,28,15,2,18,34,10,19,22,13,3,8,11,6,1,4,9)_{\adfS{1}}$,

$(0,4,2,5,9,7,11,17,21,35,24,15,38,29,1,13,16,36,26,18)_{\adfS{1}}$,

\adfLgap
$(\infty,2,1,14,24,3,9,10,38,19,11,37,20,0,4,6,12,7,5,34)_{\adfS{2}}$,

$(0,2,7,4,8,1,12,11,22,27,24,6,33,35,14,20,37,31,3,23)_{\adfS{2}}$,

\adfLgap
$(\infty,36,30,23,34,16,7,10,14,12,26,4,19,0,1,3,6,5,2,8)_{\adfS{3}}$,

$(0,3,7,8,1,2,9,10,30,23,13,27,12,20,22,35,19,32,21,11)_{\adfS{3}}$,

\adfLgap
$(\infty,30,37,22,17,32,5,21,31,7,13,0,24,1,2,4,6,8,3,9)_{\adfS{4}}$,

$(0,4,5,1,3,11,14,8,15,24,28,2,9,33,19,37,29,20,22,38)_{\adfS{4}}$,

\adfLgap
$(\infty,22,2,33,23,14,20,21,35,3,6,25,15,0,5,1,10,11,4,19)_{\adfS{5}}$,

$(0,3,1,12,2,4,15,16,5,28,11,37,23,24,19,13,36,8,35,18)_{\adfS{5}}$

\noindent and

$(\infty,29,32,12,9,16,4,20,27,37,3,7,35,0,1,2,6,8,10,18)_{\adfS{6}}$,

$(0,2,3,5,1,4,11,6,20,12,14,32,31,37,13,17,27,21,28,10)_{\adfS{6}}$

\noindent under the action of the mapping $\infty \mapsto \infty$, $x \mapsto x + 3$ (mod 39).

\adfVfy{40, \{\{39,13,3\},\{1,1,1\}\}, 39, -1, -1} 

Let the vertex set of $K_{61}$ be $Z_{61}$. The decompositions consist of

$(34,21,56,58,43,46,14,41,44,51,0,28,1,3,7,19,9,40,52,29)_{\adfS{1}}$,

\adfLgap
$(38,40,42,7,0,22,6,25,48,11,29,58,2,1,4,24,56,41,14,49)_{\adfS{2}}$,

\adfLgap
$(36,56,34,15,40,22,20,10,45,26,23,59,2,4,7,16,44,8,3,37)_{\adfS{3}}$,

\adfLgap
$(5,28,19,7,24,54,57,41,8,39,59,37,2,4,23,48,50,14,56,21)_{\adfS{4}}$,

\adfLgap
$(12,54,44,56,7,34,45,39,55,30,5,9,0,2,23,41,25,38,53,31)_{\adfS{5}}$

\noindent and

$(3,47,29,2,34,24,16,28,43,32,12,57,1,5,11,35,13,27,39,60)_{\adfS{6}}$

\noindent under the action of the mapping $x \mapsto x + 1$ (mod 61).

\adfVfy{61, \{\{61,61,1\}\}, -1, -1, -1} 

Let the vertex set of $K_{76}$ be $Z_{76}$. The decompositions consist of

$(21,69,20,4,7,63,67,70,30,18,14,59,38,48,9,0,51,3,47,6)_{\adfS{1}}$,

$(29,66,50,33,73,26,53,25,71,52,67,60,49,54,44,56,72,21,15,9)_{\adfS{1}}$,

$(41,3,27,43,1,15,6,26,19,13,5,12,28,57,70,60,30,66,20,11)_{\adfS{1}}$,

$(26,58,63,6,33,57,44,61,16,36,31,2,43,74,21,52,55,4,0,30)_{\adfS{1}}$,

$(59,10,32,38,60,67,8,64,45,61,52,39,35,9,54,70,25,56,19,37)_{\adfS{1}}$,

\adfLgap
$(3,14,1,55,7,15,60,49,21,4,42,26,28,61,58,69,62,47,40,54)_{\adfS{2}}$,

$(20,24,69,23,9,13,4,1,37,5,54,11,68,28,48,70,18,12,21,62)_{\adfS{2}}$,

$(2,24,25,15,27,63,48,20,49,6,50,74,36,0,47,46,5,23,42,54)_{\adfS{2}}$,

$(12,48,23,42,60,8,52,63,1,37,32,59,71,51,7,15,22,6,14,5)_{\adfS{2}}$,

$(44,62,57,59,19,43,60,46,36,75,33,13,37,67,1,22,71,29,53,66)_{\adfS{2}}$,

\adfLgap
$(32,63,49,62,19,26,65,47,30,48,10,53,12,14,67,29,57,50,11,5)_{\adfS{3}}$,

$(60,7,70,73,42,72,12,29,32,18,51,40,25,14,13,6,45,48,64,50)_{\adfS{3}}$,

$(34,60,55,28,70,35,36,69,21,64,72,33,44,31,41,59,53,19,56,50)_{\adfS{3}}$,

$(46,38,50,37,20,16,24,25,67,71,60,9,10,45,74,21,59,43,6,63)_{\adfS{3}}$,

$(65,69,19,8,7,34,35,49,20,27,59,0,55,30,10,15,45,58,67,24)_{\adfS{3}}$,

\adfLgap
$(46,13,64,4,33,26,5,14,19,2,7,34,48,28,40,18,66,25,37,27)_{\adfS{4}}$,

$(14,53,41,12,48,61,45,16,72,71,66,18,44,57,23,38,11,67,74,63)_{\adfS{4}}$,

$(24,59,39,46,38,71,45,51,22,0,23,5,19,41,58,36,48,63,7,49)_{\adfS{4}}$,

$(73,61,70,36,4,53,45,5,64,9,66,16,3,14,65,25,32,59,21,19)_{\adfS{4}}$,

$(23,49,39,55,64,36,28,68,30,72,54,3,57,42,35,34,14,46,31,56)_{\adfS{4}}$,

\adfLgap
$(10,23,52,32,62,60,24,37,9,0,49,50,5,56,74,51,19,25,66,20)_{\adfS{5}}$,

$(60,23,68,45,69,26,71,33,63,72,12,25,64,50,7,55,54,52,30,46)_{\adfS{5}}$,

$(72,56,15,22,49,7,3,43,17,42,2,6,54,57,0,55,66,24,63,41)_{\adfS{5}}$,

$(50,69,2,21,61,54,13,65,71,72,20,49,73,0,28,55,37,51,74,34)_{\adfS{5}}$,

$(15,0,21,71,24,16,59,50,47,35,11,28,31,22,10,19,42,33,57,69)_{\adfS{5}}$

\noindent and

$(12,60,70,56,45,18,55,23,17,13,53,32,67,35,24,64,25,41,58,61)_{\adfS{6}}$,

$(19,3,71,57,20,2,10,27,72,8,1,14,15,24,7,5,34,49,12,67)_{\adfS{6}}$,

$(14,7,22,60,74,8,28,55,34,5,12,26,1,72,64,61,3,13,69,45)_{\adfS{6}}$,

$(20,55,8,29,62,18,22,57,24,30,49,10,21,60,73,43,59,26,3,14)_{\adfS{6}}$,

$(36,70,1,59,48,43,63,10,60,66,11,15,18,23,49,0,6,33,46,47)_{\adfS{6}}$

\noindent under the action of the mapping $x \mapsto x + 4$ (mod 76).

\adfVfy{76, \{\{76,19,4\}\}, -1, -1, -1} 

Let the vertex set of $K_{85}$ be $Z_{85}$. The decompositions consist of

$(13,40,57,34,24,1,14,16,10,6,64,27,62,52,46,69,59,77,81,4)_{\adfS{1}}$,

$(43,74,47,64,21,49,34,81,39,13,6,9,56,32,73,68,75,26,62,42)_{\adfS{1}}$,

$(15,17,82,2,53,18,27,32,26,52,61,30,71,58,55,11,28,9,20,45)_{\adfS{1}}$,

$(42,15,80,19,12,51,40,10,5,26,83,46,67,34,48,62,4,65,70,37)_{\adfS{1}}$,

$(21,22,18,64,5,83,27,38,3,42,43,80,28,79,81,76,46,75,8,59)_{\adfS{1}}$,

$(47,28,50,12,7,80,38,24,14,66,15,5,22,35,4,65,36,41,73,11)_{\adfS{1}}$,

$(71,17,79,68,0,23,59,29,44,43,53,51,84,65,48,25,45,8,56,63)_{\adfS{1}}$,

\adfLgap
$(41,83,25,65,75,13,69,29,6,53,40,51,45,64,47,60,63,77,72,46)_{\adfS{2}}$,

$(83,8,16,74,48,17,42,47,6,38,35,26,23,14,66,20,82,39,11,19)_{\adfS{2}}$,

$(36,66,1,0,40,73,80,37,2,39,30,83,46,47,81,9,52,79,45,10)_{\adfS{2}}$,

$(37,73,61,75,51,58,18,20,72,33,46,27,7,66,77,63,64,0,54,4)_{\adfS{2}}$,

$(13,47,39,77,49,82,70,53,80,72,59,27,31,76,5,14,61,40,16,17)_{\adfS{2}}$,

$(31,75,38,47,10,56,77,24,42,63,64,57,13,51,12,2,54,69,19,14)_{\adfS{2}}$,

$(36,10,26,30,64,55,63,35,83,69,48,53,38,14,58,49,71,28,75,44)_{\adfS{2}}$,

\adfLgap
$(77,41,6,42,1,53,28,33,46,5,19,47,12,32,56,37,60,15,35,78)_{\adfS{3}}$,

$(60,7,29,65,53,55,80,72,69,77,34,67,63,49,6,82,64,44,48,8)_{\adfS{3}}$,

$(18,56,13,24,15,76,9,4,60,54,73,65,34,35,19,36,31,74,67,21)_{\adfS{3}}$,

$(16,76,21,1,13,59,81,39,70,14,28,32,41,29,43,73,22,49,3,53)_{\adfS{3}}$,

$(75,78,76,42,21,44,12,56,50,82,10,83,35,72,31,20,47,34,3,23)_{\adfS{3}}$,

$(83,60,5,46,50,23,67,41,0,12,6,18,16,70,29,19,53,11,79,20)_{\adfS{3}}$,

$(39,22,36,13,47,32,45,77,18,0,24,10,21,15,72,49,74,42,7,53)_{\adfS{3}}$,

\adfLgap
$(38,42,23,61,7,46,63,44,75,0,25,4,10,60,13,59,52,47,65,20)_{\adfS{4}}$,

$(51,70,80,54,77,2,10,52,81,15,74,75,83,62,17,12,41,82,63,34)_{\adfS{4}}$,

$(76,32,82,38,36,78,23,4,63,73,43,53,6,3,11,81,61,27,19,65)_{\adfS{4}}$,

$(52,67,30,32,46,41,19,4,6,35,80,54,79,66,50,78,57,11,34,27)_{\adfS{4}}$,

$(72,4,54,80,69,26,46,62,3,40,35,83,43,48,77,25,65,70,13,16)_{\adfS{4}}$,

$(62,20,59,22,38,29,5,31,49,28,24,4,82,3,77,33,69,81,39,53)_{\adfS{4}}$,

$(50,20,76,28,6,44,56,71,39,42,14,74,31,70,63,68,16,67,81,24)_{\adfS{4}}$,

\adfLgap
$(63,74,46,58,52,66,64,43,10,11,41,24,32,42,9,2,57,59,13,0)_{\adfS{5}}$,

$(58,51,14,73,77,39,1,43,57,10,63,72,20,53,21,12,3,30,27,15)_{\adfS{5}}$,

$(1,41,7,23,0,20,31,5,6,2,9,54,61,80,14,29,66,56,45,43)_{\adfS{5}}$,

$(3,27,80,28,29,4,1,81,19,61,17,77,60,13,40,10,42,56,37,83)_{\adfS{5}}$,

$(41,49,7,12,52,58,75,47,36,38,19,37,9,42,65,73,34,23,8,64)_{\adfS{5}}$,

$(46,24,18,25,63,35,75,73,4,68,26,70,38,44,20,2,67,0,50,79)_{\adfS{5}}$,

$(28,58,7,31,10,49,76,46,15,70,36,56,14,21,80,67,79,9,32,29)_{\adfS{5}}$

\noindent and

$(59,51,38,17,78,5,9,68,62,37,11,45,40,64,72,71,49,33,69,8)_{\adfS{6}}$,

$(38,29,23,64,9,21,56,19,82,27,68,25,10,15,62,49,30,72,22,6)_{\adfS{6}}$,

$(16,7,4,9,59,41,81,36,50,15,47,61,44,35,77,23,66,80,39,70)_{\adfS{6}}$,

$(48,75,42,63,67,62,29,82,52,22,36,32,61,39,60,74,80,66,18,55)_{\adfS{6}}$,

$(7,64,19,44,12,3,58,68,61,2,4,66,50,20,21,48,6,41,77,38)_{\adfS{6}}$,

$(36,77,53,66,17,58,48,18,38,67,6,69,80,39,76,22,50,54,40,55)_{\adfS{6}}$,

$(65,10,32,68,56,48,30,38,33,79,9,74,55,75,53,11,31,3,1,50)_{\adfS{6}}$

\noindent under the action of the mapping $x \mapsto x + 5$ (mod 85).

\adfVfy{85, \{\{85,17,5\}\}, -1, -1, -1} 

Let the vertex set of $K_{121}$ be $Z_{121}$. The decompositions consist of

$(41,27,34,87,98,81,0,92,73,2,\adfsplit 118,24,1,3,4,6,5,11,19,46)_{\adfS{1}}$,

$(0,10,25,26,18,37,5,64,41,3,\adfsplit 31,34,4,48,114,1,73,6,71,68)_{\adfS{1}}$,

\adfLgap
$(76,21,82,29,15,33,55,41,36,75,\adfsplit 49,100,0,1,2,5,4,14,11,27)_{\adfS{2}}$,

$(0,14,15,31,18,1,38,56,2,80,\adfsplit 41,74,3,49,107,45,76,89,78,9)_{\adfS{2}}$,

\adfLgap
$(79,80,5,9,56,97,64,90,18,78,\adfsplit 1,35,91,37,30,44,47,41,75,73)_{\adfS{3}}$,

$(104,74,29,46,45,67,64,111,58,90,\adfsplit 47,106,6,66,42,70,50,12,41,68)_{\adfS{3}}$,

\adfLgap
$(78,119,72,76,11,69,82,32,39,61,\adfsplit 47,116,0,1,2,5,4,9,15,27)_{\adfS{4}}$,

$(0,14,15,30,1,18,38,33,54,39,\adfsplit 43,3,4,17,93,46,120,81,97,49)_{\adfS{4}}$,

\adfLgap
$(32,14,24,92,86,40,81,34,111,78,\adfsplit 11,104,0,2,1,3,5,7,15,28)_{\adfS{5}}$,

$(0,6,14,11,29,2,24,22,61,44,\adfsplit 28,59,7,114,33,75,91,71,109,26)_{\adfS{5}}$

\noindent and

$(95,69,98,15,104,103,118,41,39,114,\adfsplit 48,66,0,1,4,8,24,11,27,60)_{\adfS{6}}$,

$(0,5,17,13,3,20,42,26,53,48,\adfsplit 52,87,107,105,18,56,70,62,100,16)_{\adfS{6}}$

\noindent under the action of the mapping $x \mapsto x + 1$ (mod 121).

\adfVfy{121, \{\{121,121,1\}\}, -1, -1, -1} 

Let the vertex set of $K_{160}$ be $Z_{159} \cup \{\infty\}$. The decompositions consist of

$(\infty,65,114,84,36,18,151,19,88,96,\adfsplit 90,98,7,4,158,44,94,25,54,63)_{\adfS{1}}$,

$(46,30,155,142,1,24,118,58,106,41,\adfsplit 59,83,62,47,80,6,101,123,34,66)_{\adfS{1}}$,

$(158,154,77,13,141,149,71,28,42,84,\adfsplit 69,34,64,75,48,78,31,1,104,113)_{\adfS{1}}$,

$(113,141,22,140,46,103,29,92,146,136,\adfsplit 9,117,122,3,13,34,142,116,86,28)_{\adfS{1}}$,

$(36,158,15,94,139,112,21,119,151,35,\adfsplit 42,149,5,153,74,132,117,69,45,2)_{\adfS{1}}$,

$(17,100,14,16,145,22,41,57,2,59,\adfsplit 86,3,43,126,82,51,26,99,112,104)_{\adfS{1}}$,

$(56,5,19,68,81,131,72,24,149,137,\adfsplit 44,116,143,145,7,141,64,25,45,53)_{\adfS{1}}$,

$(149,75,112,120,13,73,7,121,129,2,\adfsplit 63,60,3,10,72,27,88,147,81,70)_{\adfS{1}}$,

\adfLgap
$(\infty,128,3,104,16,81,113,157,133,72,\adfsplit 41,74,105,141,156,9,111,84,122,78)_{\adfS{2}}$,

$(96,28,44,13,85,120,128,21,97,22,\adfsplit 76,67,91,69,155,84,102,115,108,5)_{\adfS{2}}$,

$(111,58,18,46,53,97,62,136,106,34,\adfsplit 79,112,114,133,57,76,66,13,3,29)_{\adfS{2}}$,

$(65,82,54,32,156,99,69,62,55,74,\adfsplit 105,83,149,141,20,128,129,18,76,53)_{\adfS{2}}$,

$(156,158,79,145,122,103,137,149,59,0,\adfsplit 84,113,33,6,104,151,37,42,138,142)_{\adfS{2}}$,

$(125,7,115,14,73,62,108,39,37,11,\adfsplit 44,118,15,56,155,158,52,0,88,77)_{\adfS{2}}$,

$(87,67,97,70,112,3,17,42,106,29,\adfsplit 18,40,35,21,125,30,11,51,148,109)_{\adfS{2}}$,

$(69,126,109,129,112,125,139,11,101,54,\adfsplit 86,7,67,62,20,38,143,107,22,23)_{\adfS{2}}$,

\adfLgap
$(\infty,18,107,31,28,87,80,38,129,88,\adfsplit 158,116,127,8,69,2,19,114,144,103)_{\adfS{3}}$,

$(31,123,60,67,93,152,68,90,89,85,\adfsplit 144,139,76,11,77,96,153,146,141,156)_{\adfS{3}}$,

$(13,27,11,79,135,66,94,139,63,6,\adfsplit 53,140,30,145,97,44,7,131,74,92)_{\adfS{3}}$,

$(19,71,133,82,104,64,14,39,30,103,\adfsplit 52,12,69,11,32,81,114,115,59,124)_{\adfS{3}}$,

$(26,124,97,78,62,156,34,47,85,140,\adfsplit 75,93,56,54,136,112,103,139,125,92)_{\adfS{3}}$,

$(14,48,128,110,131,71,40,118,8,125,\adfsplit 23,157,27,11,13,57,77,91,86,117)_{\adfS{3}}$,

$(37,39,7,3,82,27,75,141,133,14,\adfsplit 70,135,140,85,130,86,6,138,114,30)_{\adfS{3}}$,

$(68,85,108,111,37,124,45,56,40,152,\adfsplit 120,83,110,58,114,54,150,100,76,41)_{\adfS{3}}$,

\adfLgap
$(\infty,44,151,13,132,101,33,5,86,21,\adfsplit 134,136,6,4,90,18,73,53,108,119)_{\adfS{4}}$,

$(73,36,100,124,66,41,9,121,93,37,\adfsplit 75,32,31,141,68,76,20,115,145,154)_{\adfS{4}}$,

$(1,53,104,70,57,16,93,33,39,134,\adfsplit 109,12,72,157,5,63,124,19,88,8)_{\adfS{4}}$,

$(56,128,111,91,155,101,127,84,5,72,\adfsplit 43,110,115,143,158,85,124,133,4,12)_{\adfS{4}}$,

$(136,89,105,112,134,155,28,26,150,44,\adfsplit 152,3,63,77,90,119,14,41,60,59)_{\adfS{4}}$,

$(158,116,52,19,123,41,82,59,66,68,\adfsplit 69,119,44,65,38,88,34,79,14,0)_{\adfS{4}}$,

$(74,112,123,55,68,156,5,42,145,121,\adfsplit 71,141,120,65,46,57,114,53,3,36)_{\adfS{4}}$,

$(150,21,108,97,93,117,153,31,43,127,\adfsplit 9,55,130,133,72,88,77,81,89,129)_{\adfS{4}}$,

\adfLgap
$(\infty,146,93,138,120,123,106,108,25,67,\adfsplit 156,40,139,87,47,64,97,137,98,82)_{\adfS{5}}$,

$(54,45,117,50,134,82,89,90,145,17,\adfsplit 35,69,12,132,15,151,113,7,157,44)_{\adfS{5}}$,

$(3,47,15,118,37,17,154,106,73,122,\adfsplit 109,55,46,7,36,41,66,147,151,74)_{\adfS{5}}$,

$(58,35,63,100,91,147,14,48,87,15,\adfsplit 19,5,151,115,140,68,82,62,45,92)_{\adfS{5}}$,

$(145,51,26,133,74,47,127,61,35,108,\adfsplit 5,123,24,121,69,18,33,89,0,14)_{\adfS{5}}$,

$(129,118,111,19,18,92,158,91,122,113,\adfsplit 4,52,142,102,55,56,114,64,11,151)_{\adfS{5}}$,

$(25,65,132,110,43,22,24,106,120,82,\adfsplit 157,141,87,91,56,41,135,29,119,128)_{\adfS{5}}$,

$(141,117,157,104,128,144,129,54,125,84,\adfsplit 95,133,15,44,50,26,14,98,48,116)_{\adfS{5}}$

\noindent and

$(\infty,142,111,29,46,140,15,129,105,106,\adfsplit 148,85,67,61,80,7,95,150,54,16)_{\adfS{6}}$,

$(72,130,55,23,28,85,109,156,43,116,\adfsplit 30,102,104,146,8,66,135,101,144,38)_{\adfS{6}}$,

$(44,134,66,132,85,63,7,16,113,118,\adfsplit 23,10,158,102,74,81,64,128,110,68)_{\adfS{6}}$,

$(92,65,34,140,32,89,90,79,150,111,\adfsplit 25,91,100,77,28,145,104,136,48,137)_{\adfS{6}}$,

$(117,125,154,25,64,82,33,13,80,73,\adfsplit 59,35,43,68,103,86,136,34,40,6)_{\adfS{6}}$,

$(0,30,65,9,42,121,16,90,24,81,\adfsplit 44,138,126,115,123,143,82,105,48,99)_{\adfS{6}}$,

$(145,0,146,144,39,89,95,33,44,40,\adfsplit 157,125,135,96,46,99,48,139,14,134)_{\adfS{6}}$,

$(149,124,127,91,101,45,8,131,47,150,\adfsplit 158,3,76,21,48,146,51,65,95,120)_{\adfS{6}}$

\noindent under the action of the mapping $\infty \mapsto \infty$, $x \mapsto x + 3$ (mod 159).

\adfVfy{160, \{\{159,53,3\},\{1,1,1\}\}, 159, -1, -1} 

Let the vertex set of $K_{181}$ be $Z_{181}$. The decompositions consist of

$(98,69,51,61,3,133,64,4,82,132,\adfsplit 23,164,91,74,49,67,56,34,104,10)_{\adfS{1}}$,

$(43,81,127,6,57,175,22,74,37,178,\adfsplit 18,7,71,167,161,154,148,16,100,104)_{\adfS{1}}$,

$(15,118,58,44,132,50,48,123,121,133,\adfsplit 0,43,23,53,114,145,82,51,77,113)_{\adfS{1}}$,

\adfLgap
$(114,15,165,106,100,85,94,67,93,55,\adfsplit 87,30,115,131,10,54,112,18,45,107)_{\adfS{2}}$,

$(109,89,116,102,3,57,169,144,173,114,\adfsplit 62,136,119,25,78,145,59,32,150,67)_{\adfS{2}}$,

$(77,59,53,139,54,111,70,5,80,134,\adfsplit 1,30,8,47,155,112,79,101,102,167)_{\adfS{2}}$,

\adfLgap
$(146,52,106,104,150,117,45,12,162,30,\adfsplit 4,1,177,87,9,57,55,173,38,56)_{\adfS{3}}$,

$(6,110,72,118,12,127,161,116,69,7,\adfsplit 97,129,170,60,101,51,80,45,35,55)_{\adfS{3}}$,

$(149,97,102,65,73,82,60,121,164,63,\adfsplit 1,46,6,64,137,77,0,99,113,85)_{\adfS{3}}$,

\adfLgap
$(97,25,150,14,33,122,74,87,62,20,\adfsplit 27,128,24,171,54,148,136,105,110,103)_{\adfS{4}}$,

$(74,86,40,164,48,80,134,30,58,119,\adfsplit 176,141,114,56,98,67,101,1,91,85)_{\adfS{4}}$,

$(35,26,172,66,170,16,74,161,6,5,\adfsplit 3,55,4,154,63,177,137,109,157,88)_{\adfS{4}}$,

\adfLgap
$(151,157,161,162,176,24,106,8,171,42,\adfsplit 38,11,81,126,175,120,121,9,25,77)_{\adfS{5}}$,

$(31,82,139,125,72,179,156,130,22,137,\adfsplit 5,95,70,52,98,167,9,44,16,166)_{\adfS{5}}$,

$(36,107,30,58,11,164,19,65,45,68,\adfsplit 4,5,80,37,100,176,166,138,160,22)_{\adfS{5}}$

\noindent and

$(46,111,108,22,161,90,102,157,114,23,\adfsplit 14,93,179,133,170,79,129,178,47,20)_{\adfS{6}}$,

$(108,77,178,112,11,164,110,140,18,145,\adfsplit 61,160,12,70,51,85,134,176,74,124)_{\adfS{6}}$,

$(168,36,42,130,93,16,100,175,7,58,\adfsplit 8,29,14,77,51,165,90,138,160,72)_{\adfS{6}}$

\noindent under the action of the mapping $x \mapsto x + 1$ (mod 181).

\adfVfy{181, \{\{181,181,1\}\}, -1, -1, -1} 

Let the vertex set of $K_{220}$ be $Z_{219} \cup \{\infty\}$. The decompositions consist of

$(\infty,69,185,32,110,27,167,122,196,123,\adfsplit 7,12,159,5,34,86,67,207,1,105)_{\adfS{1}}$,

$(58,59,15,207,17,35,40,49,14,186,\adfsplit 61,184,115,135,122,6,159,199,127,26)_{\adfS{1}}$,

$(40,34,30,91,57,51,53,100,36,96,\adfsplit 149,195,179,180,6,205,58,123,124,115)_{\adfS{1}}$,

$(40,205,141,181,4,124,74,82,59,106,\adfsplit 103,201,93,136,97,150,109,197,69,13)_{\adfS{1}}$,

$(193,85,101,137,52,50,30,179,155,81,\adfsplit 26,203,131,45,127,9,17,32,0,178)_{\adfS{1}}$,

$(37,126,40,184,163,148,23,53,49,123,\adfsplit 5,67,194,96,72,13,188,121,61,159)_{\adfS{1}}$,

$(68,110,190,80,101,133,20,62,66,37,\adfsplit 48,10,116,201,140,0,14,158,89,58)_{\adfS{1}}$,

$(204,146,33,119,154,52,129,75,137,115,\adfsplit 23,205,125,102,177,216,176,79,26,14)_{\adfS{1}}$,

$(49,171,109,156,53,60,123,21,215,36,\adfsplit 9,185,176,140,76,160,56,150,134,146)_{\adfS{1}}$,

$(173,194,189,166,65,76,10,116,197,37,\adfsplit 6,63,159,18,186,171,128,109,99,154)_{\adfS{1}}$,

$(33,8,193,75,164,28,201,186,96,93,\adfsplit 24,152,202,44,144,159,211,188,116,163)_{\adfS{1}}$,

\adfLgap
$(\infty,183,191,105,193,20,99,107,181,120,\adfsplit 165,153,128,82,15,204,73,9,45,60)_{\adfS{2}}$,

$(135,24,177,92,204,25,39,138,68,86,\adfsplit 188,139,122,170,190,9,155,110,19,30)_{\adfS{2}}$,

$(76,64,192,109,169,65,23,37,194,106,\adfsplit 146,70,201,97,55,139,199,43,151,203)_{\adfS{2}}$,

$(197,78,171,175,85,176,74,205,36,39,\adfsplit 193,7,174,111,122,61,102,151,207,13)_{\adfS{2}}$,

$(13,169,8,60,10,139,186,124,218,149,\adfsplit 56,215,136,20,71,100,201,157,154,180)_{\adfS{2}}$,

$(158,189,72,121,173,43,150,10,35,56,\adfsplit 208,98,48,137,163,65,153,193,107,54)_{\adfS{2}}$,

$(199,20,8,28,110,196,11,164,133,9,\adfsplit 124,148,203,190,88,119,39,3,45,187)_{\adfS{2}}$,

$(90,171,7,14,50,177,107,28,91,26,\adfsplit 44,82,138,172,17,9,197,48,131,5)_{\adfS{2}}$,

$(108,45,102,72,20,167,95,194,85,98,\adfsplit 35,112,114,59,146,155,140,111,78,213)_{\adfS{2}}$,

$(68,158,104,115,160,147,42,181,203,6,\adfsplit 208,168,49,169,96,48,47,185,15,75)_{\adfS{2}}$,

$(6,13,132,21,63,66,138,67,154,115,\adfsplit 185,161,43,80,68,102,148,179,103,72)_{\adfS{2}}$,

\adfLgap
$(\infty,125,71,123,131,43,61,36,5,148,\adfsplit 167,140,35,68,178,171,117,173,113,105)_{\adfS{3}}$,

$(22,24,194,178,137,204,51,73,1,43,\adfsplit 145,74,160,52,77,32,60,66,188,101)_{\adfS{3}}$,

$(40,23,185,163,132,188,202,60,52,119,\adfsplit 138,152,37,50,21,142,109,2,164,182)_{\adfS{3}}$,

$(89,198,37,78,62,26,187,200,195,138,\adfsplit 176,9,113,116,119,6,174,178,160,8)_{\adfS{3}}$,

$(1,87,209,184,169,82,78,133,33,199,\adfsplit 130,121,53,137,18,44,69,187,208,5)_{\adfS{3}}$,

$(153,146,189,1,150,110,183,21,169,151,\adfsplit 190,130,46,24,127,121,174,20,80,0)_{\adfS{3}}$,

$(34,210,19,95,66,39,214,28,61,206,\adfsplit 60,216,16,159,33,211,141,78,12,183)_{\adfS{3}}$,

$(204,179,81,50,144,109,17,36,120,65,\adfsplit 154,37,118,47,43,137,172,131,215,87)_{\adfS{3}}$,

$(5,160,138,11,32,174,25,130,133,151,\adfsplit 34,141,37,115,29,202,178,26,23,79)_{\adfS{3}}$,

$(189,9,47,155,63,72,195,25,125,108,\adfsplit 136,107,115,70,113,41,137,11,218,159)_{\adfS{3}}$,

$(39,184,21,3,60,138,115,150,204,90,\adfsplit 171,58,41,99,173,185,78,154,68,81)_{\adfS{3}}$,

\adfLgap
$(\infty,16,141,118,209,8,206,197,202,211,\adfsplit 34,105,160,69,199,100,117,46,84,72)_{\adfS{4}}$,

$(14,21,159,179,40,37,99,164,112,103,\adfsplit 106,140,197,151,44,107,180,51,1,163)_{\adfS{4}}$,

$(171,9,84,95,2,137,190,153,101,131,\adfsplit 178,14,0,55,189,108,72,38,141,192)_{\adfS{4}}$,

$(68,33,110,14,12,141,38,186,147,1,\adfsplit 78,218,17,187,121,103,197,29,27,130)_{\adfS{4}}$,

$(50,211,23,216,121,75,60,152,186,118,\adfsplit 44,104,76,101,100,9,170,105,82,57)_{\adfS{4}}$,

$(13,154,142,90,204,100,63,18,188,110,\adfsplit 70,7,153,103,176,99,216,53,8,130)_{\adfS{4}}$,

$(137,181,106,43,175,189,39,50,146,130,\adfsplit 199,207,78,110,38,157,170,138,92,125)_{\adfS{4}}$,

$(172,89,116,21,55,76,2,194,180,107,\adfsplit 96,174,66,58,166,176,168,97,150,171)_{\adfS{4}}$,

$(91,176,170,23,77,181,127,98,69,20,\adfsplit 44,55,45,66,25,182,19,39,34,46)_{\adfS{4}}$,

$(17,197,4,192,23,43,68,194,135,84,\adfsplit 28,179,25,163,12,201,70,216,24,15)_{\adfS{4}}$,

$(78,114,201,70,79,107,182,118,71,22,\adfsplit 188,166,123,67,159,65,89,99,176,112)_{\adfS{4}}$,

\adfLgap
$(\infty,2,129,162,83,82,40,84,146,109,\adfsplit 216,211,184,77,150,200,159,17,113,164)_{\adfS{5}}$,

$(114,45,100,84,47,115,17,1,39,193,\adfsplit 87,104,153,31,67,207,205,199,172,212)_{\adfS{5}}$,

$(212,147,54,25,125,136,105,120,19,143,\adfsplit 201,77,82,110,180,84,113,112,121,199)_{\adfS{5}}$,

$(186,157,81,187,140,21,44,40,162,116,\adfsplit 193,103,182,129,151,108,79,143,210,2)_{\adfS{5}}$,

$(109,50,90,185,62,79,29,154,211,128,\adfsplit 64,87,2,119,153,121,134,195,30,124)_{\adfS{5}}$,

$(29,98,70,83,57,63,73,197,16,128,\adfsplit 27,30,116,174,91,182,152,23,4,58)_{\adfS{5}}$,

$(69,7,181,170,92,74,189,206,123,173,\adfsplit 158,216,204,171,153,149,15,114,101,148)_{\adfS{5}}$,

$(14,12,52,70,134,114,51,68,46,78,\adfsplit 155,151,200,32,55,89,162,215,183,88)_{\adfS{5}}$,

$(175,2,212,192,110,213,105,112,83,168,\adfsplit 69,30,118,73,177,18,126,15,37,48)_{\adfS{5}}$,

$(177,206,185,10,218,176,71,174,94,95,\adfsplit 151,21,126,211,104,111,169,67,209,166)_{\adfS{5}}$,

$(121,72,212,0,137,199,217,58,106,7,\adfsplit 19,172,181,52,129,108,104,73,156,92)_{\adfS{5}}$

\noindent and

$(\infty,98,56,51,152,203,25,186,102,201,\adfsplit 154,22,189,143,61,141,155,13,1,147)_{\adfS{6}}$,

$(143,182,7,216,84,140,112,207,164,201,\adfsplit 53,81,205,122,115,82,200,51,1,173)_{\adfS{6}}$,

$(17,212,124,137,55,167,8,25,141,166,\adfsplit 11,164,108,139,91,105,37,46,112,38)_{\adfS{6}}$,

$(49,89,207,36,47,174,27,35,6,26,\adfsplit 67,211,46,76,19,215,122,87,199,132)_{\adfS{6}}$,

$(28,74,96,212,110,78,89,62,140,43,\adfsplit 156,95,64,166,106,162,190,67,198,49)_{\adfS{6}}$,

$(209,39,92,74,47,150,84,44,130,190,\adfsplit 79,73,137,3,180,140,204,23,14,29)_{\adfS{6}}$,

$(108,204,20,62,75,136,91,0,124,112,\adfsplit 175,60,96,149,187,202,170,5,101,177)_{\adfS{6}}$,

$(140,57,84,205,49,196,21,207,186,52,\adfsplit 69,75,81,198,164,178,62,79,115,152)_{\adfS{6}}$,

$(170,32,194,201,146,147,46,142,91,99,\adfsplit 72,206,116,202,172,95,84,112,18,75)_{\adfS{6}}$,

$(163,151,85,18,17,53,27,150,36,111,\adfsplit 204,145,174,24,205,141,83,176,164,16)_{\adfS{6}}$,

$(31,99,180,137,143,3,214,145,118,211,\adfsplit 111,5,120,57,37,107,142,218,45,200)_{\adfS{6}}$

\noindent under the action of the mapping $\infty \mapsto \infty$, $x \mapsto x + 3$ (mod 219).\eproof

\adfVfy{220, \{\{219,73,3\},\{1,1,1\}\}, 219, -1, -1} 

{\lemma \label{lem:snark20 multipartite}
There exist decompositions of the complete multipartite graphs
$K_{5,5,5,5}$, $K_{10,10,10,10}$, $K_{6,6,6,6,6}$, $K_{60,60,75}$, $K_{60,60,60,75}$,
$K_{15,15,15,21}$, $K_{24,24,24,24,39}$ and $K_{24,24,24,24,24,24,60}$
into each of the six snarks on $20$ vertices.}\\


\noindent\textbf{Proof.}
Let the vertex set of $K_{5,5,5,5}$ be $Z_{20}$ partitioned according to residue classes modulo 4.
The decomposition consists of

$(14,3,17,9,15,13,18,0,8,2,1,5,11,6,16,19,12,10,4,7)_{\adfS{1}}$,

\adfLgap
$(3,4,16,5,13,7,2,8,0,6,17,11,18,12,1,19,14,10,15,9)_{\adfS{2}}$,

\adfLgap
$(2,15,1,9,3,6,0,8,18,4,7,17,12,5,10,13,19,11,14,16)_{\adfS{3}}$,

\adfLgap
$(12,1,14,6,15,8,3,11,13,4,7,10,18,2,5,0,19,9,17,16)_{\adfS{4}}$,

\adfLgap
$(12,9,2,18,0,3,1,8,13,4,17,15,16,14,5,6,19,7,11,10)_{\adfS{5}}$

\noindent and

$(14,7,1,13,8,6,17,2,11,16,4,3,10,15,18,5,12,0,19,9)_{\adfS{6}}$

\noindent under the action of the mapping $x \mapsto x + 4$ (mod 20).

\adfVfy{20, \{\{20,5,4\}\}, -1, \{\{5,\{0,1,2,3\}\}\}, -1} 

Let the vertex set of $K_{10,10,10,10}$ be $Z_{40}$ partitioned according to residue classes modulo 4.
The decomposition consists of

$(16,35,30,32,26,37,14,7,11,17,18,15,29,20,34,33,3,36,38,27)_{\adfS{1}}$,

\adfLgap
$(13,0,11,33,7,8,30,36,17,10,5,26,1,3,25,15,16,12,20,22)_{\adfS{2}}$,

\adfLgap
$(12,35,5,26,25,4,17,24,11,3,33,34,6,27,32,7,14,9,10,16)_{\adfS{3}}$,

\adfLgap
$(4,18,14,23,36,33,7,29,27,3,26,28,8,11,12,31,5,2,1,24)_{\adfS{4}}$,

\adfLgap
$(29,11,25,23,30,8,6,4,21,26,22,31,28,3,14,7,27,18,17,12)_{\adfS{5}}$

\noindent and

$(20,6,13,30,28,29,35,1,11,14,16,3,8,34,36,5,7,39,2,21)_{\adfS{6}}$

\noindent under the action of the mapping $x \mapsto x + 2$ (mod 40).

\adfVfy{40, \{\{40,20,2\}\}, -1, \{\{10,\{0,1,2,3\}\}\}, -1} 

Let the vertex set of $K_{6,6,6,6,6}$ be $Z_{30}$ partitioned according to residue classes modulo 5.
The decomposition consists of

$(7,26,22,28,6,15,21,3,24,2,25,5,11,8,17,4,23,9,27,1)_{\adfS{1}}$,

$(5,1,20,4,9,15,28,13,25,26,24,19,12,29,0,27,8,16,3,22)_{\adfS{1}}$,

\adfLgap
$(23,20,25,4,14,12,28,2,16,3,22,18,27,0,21,10,9,19,1,8)_{\adfS{2}}$,

$(28,1,15,12,24,5,11,18,4,7,17,26,29,21,13,25,20,14,22,16)_{\adfS{2}}$,

\adfLgap
$(24,5,21,12,3,14,7,17,1,18,11,20,16,10,4,8,25,23,15,2)_{\adfS{3}}$,

$(25,3,2,28,23,11,9,15,14,4,12,21,24,8,22,20,26,27,1,19)_{\adfS{3}}$,

\adfLgap
$(21,18,0,19,17,5,13,6,4,16,24,10,23,8,20,7,12,27,26,9)_{\adfS{4}}$,

$(25,28,4,0,2,8,19,24,1,12,14,6,26,23,10,17,11,22,20,13)_{\adfS{4}}$,

\adfLgap
$(4,23,1,16,10,9,6,22,3,27,26,8,20,15,21,28,25,18,7,19)_{\adfS{5}}$,

$(0,22,16,7,15,2,4,1,24,13,20,17,14,9,25,8,11,23,12,19)_{\adfS{5}}$

\noindent and

$(11,2,19,17,26,27,18,14,15,23,13,12,1,24,16,0,10,3,4,6)_{\adfS{6}}$,

$(25,18,20,28,3,15,12,1,4,14,23,17,26,11,24,5,2,22,19,10)_{\adfS{6}}$

\noindent under the action of the mapping $x \mapsto x + 5$ (mod 30).

\adfVfy{30, \{\{30,6,5\}\}, -1, \{\{6,\{0,1,2,3,4\}\}\}, -1} 

Let the vertex set of $K_{60,60,75}$ be $\{0, 1, \dots, 194\}$ partitioned into
$\{3j:     j = 0,1,\dots,59\}$,
$\{3j + 1: j = 0,1,\dots,59\}$ and
$\{3j + 2: j = 0,1,\dots,59\} \cup \{180, 181, \dots, 194\}$.
The decomposition consists of

$(81,8,27,0,186,151,86,107,169,89,\adfsplit 124,37,9,173,6,105,13,68,166,115)_{\adfS{1}}$,

$(9,190,169,66,113,78,76,53,127,65,\adfsplit 114,21,84,186,163,121,143,19,131,135)_{\adfS{1}}$,

$(127,27,85,121,156,17,18,15,106,183,\adfsplit 55,124,47,153,149,71,138,130,174,48)_{\adfS{1}}$,

$(55,62,1,115,95,124,193,92,25,108,\adfsplit 91,157,104,144,137,189,9,14,174,165)_{\adfS{1}}$,

$(42,65,39,54,118,152,139,130,182,174,\adfsplit 107,47,102,189,19,61,80,66,11,56)_{\adfS{1}}$,

$(39,131,172,157,106,125,42,177,95,64,\adfsplit 107,158,19,92,100,36,161,163,86,167)_{\adfS{1}}$,

$(146,0,149,47,3,180,15,87,148,30,\adfsplit 31,187,75,88,66,63,49,90,142,43)_{\adfS{1}}$,

\adfLgap
$(61,153,134,4,191,21,13,18,159,82,\adfsplit 148,189,81,121,39,26,92,129,143,88)_{\adfS{2}}$,

$(154,146,186,175,93,152,94,137,159,73,\adfsplit 28,187,99,12,121,47,26,144,119,124)_{\adfS{2}}$,

$(147,10,17,69,137,13,39,34,98,117,\adfsplit 72,4,183,140,162,178,40,168,1,80)_{\adfS{2}}$,

$(144,191,53,91,10,92,135,26,29,18,\adfsplit 172,96,192,55,165,142,128,6,15,28)_{\adfS{2}}$,

$(113,159,139,46,3,84,107,35,61,9,\adfsplit 161,42,13,191,135,104,133,111,168,83)_{\adfS{2}}$,

$(14,177,106,188,129,88,180,136,162,23,\adfsplit 82,50,121,113,64,140,125,102,42,91)_{\adfS{2}}$,

$(14,30,7,176,108,31,164,55,52,38,\adfsplit 146,15,11,105,26,12,136,167,178,45)_{\adfS{2}}$,

\adfLgap
$(71,93,29,25,11,154,138,56,21,108,\adfsplit 168,55,15,64,18,151,67,183,137,177)_{\adfS{3}}$,

$(48,40,191,92,119,114,173,63,136,10,\adfsplit 154,104,126,167,28,137,60,124,117,192)_{\adfS{3}}$,

$(18,178,89,26,27,61,41,104,79,148,\adfsplit 68,120,149,168,77,154,3,139,39,180)_{\adfS{3}}$,

$(32,157,96,108,117,47,94,192,139,35,\adfsplit 115,78,86,92,64,9,125,111,191,73)_{\adfS{3}}$,

$(187,156,67,55,173,94,165,128,90,66,\adfsplit 83,79,159,61,5,154,93,146,50,31)_{\adfS{3}}$,

$(88,81,13,129,26,61,78,128,28,187,\adfsplit 36,94,111,139,56,16,174,29,32,58)_{\adfS{3}}$,

$(0,17,25,190,113,162,106,147,49,105,\adfsplit 128,55,81,160,66,134,92,188,63,4)_{\adfS{3}}$,

\adfLgap
$(38,84,64,166,2,55,138,131,120,93,\adfsplit 183,94,14,16,92,163,150,156,188,61)_{\adfS{4}}$,

$(158,150,160,44,45,193,141,22,143,129,\adfsplit 66,92,149,14,54,82,73,115,50,105)_{\adfS{4}}$,

$(189,7,90,71,178,20,36,81,47,65,\adfsplit 141,154,135,175,95,91,10,131,92,78)_{\adfS{4}}$,

$(95,105,43,8,49,63,67,148,168,32,\adfsplit 70,45,39,65,7,20,57,22,171,62)_{\adfS{4}}$,

$(12,91,169,110,160,104,49,45,182,71,\adfsplit 0,109,148,116,114,181,22,21,171,140)_{\adfS{4}}$,

$(50,144,118,77,124,104,58,14,129,185,\adfsplit 69,4,13,101,24,189,163,60,117,184)_{\adfS{4}}$,

$(0,22,49,8,120,7,152,57,164,135,\adfsplit 9,127,183,101,4,84,30,76,46,53)_{\adfS{4}}$,

\adfLgap
$(36,52,168,43,155,103,77,98,60,65,\adfsplit 48,70,93,20,31,32,51,85,21,41)_{\adfS{5}}$,

$(133,164,61,167,29,120,119,39,79,105,\adfsplit 122,124,189,142,62,49,55,69,3,95)_{\adfS{5}}$,

$(164,18,22,171,27,155,160,124,183,190,\adfsplit 136,54,121,113,7,159,66,108,122,4)_{\adfS{5}}$,

$(185,118,110,30,100,167,36,69,65,157,\adfsplit 120,37,62,49,126,25,59,138,95,163)_{\adfS{5}}$,

$(95,63,22,141,36,70,126,61,158,181,\adfsplit 117,5,135,175,104,145,177,180,86,28)_{\adfS{5}}$,

$(81,161,132,133,97,123,188,177,183,38,\adfsplit 176,112,186,15,125,42,117,4,79,95)_{\adfS{5}}$,

$(0,1,87,88,148,95,112,162,194,131,\adfsplit 46,119,21,125,31,86,60,64,12,35)_{\adfS{5}}$

\noindent and

$(67,102,152,161,115,44,75,187,114,85,\adfsplit 76,51,96,137,21,173,4,157,37,33)_{\adfS{6}}$,

$(187,174,112,90,6,184,22,14,70,185,\adfsplit 4,171,13,186,166,138,170,150,11,118)_{\adfS{6}}$,

$(65,130,116,84,73,99,115,30,143,42,\adfsplit 148,44,185,33,20,105,178,55,14,132)_{\adfS{6}}$,

$(57,122,18,7,47,117,83,160,129,181,\adfsplit 65,81,99,142,66,98,41,140,91,84)_{\adfS{6}}$,

$(91,105,46,11,77,124,173,149,103,143,\adfsplit 83,132,3,12,148,193,57,85,157,140)_{\adfS{6}}$,

$(41,96,104,28,109,181,115,73,156,126,\adfsplit 97,38,13,189,16,18,69,174,125,106)_{\adfS{6}}$,

$(0,5,4,119,32,148,86,84,172,95,\adfsplit 39,124,163,161,136,41,126,93,51,89)_{\adfS{6}}$

\noindent under the action of the mapping $x \mapsto x + 3$ (mod 180) for $x < 180$,
$x \mapsto 180 + (x + 1 \mathrm{~(mod~15)})$ for $x \ge 180$.

\adfVfy{195, \{\{180,60,3\},\{15,15,1\}\}, -1, \{\{60,\{0,1,2\}\},\{15,\{2\}\}\}, -1} 

Let the vertex set of $K_{60,60,60,75}$ be $\{0, 1, \dots, 254\}$ partitioned into
$\{3j + i: j = 0, 1, \dots, 59\}$, $i = 0, 1, 2$, and $\{180, 181, \dots, 254\}$.
The decomposition consists of

$(251,25,138,197,92,73,171,23,0,244,\adfsplit 157,58,44,232,122,29,21,164,109,127)_{\adfS{1}}$,

$(46,153,76,236,227,32,42,135,78,97,\adfsplit 161,98,247,175,213,230,30,180,4,149)_{\adfS{1}}$,

$(167,21,110,104,230,171,157,76,66,31,\adfsplit 36,249,253,176,168,233,160,214,165,172)_{\adfS{1}}$,

$(247,29,238,27,177,225,131,114,186,93,\adfsplit 146,151,166,179,15,21,107,84,222,85)_{\adfS{1}}$,

$(187,175,5,42,82,199,1,96,209,63,\adfsplit 76,238,101,240,79,74,69,143,190,221)_{\adfS{1}}$,

$(45,227,178,114,101,16,128,68,186,137,\adfsplit 184,229,121,3,140,97,53,183,125,171)_{\adfS{1}}$,

$(118,206,31,158,147,204,64,119,214,81,\adfsplit 192,201,36,61,78,140,251,124,233,220)_{\adfS{1}}$,

$(3,145,5,36,104,220,83,70,46,117,\adfsplit 215,203,0,119,162,115,28,177,248,136)_{\adfS{1}}$,

$(47,215,146,6,211,28,191,66,179,48,\adfsplit 83,26,189,50,242,0,76,238,104,72)_{\adfS{1}}$,

\adfLgap
$(100,96,187,53,251,145,87,174,59,46,\adfsplit 55,229,159,240,28,228,119,143,107,201)_{\adfS{2}}$,

$(142,174,171,73,181,197,175,113,125,207,\adfsplit 161,214,30,8,168,248,61,62,58,99)_{\adfS{2}}$,

$(136,227,75,143,223,237,70,138,108,169,\adfsplit 71,97,246,67,59,144,115,134,135,200)_{\adfS{2}}$,

$(68,21,233,50,209,226,126,97,224,51,\adfsplit 103,119,185,15,130,137,32,118,5,54)_{\adfS{2}}$,

$(116,132,192,201,36,13,80,82,194,68,\adfsplit 40,6,41,35,213,27,89,154,18,233)_{\adfS{2}}$,

$(234,2,59,10,125,150,82,132,184,79,\adfsplit 163,96,235,113,233,167,45,135,24,110)_{\adfS{2}}$,

$(15,112,185,33,219,158,159,31,89,114,\adfsplit 211,163,217,169,135,10,139,134,84,250)_{\adfS{2}}$,

$(241,84,120,158,127,157,162,228,247,89,\adfsplit 181,68,171,252,108,229,46,146,94,220)_{\adfS{2}}$,

$(95,186,112,52,202,248,85,122,108,204,\adfsplit 19,198,13,12,10,86,230,250,190,9)_{\adfS{2}}$,

\adfLgap
$(18,220,162,134,39,219,133,31,253,127,\adfsplit 6,149,189,229,79,150,239,14,234,63)_{\adfS{3}}$,

$(216,16,114,31,185,152,93,235,27,248,\adfsplit 174,65,155,59,42,208,28,251,173,84)_{\adfS{3}}$,

$(155,156,56,253,190,100,223,106,165,148,\adfsplit 204,114,136,107,40,47,9,69,37,176)_{\adfS{3}}$,

$(234,173,182,14,164,241,99,172,45,166,\adfsplit 114,242,184,29,10,6,54,0,91,232)_{\adfS{3}}$,

$(124,27,250,23,3,85,149,16,202,79,\adfsplit 186,53,208,74,178,159,95,84,246,100)_{\adfS{3}}$,

$(230,129,127,64,71,111,97,157,202,41,\adfsplit 224,54,25,61,9,59,227,246,117,5)_{\adfS{3}}$,

$(187,112,164,45,203,48,8,207,178,94,\adfsplit 248,41,86,68,235,40,121,242,188,153)_{\adfS{3}}$,

$(198,70,12,63,232,147,92,177,107,74,\adfsplit 43,221,55,68,180,20,118,194,250,141)_{\adfS{3}}$,

$(10,149,144,78,66,5,111,165,40,216,\adfsplit 206,2,109,64,62,235,221,32,23,52)_{\adfS{3}}$,

\adfLgap
$(165,253,50,99,58,9,17,57,248,226,\adfsplit 207,70,44,122,108,202,8,200,100,162)_{\adfS{4}}$,

$(116,144,7,154,189,84,211,234,69,177,\adfsplit 118,161,11,176,124,210,101,40,57,205)_{\adfS{4}}$,

$(75,198,134,99,104,210,10,156,79,154,\adfsplit 113,221,220,188,96,143,204,97,60,20)_{\adfS{4}}$,

$(170,64,243,244,108,145,189,66,164,19,\adfsplit 239,47,17,22,192,141,29,208,46,24)_{\adfS{4}}$,

$(32,232,31,61,211,83,118,129,65,132,\adfsplit 54,248,247,223,124,3,131,99,191,157)_{\adfS{4}}$,

$(238,88,97,107,0,197,44,113,244,141,\adfsplit 40,232,47,32,251,94,55,121,229,98)_{\adfS{4}}$,

$(40,186,150,63,17,204,53,87,34,215,\adfsplit 49,137,143,183,147,190,130,70,18,119)_{\adfS{4}}$,

$(211,111,52,222,121,190,162,11,219,140,\adfsplit 104,171,102,9,149,220,0,113,28,187)_{\adfS{4}}$,

$(28,140,44,36,231,55,167,103,48,76,\adfsplit 180,5,77,20,13,246,155,81,9,212)_{\adfS{4}}$,

\adfLgap
$(166,182,20,60,148,193,81,155,234,181,\adfsplit 52,200,118,30,252,132,97,183,223,165)_{\adfS{5}}$,

$(146,55,18,177,67,237,78,171,80,196,\adfsplit 136,89,46,219,94,53,63,122,235,166)_{\adfS{5}}$,

$(161,163,247,30,26,168,151,195,109,77,\adfsplit 82,42,201,244,73,131,227,81,208,22)_{\adfS{5}}$,

$(177,16,0,166,201,18,106,83,123,146,\adfsplit 44,63,19,237,92,213,67,143,78,215)_{\adfS{5}}$,

$(94,164,219,21,113,117,152,247,4,55,\adfsplit 178,123,217,2,235,85,77,67,72,253)_{\adfS{5}}$,

$(233,13,14,113,22,17,91,147,201,196,\adfsplit 76,236,83,65,180,120,112,60,204,110)_{\adfS{5}}$,

$(212,85,216,75,92,224,160,155,196,140,\adfsplit 167,238,150,229,113,52,63,2,170,190)_{\adfS{5}}$,

$(96,43,137,40,24,1,209,95,203,33,\adfsplit 192,10,235,76,243,4,178,101,5,230)_{\adfS{5}}$,

$(13,167,126,101,155,229,21,94,254,1,\adfsplit 112,48,221,165,151,9,119,14,173,210)_{\adfS{5}}$

\noindent and

$(174,137,58,181,0,1,182,117,46,150,\adfsplit 50,22,207,61,245,152,34,38,126,221)_{\adfS{6}}$,

$(174,32,31,151,156,193,145,171,248,99,\adfsplit 83,246,194,195,139,177,211,60,190,104)_{\adfS{6}}$,

$(79,114,4,191,12,187,14,183,37,24,\adfsplit 80,121,46,222,32,51,72,159,175,239)_{\adfS{6}}$,

$(61,143,81,218,208,62,177,192,103,160,\adfsplit 138,47,174,184,135,149,178,97,145,216)_{\adfS{6}}$,

$(225,169,187,175,129,206,96,59,215,43,\adfsplit 38,202,229,94,253,84,6,179,178,223)_{\adfS{6}}$,

$(201,29,90,134,124,141,163,82,117,127,\adfsplit 253,2,229,176,250,148,60,27,135,183)_{\adfS{6}}$,

$(218,92,210,30,141,37,143,160,225,91,\adfsplit 196,96,62,166,116,209,163,229,130,120)_{\adfS{6}}$,

$(250,87,34,54,117,148,164,167,58,202,\adfsplit 89,37,5,231,62,94,184,7,180,147)_{\adfS{6}}$,

$(6,55,249,67,177,131,124,45,212,9,\adfsplit 165,1,229,43,129,22,41,146,232,148)_{\adfS{6}}$

\noindent under the action of the mapping $x \mapsto x + 2$ (mod 180) for $x < 180$,
$x \mapsto 180 + (x - 180 + 5 \mathrm{~(mod~75)})$ for $x \ge 180$.

\adfVfy{255, \{\{180,90,2\},\{75,15,5\}\}, -1, \{\{60,\{0,1,2\}\},\{75,\{3\}\}\}, -1} 

Let the vertex set of $K_{15,15,15,21}$ be $\{0, 1, \dots, 65\}$ partitioned into
$\{3j + i: j = 0, 1, \dots, 14\}$, $i = 0, 1, 2$, and $\{45, 46, \dots, 65\}$.
The decomposition consists of

$(57,0,48,51,18,11,31,6,28,3,47,50,55,36,7,40,1,64,24,14)_{\adfS{1}}$,

$(56,37,11,52,0,61,21,33,60,40,29,59,35,51,7,19,34,39,41,49)_{\adfS{1}}$,

$(50,12,20,23,2,4,36,21,60,41,51,56,39,22,17,38,8,6,37,7)_{\adfS{1}}$,

$(21,46,9,3,13,29,60,52,5,62,4,23,45,38,0,12,10,17,48,47)_{\adfS{1}}$,

$(0,62,44,2,28,33,58,51,11,19,32,15,21,5,40,36,14,3,54,63)_{\adfS{1}}$,

$(37,47,14,19,35,57,44,40,12,52,38,18,1,20,49,25,53,17,42,13)_{\adfS{1}}$,

\adfLgap
$(14,16,47,45,0,38,63,44,7,9,20,3,4,15,58,32,28,31,56,24)_{\adfS{2}}$,

$(54,24,26,8,28,64,15,36,7,48,1,51,5,2,33,22,17,58,60,4)_{\adfS{2}}$,

$(61,43,8,51,27,3,2,4,46,18,57,5,15,53,9,26,17,28,13,6)_{\adfS{2}}$,

$(37,47,60,20,56,39,5,42,62,40,12,38,19,48,27,45,17,6,8,59)_{\adfS{2}}$,

$(17,45,24,39,30,50,10,62,15,44,25,59,28,16,48,33,32,49,64,36)_{\adfS{2}}$,

$(0,4,1,65,47,6,11,40,63,33,44,36,10,31,12,26,64,60,41,9)_{\adfS{2}}$,

\adfLgap
$(21,35,27,55,57,22,54,63,30,37,8,6,43,61,34,5,23,39,15,62)_{\adfS{3}}$,

$(47,20,51,21,0,17,12,40,60,41,61,6,43,33,53,38,36,45,25,44)_{\adfS{3}}$,

$(34,39,8,35,57,27,48,46,4,18,55,38,17,13,41,63,64,49,9,19)_{\adfS{3}}$,

$(35,55,31,39,60,9,45,26,7,23,20,57,16,58,38,27,51,12,4,29)_{\adfS{3}}$,

$(54,39,57,38,29,27,5,26,31,16,51,9,22,19,41,59,25,63,10,33)_{\adfS{3}}$,

$(14,33,1,13,10,6,22,32,51,63,45,25,36,37,18,59,2,35,42,60)_{\adfS{3}}$,

\adfLgap
$(28,47,42,6,44,19,45,30,40,10,36,23,5,17,16,64,62,25,24,32)_{\adfS{4}}$,

$(57,5,8,19,20,45,21,52,1,43,10,27,61,49,33,32,60,24,53,40)_{\adfS{4}}$,

$(6,2,43,22,64,17,55,34,44,59,27,11,9,61,35,51,56,37,30,26)_{\adfS{4}}$,

$(34,8,36,64,0,23,18,4,58,63,45,2,11,12,62,40,3,51,29,22)_{\adfS{4}}$,

$(32,28,46,62,4,64,27,11,36,29,6,61,60,16,53,38,30,2,49,13)_{\adfS{4}}$,

$(43,23,41,9,35,12,65,15,1,18,34,61,58,48,44,39,63,28,31,30)_{\adfS{4}}$,

\adfLgap
$(51,17,7,40,20,15,41,64,9,61,1,33,60,62,21,19,11,13,52,24)_{\adfS{5}}$,

$(61,17,28,18,58,8,24,4,20,56,36,63,41,47,43,64,42,25,7,45)_{\adfS{5}}$,

$(31,9,14,63,46,40,36,5,30,28,18,43,39,41,59,13,33,53,11,7)_{\adfS{5}}$,

$(62,12,35,37,43,50,15,28,36,51,19,48,29,0,59,27,21,7,32,60)_{\adfS{5}}$,

$(46,35,52,43,44,54,1,0,61,33,7,56,19,15,14,58,34,49,17,42)_{\adfS{5}}$,

$(2,33,5,36,16,48,0,37,44,64,25,57,29,26,12,55,31,19,13,51)_{\adfS{5}}$

\noindent and

$(6,2,13,14,48,21,35,59,16,17,63,12,30,5,28,44,64,47,56,25)_{\adfS{6}}$,

$(54,33,40,23,58,3,11,32,1,14,43,9,25,49,21,57,53,27,18,16)_{\adfS{6}}$,

$(59,8,4,36,53,15,13,60,40,55,41,12,43,1,39,35,26,18,42,48)_{\adfS{6}}$,

$(36,56,8,64,21,54,14,27,1,58,37,44,26,24,43,15,33,25,55,20)_{\adfS{6}}$,

$(63,33,61,29,21,4,39,34,52,44,20,42,47,58,43,12,53,22,64,32)_{\adfS{6}}$,

$(0,58,9,20,50,12,37,27,10,65,53,26,7,39,19,33,45,55,59,25)_{\adfS{6}}$

\noindent under the action of the mapping $x \mapsto x + 5$ (mod 45) for $x < 45$,
$x \mapsto 45 + (x - 45 + 7 \mathrm{~(mod~21)})$ for $x \ge 45$.

\adfVfy{66, \{\{45,9,5\},\{21,3,7\}\}, -1, \{\{15,\{0,1,2\}\},\{21,\{3\}\}\}, -1} 

Let the vertex set of $K_{24,24,24,24,39}$ be $\{0, 1, \dots, 134\}$ partitioned into
$\{4j + i: j = 0, 1, \dots, 23\}$, $i = 0, 1, 2, 3$, and $\{96, 97, \dots, 134\}$.
The decomposition consists of

$(0,124,68,29,116,45,63,47,35,132,\adfsplit 38,93,133,44,107,67,87,56,46,30)_{\adfS{1}}$,

$(75,17,112,60,49,121,82,51,48,123,\adfsplit 10,33,18,87,36,102,125,92,29,38)_{\adfS{1}}$,

$(94,43,30,90,23,92,1,122,42,57,\adfsplit 47,124,31,126,70,54,97,87,93,129)_{\adfS{1}}$,

$(26,119,88,5,105,0,34,117,41,95,\adfsplit 49,84,109,19,113,72,44,31,127,9)_{\adfS{1}}$,

$(14,33,74,19,117,63,78,68,29,64,\adfsplit 31,39,101,4,105,134,0,102,20,13)_{\adfS{1}}$,

\adfLgap
$(128,47,26,102,6,24,59,15,12,131,\adfsplit 62,7,16,107,65,124,96,57,38,55)_{\adfS{2}}$,

$(68,30,5,89,34,133,80,112,73,18,\adfsplit 35,100,28,25,58,57,110,88,42,19)_{\adfS{2}}$,

$(93,82,74,4,68,43,130,94,30,114,\adfsplit 119,75,54,108,21,106,55,64,24,121)_{\adfS{2}}$,

$(93,23,15,62,8,69,118,85,65,101,\adfsplit 129,14,98,39,74,76,21,96,81,6)_{\adfS{2}}$,

$(75,130,100,83,1,33,91,99,5,30,\adfsplit 4,116,62,32,29,15,58,74,35,123)_{\adfS{2}}$,

\adfLgap
$(48,57,92,39,124,20,131,80,85,122,\adfsplit 121,29,34,87,44,75,94,114,69,68)_{\adfS{3}}$,

$(81,46,79,97,2,129,27,68,83,14,\adfsplit 12,39,107,6,45,7,92,57,44,54)_{\adfS{3}}$,

$(89,79,4,34,10,3,113,26,47,33,\adfsplit 104,44,87,22,84,110,128,76,86,83)_{\adfS{3}}$,

$(103,67,132,72,130,81,15,98,95,39,\adfsplit 86,20,70,118,59,107,2,45,42,96)_{\adfS{3}}$,

$(43,121,28,21,34,99,61,78,125,101,\adfsplit 3,29,113,128,44,30,119,79,47,45)_{\adfS{3}}$,

\adfLgap
$(77,7,115,97,32,69,23,60,70,6,\adfsplit 48,79,119,53,122,37,24,71,94,113)_{\adfS{4}}$,

$(110,60,71,21,32,39,61,26,5,114,\adfsplit 105,12,13,85,40,120,15,76,38,103)_{\adfS{4}}$,

$(26,25,37,42,83,92,73,112,24,107,\adfsplit 131,51,40,30,125,58,13,132,56,23)_{\adfS{4}}$,

$(126,40,71,93,117,7,45,10,16,84,\adfsplit 69,104,81,98,59,108,66,67,44,101)_{\adfS{4}}$,

$(75,26,121,116,57,42,77,98,72,44,\adfsplit 20,70,86,15,82,51,37,122,128,34)_{\adfS{4}}$,

\adfLgap
$(15,29,97,130,58,36,61,18,48,13,\adfsplit 37,7,70,127,27,106,109,16,19,53)_{\adfS{5}}$,

$(51,32,59,53,41,42,132,47,0,36,\adfsplit 48,102,27,104,1,34,112,46,52,17)_{\adfS{5}}$,

$(125,59,64,48,50,23,16,72,93,98,\adfsplit 103,7,45,0,94,118,32,105,42,71)_{\adfS{5}}$,

$(127,47,110,90,46,109,44,88,43,32,\adfsplit 24,126,38,133,5,95,128,15,120,12)_{\adfS{5}}$,

$(74,53,63,83,72,55,121,40,116,111,\adfsplit 20,108,9,66,51,46,73,35,21,100)_{\adfS{5}}$

\noindent and

$(78,116,74,43,126,89,115,23,60,77,\adfsplit 22,48,80,24,54,67,110,62,20,53)_{\adfS{6}}$,

$(9,118,31,95,107,59,108,50,84,54,\adfsplit 124,39,7,60,82,119,1,131,12,2)_{\adfS{6}}$,

$(25,46,115,64,21,39,4,55,14,91,\adfsplit 61,36,28,7,66,98,117,126,57,51)_{\adfS{6}}$,

$(23,103,13,9,125,8,99,80,35,70,\adfsplit 110,75,107,4,49,68,122,21,54,48)_{\adfS{6}}$,

$(19,96,0,69,43,101,18,40,130,31,\adfsplit 2,108,61,119,11,114,42,3,72,65)_{\adfS{6}}$

\noindent under the action of the mapping $x \mapsto x + 2$ (mod 96) for $x < 96$,
$x \mapsto 96 + (x - 96 + 13 \mathrm{~(mod~39)} )$ for $x \ge 96$.

\adfVfy{135, \{\{96,48,2\},\{39,3,13\}\}, -1, \{\{24,\{0,1,2,3\}\},\{39,\{4\}\}\}, -1} 

Let the vertex set of $K_{24,24,24,24,24,24,60}$ be $\{0, 1, \dots, 203\}$ partitioned into
$\{6j + i: j = 0, 1, \dots, 23\}$, $i = 0, 1, \dots, 5$, and $\{144, 145, \dots, 203\}$.
The decomposition consists of

$(130,126,170,37,198,66,131,68,29,169,\adfsplit 117,75,171,110,17,179,19,201,16,134)_{\adfS{1}}$,

$(62,184,85,103,22,75,23,104,109,155,\adfsplit 28,122,191,99,157,198,70,129,60,59)_{\adfS{1}}$,

$(30,50,198,132,74,33,142,47,27,12,\adfsplit 15,75,180,86,202,61,101,40,172,49)_{\adfS{1}}$,

$(134,166,119,15,168,133,128,194,51,23,\adfsplit 6,0,172,10,143,44,95,175,22,12)_{\adfS{1}}$,

\adfLgap
$(131,10,14,41,135,114,185,106,95,191,\adfsplit 102,88,109,64,123,173,22,29,122,152)_{\adfS{2}}$,

$(157,33,34,59,19,193,14,199,123,121,\adfsplit 192,143,110,128,169,54,39,131,73,164)_{\adfS{2}}$,

$(90,135,124,175,176,131,120,142,25,125,\adfsplit 72,166,83,173,24,14,33,96,80,58)_{\adfS{2}}$,

$(169,57,131,17,22,78,181,63,101,190,\adfsplit 4,170,126,193,61,162,80,38,142,6)_{\adfS{2}}$,

\adfLgap
$(168,82,154,2,159,6,43,36,89,33,\adfsplit 35,44,32,104,94,111,109,19,151,0)_{\adfS{3}}$,

$(122,16,127,163,75,180,70,164,41,111,\adfsplit 64,61,101,140,97,74,8,69,190,71)_{\adfS{3}}$,

$(92,157,63,148,179,59,186,40,137,91,\adfsplit 54,192,49,156,15,71,22,48,75,115)_{\adfS{3}}$,

$(80,201,4,124,37,95,165,115,187,103,\adfsplit 15,54,57,29,99,160,147,188,2,84)_{\adfS{3}}$,

\adfLgap
$(13,98,32,55,6,155,122,145,34,103,\adfsplit 92,35,149,156,51,4,82,81,202,104)_{\adfS{4}}$,

$(101,39,51,77,36,23,98,130,181,161,\adfsplit 156,49,125,102,119,194,85,88,31,153)_{\adfS{4}}$,

$(22,66,18,174,89,94,54,86,190,8,\adfsplit 48,183,29,187,79,118,46,117,119,153)_{\adfS{4}}$,

$(61,158,10,112,123,155,45,42,71,147,\adfsplit 17,116,6,81,90,64,97,164,157,122)_{\adfS{4}}$,

\adfLgap
$(170,5,54,76,157,9,74,195,116,8,\adfsplit 112,120,113,189,34,172,146,97,135,78)_{\adfS{5}}$,

$(105,58,158,189,50,71,178,148,83,32,\adfsplit 170,45,82,133,44,154,102,135,79,130)_{\adfS{5}}$,

$(152,88,127,46,113,50,33,161,36,3,\adfsplit 82,55,173,120,51,48,26,21,133,156)_{\adfS{5}}$,

$(145,13,2,15,146,63,41,4,77,28,\adfsplit 7,3,34,36,71,192,149,86,56,22)_{\adfS{5}}$

\noindent and

$(96,154,97,140,164,112,127,19,101,201,\adfsplit 36,27,114,145,32,85,141,136,143,151)_{\adfS{6}}$,

$(27,73,168,46,131,129,133,105,196,125,\adfsplit 153,36,102,3,134,107,113,91,62,202)_{\adfS{6}}$,

$(116,106,159,142,9,72,123,162,92,35,\adfsplit 189,109,156,95,190,68,114,31,202,11)_{\adfS{6}}$,

$(168,122,79,46,32,163,127,119,192,65,\adfsplit 3,40,8,37,20,112,131,180,190,36)_{\adfS{6}}$

\noindent under the action of the mapping $x \mapsto x + 1$ (mod 144) for $x < 144$,
$x \mapsto 144 + (x - 144 + 5 \mathrm{~(mod~60)})$ for $x \ge 144$.\eproof

\adfVfy{204, \{\{144,144,1\},\{60,12,5\}\}, -1, \{\{24,\{0,1,2,3,4,5\}\},\{60,\{6\}\}\}, -1} 

\vskip 2mm
Theorem~\ref{thm:snark20 20} follows from
Lemmas~\ref{lem:snark20 designs} and \ref{lem:snark20 multipartite}, and Proposition~\ref{prop:d=3, v=20}.


\section{The twenty non-trivial snarks on 22 vertices}
\label{sec:Loupekine}
\newcommand{\adfLa}{\mathrm{L1}}
\newcommand{\adfLb}{\mathrm{L2}}
\newcommand{\adfLc}{\mathrm{L3}}
\newcommand{\adfLd}{\mathrm{L4}}
\newcommand{\adfLe}{\mathrm{L5}}
\newcommand{\adfLf}{\mathrm{L6}}
\newcommand{\adfLg}{\mathrm{L7}}
\newcommand{\adfLh}{\mathrm{L8}}
\newcommand{\adfLi}{\mathrm{L9}}
\newcommand{\adfLj}{\mathrm{L10}}
\newcommand{\adfLk}{\mathrm{L11}}
\newcommand{\adfLl}{\mathrm{L12}}
\newcommand{\adfLm}{\mathrm{L13}}
\newcommand{\adfLn}{\mathrm{L14}}
\newcommand{\adfLo}{\mathrm{L15}}
\newcommand{\adfLp}{\mathrm{L16}}
\newcommand{\adfLq}{\mathrm{L17}}
\newcommand{\adfLr}{\mathrm{L18}}
\newcommand{\adfLs}{\mathrm{L19}}
\newcommand{\adfLt}{\mathrm{L20}}
There are twenty non-trivial snarks on 22 vertices. We represent them by ordered 22-tuples of vertices:
(1, 2, \dots, 22)${}_{\adfLa}$, (1, 2, \dots, 22)${}_{\adfLb}$, \dots, (1, 2, \dots, 22)${}_{\adfLt}$.
The edge sets are respectively

$\adfLa$:
  \{$\{1,2\}$, $\{1,3\}$, $\{1,10\}$, $\{2,4\}$, $\{2,21\}$, $\{3,6\}$, $\{3,7\}$, $\{4,5\}$, $\{4,7\}$, $\{5,6\}$, $\{5,8\}$, $\{6,17\}$, $\{7,11\}$,
   $\{8,13\}$, $\{8,15\}$, $\{9,10\}$, $\{9,12\}$, $\{9,15\}$, $\{10,14\}$, $\{11,12\}$, $\{11,16\}$, $\{12,13\}$, $\{13,14\}$, $\{14,22\}$,
   $\{15,18\}$, $\{16,19\}$, $\{16,20\}$, $\{17,18\}$, $\{17,20\}$, $\{18,19\}$, $\{19,21\}$, $\{20,22\}$, $\{21,22\}$\},

$\adfLb$:
  \{$\{1,2\}$, $\{1,3\}$, $\{1,10\}$, $\{2,4\}$, $\{2,17\}$, $\{3,6\}$, $\{3,7\}$, $\{4,5\}$, $\{4,7\}$, $\{5,6\}$, $\{5,8\}$, $\{6,21\}$, $\{7,11\}$,
   $\{8,13\}$, $\{8,15\}$, $\{9,10\}$, $\{9,12\}$, $\{9,15\}$, $\{10,14\}$, $\{11,12\}$, $\{11,16\}$, $\{12,13\}$, $\{13,14\}$, $\{14,22\}$,
   $\{15,18\}$, $\{16,19\}$, $\{16,20\}$, $\{17,18\}$, $\{17,20\}$, $\{18,19\}$, $\{19,21\}$, $\{20,22\}$, $\{21,22\}$\},

$\adfLc$:
  \{$\{1,2\}$, $\{1,3\}$, $\{1,4\}$, $\{2,5\}$, $\{2,6\}$, $\{3,7\}$, $\{3,8\}$, $\{4,9\}$, $\{4,10\}$, $\{5,7\}$, $\{5,9\}$,
   $\{6,8\}$, $\{6,10\}$, $\{7,11\}$, $\{8,12\}$, $\{9,13\}$, $\{10,14\}$, $\{11,14\}$, $\{11,15\}$, $\{12,13\}$, $\{12,16\}$,
   $\{13,17\}$, $\{14,18\}$, $\{15,19\}$, $\{15,20\}$, $\{16,19\}$, $\{16,21\}$, $\{17,20\}$, $\{17,22\}$, $\{18,21\}$, $\{18,22\}$,
   $\{19,22\}$, $\{20,21\}$\},

$\adfLd$:
  \{$\{1,2\}$, $\{1,3\}$, $\{1,4\}$, $\{2,5\}$, $\{2,6\}$, $\{3,7\}$, $\{3,8\}$, $\{4,9\}$, $\{4,10\}$, $\{5,7\}$, $\{5,9\}$,
   $\{6,8\}$, $\{6,11\}$, $\{7,10\}$, $\{8,12\}$, $\{9,13\}$, $\{10,14\}$, $\{11,14\}$, $\{11,15\}$, $\{12,13\}$, $\{12,16\}$,
   $\{13,17\}$, $\{14,18\}$, $\{15,19\}$, $\{15,20\}$, $\{16,19\}$, $\{16,21\}$, $\{17,20\}$, $\{17,22\}$, $\{18,21\}$, $\{18,22\}$,
   $\{19,22\}$, $\{20,21\}$\},

$\adfLe$:
  \{$\{1,2\}$, $\{1,3\}$, $\{1,4\}$, $\{2,5\}$, $\{2,6\}$, $\{3,7\}$, $\{3,8\}$, $\{4,9\}$, $\{4,10\}$, $\{5,7\}$, $\{5,9\}$,
   $\{6,8\}$, $\{6,11\}$, $\{7,12\}$, $\{8,13\}$, $\{9,14\}$, $\{10,11\}$, $\{10,12\}$, $\{11,15\}$, $\{12,16\}$,
   $\{13,14\}$, $\{13,17\}$, $\{14,18\}$, $\{15,19\}$, $\{15,20\}$, $\{16,21\}$, $\{16,22\}$, $\{17,19\}$, $\{17,21\}$,
   $\{18,20\}$, $\{18,22\}$, $\{19,22\}$, $\{20,21\}$\},

$\adfLf$:
  \{$\{1,2\}$, $\{1,3\}$, $\{1,4\}$, $\{2,5\}$, $\{2,6\}$, $\{3,7\}$, $\{3,8\}$, $\{4,9\}$, $\{4,10\}$, $\{5,7\}$, $\{5,9\}$,
   $\{6,8\}$, $\{6,11\}$, $\{7,12\}$, $\{8,13\}$, $\{9,14\}$, $\{10,11\}$, $\{10,12\}$, $\{11,15\}$, $\{12,16\}$,
   $\{13,17\}$, $\{13,18\}$, $\{14,19\}$, $\{14,20\}$, $\{15,16\}$, $\{15,21\}$, $\{16,22\}$, $\{17,19\}$, $\{17,21\}$,
   $\{18,20\}$, $\{18,22\}$, $\{19,22\}$, $\{20,21\}$\},

$\adfLg$:
  \{$\{1,2\}$, $\{1,3\}$, $\{1,4\}$, $\{2,5\}$, $\{2,6\}$, $\{3,7\}$, $\{3,8\}$, $\{4,9\}$, $\{4,10\}$, $\{5,7\}$, $\{5,9\}$,
   $\{6,8\}$, $\{6,11\}$, $\{7,12\}$, $\{8,13\}$, $\{9,14\}$, $\{10,11\}$, $\{10,15\}$, $\{11,16\}$,
   $\{12,15\}$, $\{12,16\}$, $\{13,17\}$, $\{13,18\}$, $\{14,19\}$, $\{14,20\}$, $\{15,21\}$, $\{16,22\}$, $\{17,19\}$, $\{17,21\}$,
   $\{18,20\}$, $\{18,22\}$, $\{19,22\}$, $\{20,21\}$\},

$\adfLh$:
  \{$\{1,2\}$, $\{1,3\}$, $\{1,4\}$, $\{2,5\}$, $\{2,6\}$, $\{3,7\}$, $\{3,8\}$, $\{4,9\}$, $\{4,10\}$, $\{5,7\}$, $\{5,9\}$,
   $\{6,8\}$, $\{6,11\}$, $\{7,12\}$, $\{8,13\}$, $\{9,14\}$, $\{10,12\}$, $\{10,15\}$, $\{11,12\}$, $\{11,16\}$,
   $\{13,14\}$, $\{13,17\}$, $\{14,18\}$, $\{15,19\}$, $\{15,20\}$, $\{16,21\}$, $\{16,22\}$, $\{17,19\}$, $\{17,21\}$,
   $\{18,20\}$, $\{18,22\}$, $\{19,22\}$, $\{20,21\}$\},

$\adfLi$:
  \{$\{1,2\}$, $\{1,3\}$, $\{1,4\}$, $\{2,5\}$, $\{2,6\}$, $\{3,7\}$, $\{3,8\}$, $\{4,9\}$, $\{4,10\}$, $\{5,7\}$, $\{5,9\}$,
   $\{6,8\}$, $\{6,11\}$, $\{7,12\}$, $\{8,13\}$, $\{9,14\}$, $\{10,12\}$, $\{10,15\}$, $\{11,12\}$, $\{11,16\}$,
   $\{13,17\}$, $\{13,18\}$, $\{14,19\}$, $\{14,20\}$, $\{15,16\}$, $\{15,21\}$, $\{16,22\}$, $\{17,19\}$, $\{17,21\}$,
   $\{18,20\}$, $\{18,22\}$, $\{19,22\}$, $\{20,21\}$\},

$\adfLj$:
  \{$\{1,2\}$, $\{1,3\}$, $\{1,4\}$, $\{2,5\}$, $\{2,6\}$, $\{3,7\}$, $\{3,8\}$, $\{4,9\}$, $\{4,10\}$, $\{5,7\}$, $\{5,9\}$,
  $\{6,8\}$, $\{6,11\}$, $\{7,12\}$, $\{8,13\}$, $\{9,14\}$, $\{10,12\}$, $\{10,15\}$, $\{11,15\}$, $\{11,16\}$, $\{12,16\}$,
  $\{13,17\}$, $\{13,18\}$, $\{14,19\}$, $\{14,20\}$, $\{15,21\}$, $\{16,22\}$, $\{17,19\}$, $\{17,21\}$,
  $\{18,20\}$, $\{18,22\}$, $\{19,22\}$, $\{20,21\}$\},

$\adfLk$:
  \{$\{1,2\}$, $\{1,3\}$, $\{1,4\}$, $\{2,5\}$, $\{2,6\}$, $\{3,7\}$, $\{3,8\}$, $\{4,9\}$, $\{4,10\}$, $\{5,7\}$, $\{5,9\}$,
   $\{6,11\}$, $\{6,12\}$, $\{7,10\}$, $\{8,11\}$, $\{8,13\}$, $\{9,13\}$, $\{10,14\}$, $\{11,15\}$,
   $\{12,14\}$, $\{12,16\}$, $\{13,17\}$, $\{14,18\}$, $\{15,19\}$, $\{15,20\}$, $\{16,19\}$, $\{16,21\}$, $\{17,21\}$, $\{17,22\}$,
   $\{18,20\}$, $\{18,22\}$, $\{19,22\}$, $\{20,21\}$\},

$\adfLl$:
  \{$\{1,2\}$, $\{1,3\}$, $\{1,4\}$, $\{2,5\}$, $\{2,6\}$, $\{3,7\}$, $\{3,8\}$, $\{4,9\}$, $\{4,10\}$, $\{5,7\}$, $\{5,9\}$,
   $\{6,11\}$, $\{6,12\}$, $\{7,10\}$, $\{8,11\}$, $\{8,13\}$, $\{9,13\}$, $\{10,14\}$, $\{11,15\}$,
   $\{12,16\}$, $\{12,17\}$, $\{13,18\}$, $\{14,19\}$, $\{14,20\}$, $\{15,18\}$, $\{15,21\}$, $\{16,19\}$, $\{16,21\}$, $\{17,20\}$, $\{17,22\}$,
   $\{18,22\}$, $\{19,22\}$, $\{20,21\}$\},

$\adfLm$:
  \{$\{1,2\}$, $\{1,3\}$, $\{1,4\}$, $\{2,5\}$, $\{2,6\}$, $\{3,7\}$, $\{3,8\}$, $\{4,9\}$, $\{4,10\}$, $\{5,7\}$, $\{5,9\}$,
   $\{6,11\}$, $\{6,12\}$, $\{7,10\}$, $\{8,11\}$, $\{8,13\}$, $\{9,14\}$, $\{10,15\}$, $\{11,14\}$,
   $\{12,15\}$, $\{12,16\}$, $\{13,17\}$, $\{13,18\}$, $\{14,19\}$, $\{15,20\}$, $\{16,17\}$, $\{16,21\}$, $\{17,22\}$,
   $\{18,20\}$, $\{18,21\}$, $\{19,21\}$, $\{19,22\}$, $\{20,22\}$\},

$\adfLn$:
  \{$\{1,2\}$, $\{1,3\}$, $\{1,4\}$, $\{2,5\}$, $\{2,6\}$, $\{3,7\}$, $\{3,8\}$, $\{4,9\}$, $\{4,10\}$, $\{5,7\}$, $\{5,9\}$,
   $\{6,11\}$, $\{6,12\}$, $\{7,10\}$, $\{8,11\}$, $\{8,13\}$, $\{9,14\}$, $\{10,15\}$, $\{11,14\}$,
   $\{12,16\}$, $\{12,17\}$, $\{13,18\}$, $\{13,19\}$, $\{14,18\}$, $\{15,20\}$, $\{15,21\}$, $\{16,19\}$, $\{16,20\}$, $\{17,21\}$, $\{17,22\}$,
   $\{18,22\}$, $\{19,21\}$, $\{20,22\}$\},

$\adfLo$:
  \{$\{1,2\}$, $\{1,3\}$, $\{1,4\}$, $\{2,5\}$, $\{2,6\}$, $\{3,7\}$, $\{3,8\}$, $\{4,9\}$, $\{4,10\}$, $\{5,7\}$, $\{5,9\}$,
   $\{6,11\}$, $\{6,12\}$, $\{7,10\}$, $\{8,11\}$, $\{8,13\}$, $\{9,14\}$, $\{10,15\}$, $\{11,16\}$,
   $\{12,15\}$, $\{12,17\}$, $\{13,14\}$, $\{13,18\}$, $\{14,19\}$, $\{15,20\}$, $\{16,21\}$, $\{16,22\}$, $\{17,18\}$, $\{17,21\}$,
   $\{18,22\}$, $\{19,20\}$, $\{19,21\}$, $\{20,22\}$\},

$\adfLp$:
  \{$\{1,2\}$, $\{1,3\}$, $\{1,4\}$, $\{2,5\}$, $\{2,6\}$, $\{3,7\}$, $\{3,8\}$, $\{4,9\}$, $\{4,10\}$, $\{5,7\}$, $\{5,9\}$,
   $\{6,11\}$, $\{6,12\}$, $\{7,10\}$, $\{8,11\}$, $\{8,13\}$, $\{9,14\}$, $\{10,15\}$, $\{11,16\}$,
   $\{12,15\}$, $\{12,17\}$, $\{13,14\}$, $\{13,18\}$, $\{14,19\}$, $\{15,20\}$, $\{16,21\}$, $\{16,22\}$, $\{17,19\}$, $\{17,21\}$,
   $\{18,20\}$, $\{18,21\}$, $\{19,22\}$, $\{20,22\}$\},

$\adfLq$:
  \{$\{1,2\}$, $\{1,3\}$, $\{1,4\}$, $\{2,5\}$, $\{2,6\}$, $\{3,7\}$, $\{3,8\}$, $\{4,9\}$, $\{4,10\}$, $\{5,7\}$, $\{5,9\}$,
   $\{6,11\}$, $\{6,12\}$, $\{7,10\}$, $\{8,11\}$, $\{8,13\}$, $\{9,14\}$, $\{10,15\}$, $\{11,16\}$,
   $\{12,15\}$, $\{12,17\}$, $\{13,18\}$, $\{13,19\}$, $\{14,20\}$, $\{14,21\}$, $\{15,16\}$, $\{16,22\}$, $\{17,18\}$, $\{17,20\}$,
   $\{18,21\}$, $\{19,20\}$, $\{19,22\}$, $\{21,22\}$\},

$\adfLr$:
  \{$\{1,2\}$, $\{1,3\}$, $\{1,4\}$, $\{2,5\}$, $\{2,6\}$, $\{3,7\}$, $\{3,8\}$, $\{4,9\}$, $\{4,10\}$, $\{5,7\}$, $\{5,9\}$,
   $\{6,11\}$, $\{6,12\}$, $\{7,10\}$, $\{8,13\}$, $\{8,14\}$, $\{9,15\}$, $\{10,16\}$, $\{11,13\}$, $\{11,16\}$,
   $\{12,17\}$, $\{12,18\}$, $\{13,15\}$, $\{14,17\}$, $\{14,19\}$, $\{15,20\}$, $\{16,21\}$, $\{17,22\}$,
   $\{18,19\}$, $\{18,20\}$, $\{19,21\}$, $\{20,22\}$, $\{21,22\}$\},

$\adfLs$:
  \{$\{1,2\}$, $\{1,3\}$, $\{1,4\}$, $\{2,5\}$, $\{2,6\}$, $\{3,7\}$, $\{3,8\}$, $\{4,9\}$, $\{4,10\}$, $\{5,7\}$, $\{5,9\}$,
   $\{6,11\}$, $\{6,12\}$, $\{7,10\}$, $\{8,13\}$, $\{8,14\}$, $\{9,15\}$, $\{10,16\}$, $\{11,13\}$, $\{11,16\}$, $\{12,17\}$,
   $\{12,18\}$, $\{13,19\}$, $\{14,15\}$, $\{14,17\}$, $\{15,20\}$, $\{16,21\}$, $\{17,22\}$,
   $\{18,19\}$, $\{18,20\}$, $\{19,22\}$, $\{20,21\}$, $\{21,22\}$\}
and

$\adfLt$:
  \{$\{1,2\}$, $\{1,3\}$, $\{1,4\}$, $\{2,5\}$, $\{2,6\}$, $\{3,7\}$, $\{3,8\}$, $\{4,9\}$, $\{4,10\}$, $\{5,7\}$, $\{5,9\}$,
  $\{6,11\}$, $\{6,12\}$, $\{7,10\}$, $\{8,13\}$, $\{8,14\}$, $\{9,15\}$, $\{10,16\}$, $\{11,13\}$, $\{11,17\}$,
  $\{12,16\}$, $\{12,18\}$, $\{13,19\}$, $\{14,15\}$, $\{14,18\}$, $\{15,20\}$, $\{16,21\}$, $\{17,20\}$, $\{17,22\}$,
  $\{18,22\}$, $\{19,21\}$, $\{19,22\}$, $\{20,21\}$\}.

The first two correspond to the Loupekine snarks \#1 and \#2 as supplied with the
{\sc Mathematica} system, \cite{Math-Snarks}. The remaining edge sets are precisely those of the
first eighteen 22-vertex snarks in Royale's list, \cite{Roy}.


{\lemma \label{lem:snark22 designs}
There exist designs of order $22$, $34$, $55$, $67$ and $88$ for each of the twenty non-trivial snarks on $22$ vertices.}\\


\noindent\textbf{Proof.}
Let the vertex set of $K_{22}$ be $Z_{21} \cup \{\infty\}$. The decompositions consist of

\adfLgap
$(\infty,13,0,7,3,14,18,12,4,17,5,2,11,15,6,10,8,1,9,19,16,20)_{\adfLa}$,

\adfLgap
$(\infty,5,4,1,13,7,18,11,17,6,0,15,3,19,16,20,12,8,2,10,14,9)_{\adfLb}$,

\adfLgap
$(\infty,13,3,5,17,8,2,7,19,18,11,20,6,14,1,15,16,12,0,4,10,9)_{\adfLc}$,

\adfLgap
$(\infty,9,10,11,4,14,13,12,18,3,15,16,8,0,6,17,20,1,19,2,7,5)_{\adfLd}$,

\adfLgap
$(\infty,9,5,13,19,3,18,16,14,1,8,20,12,15,11,6,0,10,7,2,17,4)_{\adfLe}$,

\adfLgap
$(\infty,16,18,2,9,0,12,10,14,20,4,3,17,6,19,1,15,11,5,7,8,13)_{\adfLf}$,

\adfLgap
$(\infty,18,10,20,16,13,7,11,2,14,9,1,19,3,5,17,8,12,6,4,15,0)_{\adfLg}$,

\adfLgap
$(\infty,20,19,0,2,4,12,6,14,9,15,11,3,8,17,16,18,7,13,10,1,5)_{\adfLh}$,

\adfLgap
$(\infty,1,14,12,15,19,16,13,4,10,7,2,17,9,20,3,11,5,0,8,6,18)_{\adfLi}$,

\adfLgap
$(\infty,12,2,13,14,1,11,7,16,15,6,5,20,3,9,18,19,0,8,17,10,4)_{\adfLj}$,

\adfLgap
$(\infty,14,16,6,15,9,0,18,8,13,5,12,19,3,4,2,17,7,20,1,10,11)_{\adfLk}$,

\adfLgap
$(\infty,9,20,1,19,10,11,0,3,7,17,13,14,15,12,8,4,16,2,18,6,5)_{\adfLl}$,

\adfLgap
$(\infty,18,17,7,8,20,11,16,13,3,15,10,19,1,2,14,5,6,12,4,0,9)_{\adfLm}$,

\adfLgap
$(\infty,12,10,2,4,20,0,7,16,9,5,3,8,1,6,18,14,17,19,13,11,15)_{\adfLn}$,

\adfLgap
$(\infty,9,10,17,4,1,3,13,16,7,20,5,19,11,18,8,2,0,6,15,12,14)_{\adfLo}$,

\adfLgap
$(\infty,20,9,19,6,2,7,1,8,4,10,16,15,3,0,5,14,12,13,11,18,17)_{\adfLp}$,

\adfLgap
$(\infty,11,19,9,4,10,13,16,17,3,12,5,20,7,8,2,6,1,18,15,0,14)_{\adfLq}$,

\adfLgap
$(\infty,6,13,17,4,2,20,14,7,11,5,12,1,0,16,3,18,19,9,15,8,10)_{\adfLr}$,

\adfLgap
$(\infty,1,3,14,15,16,7,8,12,11,19,4,0,13,18,20,2,17,10,6,5,9)_{\adfLs}$

\noindent and

$(\infty,6,7,14,20,17,19,11,12,16,8,1,5,4,10,0,3,15,9,2,13,18)_{\adfLt}$

\noindent under the action of the mapping $\infty \mapsto \infty$, $x \mapsto x + 3$ (mod 21).

\adfVfy{22, \{\{21,7,3\},\{1,1,1\}\}, 21, -1, -1} 

Let the vertex set of $K_{34}$ be $Z_{34}$. The decompositions consist of

\adfLgap
$(7,15,32,21,17,22,10,16,3,27,1,30,0,14,19,12,9,26,20,11,18,33)_{\adfLa}$,

\adfLgap
$(14,27,31,16,28,29,11,18,10,8,20,25,22,7,26,17,5,12,33,9,23,0)_{\adfLb}$,

\adfLgap
$(17,10,11,28,20,25,23,16,26,8,22,4,3,1,6,12,7,15,31,2,33,5)_{\adfLc}$,

\adfLgap
$(5,20,21,7,30,23,32,31,1,12,27,22,14,13,0,19,26,25,6,33,8,10)_{\adfLd}$,

\adfLgap
$(15,8,2,30,1,12,25,14,4,10,21,13,5,9,0,33,23,26,28,16,31,27)_{\adfLe}$,

\adfLgap
$(16,10,25,2,32,11,6,0,19,29,15,24,31,12,5,23,9,18,1,8,3,20)_{\adfLf}$,

\adfLgap
$(20,5,19,4,23,18,27,21,0,10,11,7,16,12,24,1,25,13,14,22,3,6)_{\adfLg}$,

\adfLgap
$(20,27,24,22,19,31,25,9,1,21,7,32,12,15,6,5,26,28,17,18,10,0)_{\adfLh}$,

\adfLgap
$(1,20,2,3,12,7,21,0,8,23,10,27,24,30,31,26,18,17,19,33,9,6)_{\adfLi}$,

\adfLgap
$(20,15,31,12,9,19,23,10,30,26,7,4,6,13,25,5,29,0,22,3,32,24)_{\adfLj}$,

\adfLgap
$(9,10,24,7,31,21,2,17,27,16,20,32,12,14,26,15,19,4,5,0,3,13)_{\adfLk}$,

\adfLgap
$(6,26,29,3,28,25,23,22,11,21,15,13,31,17,12,2,8,18,32,30,20,33)_{\adfLl}$,

\adfLgap
$(15,28,14,17,12,16,23,24,27,18,22,7,4,11,21,19,0,33,5,25,26,8)_{\adfLm}$,

\adfLgap
$(26,19,15,20,7,25,6,8,4,10,21,30,16,2,32,17,11,13,33,12,31,3)_{\adfLn}$,

\adfLgap
$(29,11,20,31,26,14,16,3,4,21,25,6,8,22,0,12,5,9,13,32,19,1)_{\adfLo}$,

\adfLgap
$(22,13,26,7,17,11,3,32,0,28,16,18,24,23,4,21,6,25,8,1,9,29)_{\adfLp}$,

\adfLgap
$(3,25,8,5,24,16,12,29,13,2,23,0,19,27,17,4,6,14,20,31,11,18)_{\adfLq}$,

\adfLgap
$(4,10,18,31,23,26,25,13,3,22,0,24,21,17,32,12,7,5,20,28,11,29)_{\adfLr}$,

\adfLgap
$(17,24,23,29,10,26,14,27,32,4,15,13,8,22,16,12,5,18,19,31,11,21)_{\adfLs}$

\noindent and

$(4,21,25,30,10,9,23,26,28,3,18,13,12,2,29,32,14,24,7,11,1,33)_{\adfLt}$

\noindent under the action of the mapping $x \mapsto x + 2$ (mod 34).

\adfVfy{34, \{\{34,17,2\}\}, -1, -1, -1} 

Let the vertex set of $K_{55}$ be $Z_{54} \cup \{\infty\}$. The decompositions consist of

\adfLgap
$(\infty,44,47,17,13,50,10,6,23,22,1,\adfsplit 4,3,8,31,15,53,33,45,5,24,20)_{\adfLa}$,

$(36,44,29,1,7,27,43,18,6,40,23,\adfsplit 15,31,48,30,2,37,53,0,9,19,26)_{\adfLa}$,

$(49,18,23,47,35,25,42,9,39,34,16,\adfsplit 43,20,50,24,44,26,22,53,46,36,27)_{\adfLa}$,

$(11,33,24,8,41,19,14,26,6,20,15,\adfsplit 9,44,25,0,52,27,16,45,4,29,40)_{\adfLa}$,

$(2,15,46,37,34,1,4,21,42,10,29,\adfsplit 53,35,40,\infty,16,25,7,6,38,33,0)_{\adfLa}$,

\adfLgap
$(\infty,22,27,6,0,23,18,11,48,49,17,\adfsplit 2,39,44,28,25,15,26,36,50,38,24)_{\adfLb}$,

$(40,52,30,5,51,6,25,37,1,34,26,\adfsplit 10,9,48,28,3,8,14,50,49,47,23)_{\adfLb}$,

$(19,21,17,43,26,46,27,14,13,7,0,\adfsplit 36,49,24,35,33,52,3,15,9,25,28)_{\adfLb}$,

$(47,20,30,42,9,23,49,29,48,53,31,\adfsplit 16,34,35,40,28,52,15,2,0,10,51)_{\adfLb}$,

$(3,2,7,11,46,31,23,16,43,9,21,\adfsplit 50,20,24,47,18,13,44,\infty,53,41,8)_{\adfLb}$,

\adfLgap
$(\infty,22,35,1,39,19,49,43,16,15,14,\adfsplit 50,45,5,29,2,18,21,10,11,42,6)_{\adfLc}$,

$(44,43,12,17,28,32,29,2,10,19,26,\adfsplit 3,4,14,47,35,53,37,31,33,1,24)_{\adfLc}$,

$(29,41,42,40,33,47,45,22,31,51,20,\adfsplit 50,52,53,39,9,30,16,4,13,18,46)_{\adfLc}$,

$(21,36,0,1,30,35,11,25,26,7,43,\adfsplit 33,9,34,42,44,12,48,51,50,24,4)_{\adfLc}$,

$(19,3,42,2,29,39,24,\infty,53,33,13,\adfsplit 50,40,46,1,10,20,23,36,52,14,38)_{\adfLc}$,

\adfLgap
$(\infty,0,46,14,16,53,42,47,11,18,43,\adfsplit 44,28,24,4,2,19,21,29,32,34,22)_{\adfLd}$,

$(34,46,39,7,10,13,25,19,23,36,20,\adfsplit 38,15,0,41,47,32,3,51,30,28,45)_{\adfLd}$,

$(31,39,43,21,46,6,35,14,3,47,53,\adfsplit 37,41,24,38,26,0,8,17,12,25,32)_{\adfLd}$,

$(2,46,3,52,12,49,33,19,4,31,48,\adfsplit 37,35,11,8,28,5,9,18,27,20,23)_{\adfLd}$,

$(26,31,40,32,48,24,16,9,3,50,19,\adfsplit 36,22,45,47,27,39,11,12,7,\infty,43)_{\adfLd}$,

\adfLgap
$(\infty,18,43,9,29,24,44,3,48,17,10,\adfsplit 49,20,46,41,52,15,16,27,0,21,38)_{\adfLe}$,

$(53,23,24,48,37,40,1,7,41,32,25,\adfsplit 44,9,47,0,13,43,3,27,28,21,22)_{\adfLe}$,

$(8,11,38,53,2,48,19,6,4,10,51,\adfsplit 18,30,15,21,36,32,47,39,20,45,5)_{\adfLe}$,

$(32,28,26,14,45,38,41,18,0,42,22,\adfsplit 49,31,19,12,6,37,43,34,16,50,5)_{\adfLe}$,

$(36,33,32,31,2,5,53,\infty,21,16,49,\adfsplit 17,10,23,29,52,46,25,1,50,4,51)_{\adfLe}$,

\adfLgap
$(\infty,25,45,14,39,48,7,53,11,32,27,\adfsplit 22,15,2,1,36,5,49,10,41,28,51)_{\adfLf}$,

$(6,19,24,15,2,48,45,44,16,14,29,\adfsplit 47,50,23,0,10,46,26,13,3,52,38)_{\adfLf}$,

$(46,7,35,27,28,3,48,17,39,19,38,\adfsplit 16,47,0,51,14,5,31,2,49,6,30)_{\adfLf}$,

$(43,29,19,50,4,23,12,32,1,36,30,\adfsplit 44,5,53,27,0,38,31,49,21,45,37)_{\adfLf}$,

$(33,27,4,30,23,10,0,26,\infty,41,22,\adfsplit 24,49,52,40,12,50,29,21,7,16,28)_{\adfLf}$,

\adfLgap
$(\infty,51,48,41,42,30,33,53,12,6,49,\adfsplit 4,19,8,37,52,31,32,46,43,16,44)_{\adfLg}$,

$(8,9,4,37,48,25,14,45,0,29,41,\adfsplit 36,26,13,23,38,17,12,39,15,52,16)_{\adfLg}$,

$(53,13,2,42,8,43,35,39,22,0,25,\adfsplit 31,47,10,40,24,9,23,28,38,20,21)_{\adfLg}$,

$(17,22,10,49,35,30,42,40,26,27,20,\adfsplit 4,23,51,47,13,6,33,39,38,32,3)_{\adfLg}$,

$(51,37,43,3,47,11,49,44,29,30,12,\adfsplit \infty,14,33,4,20,45,31,34,28,5,36)_{\adfLg}$,

\adfLgap
$(\infty,31,12,9,30,14,52,25,15,41,3,\adfsplit 23,33,5,35,2,21,17,19,43,28,20)_{\adfLh}$,

$(36,52,8,51,5,15,42,18,44,47,50,\adfsplit 45,6,30,49,17,40,14,27,7,12,37)_{\adfLh}$,

$(47,34,8,23,24,26,28,30,1,12,13,\adfsplit 33,49,19,45,4,42,3,15,32,35,46)_{\adfLh}$,

$(31,17,42,6,48,7,50,24,47,12,28,\adfsplit 10,33,2,4,29,52,14,43,20,25,46)_{\adfLh}$,

$(24,53,5,15,34,20,32,16,29,33,3,\adfsplit 41,26,\infty,37,4,50,22,31,51,25,28)_{\adfLh}$,

\adfLgap
$(\infty,12,32,10,16,5,7,13,14,33,37,\adfsplit 19,52,31,11,42,36,40,25,6,20,3)_{\adfLi}$,

$(52,41,9,44,53,39,37,5,47,45,43,\adfsplit 50,6,16,14,12,18,11,38,48,21,3)_{\adfLi}$,

$(51,0,4,8,13,30,23,48,2,50,22,\adfsplit 19,31,32,33,46,26,52,36,11,39,47)_{\adfLi}$,

$(44,53,26,15,23,46,40,52,49,51,6,\adfsplit 13,8,3,30,38,29,7,47,41,2,22)_{\adfLi}$,

$(41,26,43,27,23,46,42,13,37,25,45,\adfsplit 16,\infty,18,22,40,53,33,12,3,24,21)_{\adfLi}$,

\adfLgap
$(\infty,22,38,18,50,47,51,44,30,49,29,\adfsplit 6,5,48,21,39,20,3,4,1,16,23)_{\adfLj}$,

$(21,7,6,4,44,48,15,52,41,40,11,\adfsplit 20,39,35,37,34,0,17,25,13,26,22)_{\adfLj}$,

$(33,4,21,44,36,51,47,40,34,11,24,\adfsplit 49,41,1,52,27,29,8,43,18,20,7)_{\adfLj}$,

$(10,4,24,40,11,47,18,2,19,31,42,\adfsplit 37,20,46,22,17,13,32,44,15,21,48)_{\adfLj}$,

$(50,23,27,20,42,47,37,51,12,25,48,\adfsplit 21,32,33,49,35,18,7,15,26,34,\infty)_{\adfLj}$,

\adfLgap
$(\infty,18,28,31,24,29,12,2,23,34,7,\adfsplit 4,19,14,45,47,6,50,40,8,30,44)_{\adfLk}$,

$(37,44,11,15,28,43,10,25,16,18,4,\adfsplit 31,2,47,6,22,48,20,39,45,17,21)_{\adfLk}$,

$(38,41,40,17,27,21,36,30,7,32,12,\adfsplit 52,13,29,1,31,3,10,20,25,0,51)_{\adfLk}$,

$(39,25,28,43,7,18,42,31,41,14,33,\adfsplit 20,36,0,45,35,3,1,29,10,47,27)_{\adfLk}$,

$(25,38,52,33,48,3,22,18,41,44,53,\adfsplit 35,5,11,40,26,\infty,20,34,15,2,27)_{\adfLk}$,

\adfLgap
$(\infty,6,44,19,2,13,32,31,4,16,29,\adfsplit 3,10,45,53,20,7,36,48,34,38,35)_{\adfLl}$,

$(43,1,12,23,8,37,52,15,31,38,24,\adfsplit 7,35,44,25,17,39,46,11,16,34,4)_{\adfLl}$,

$(53,24,20,52,43,5,9,40,11,30,48,\adfsplit 50,6,8,31,41,13,26,15,21,45,51)_{\adfLl}$,

$(30,25,19,45,26,40,21,13,29,53,39,\adfsplit 49,2,28,34,52,20,10,15,18,36,6)_{\adfLl}$,

$(10,3,5,41,9,18,8,45,51,39,26,\adfsplit 12,24,\infty,7,35,50,36,53,28,47,6)_{\adfLl}$,

\adfLgap
$(\infty,14,31,11,47,5,53,17,30,40,29,\adfsplit 8,52,51,18,13,9,38,6,16,39,22)_{\adfLm}$,

$(45,37,32,17,41,49,5,19,30,9,29,\adfsplit 10,8,40,27,1,4,36,33,11,39,2)_{\adfLm}$,

$(28,33,23,44,47,8,3,25,48,13,30,\adfsplit 26,27,14,50,32,15,29,49,36,24,18)_{\adfLm}$,

$(37,7,4,11,53,27,46,40,50,42,32,\adfsplit 49,13,17,34,9,20,10,1,44,36,19)_{\adfLm}$,

$(17,16,34,40,39,36,9,15,0,10,6,\adfsplit 48,8,7,31,\infty,27,20,1,18,46,12)_{\adfLm}$,

\adfLgap
$(\infty,4,15,47,26,44,31,40,52,21,38,\adfsplit 42,41,51,35,9,22,34,2,27,10,48)_{\adfLn}$,

$(2,31,52,21,8,28,3,42,7,53,19,\adfsplit 9,4,41,13,16,32,6,22,49,5,24)_{\adfLn}$,

$(7,52,17,15,8,41,23,38,12,18,21,\adfsplit 32,2,34,25,3,20,19,47,30,36,6)_{\adfLn}$,

$(44,12,41,45,13,18,20,3,2,30,37,\adfsplit 38,34,40,47,14,25,1,42,49,29,51)_{\adfLn}$,

$(36,29,11,27,30,52,53,35,42,51,32,\adfsplit 28,19,\infty,45,24,41,49,7,33,13,22)_{\adfLn}$,

\adfLgap
$(\infty,38,36,47,0,7,29,10,14,43,20,\adfsplit 45,39,2,22,15,25,37,26,8,51,53)_{\adfLo}$,

$(50,51,9,30,23,5,29,48,0,12,52,\adfsplit 13,36,3,10,16,47,15,26,32,35,25)_{\adfLo}$,

$(23,33,47,41,34,42,32,51,27,28,37,\adfsplit 10,40,46,7,36,15,3,17,18,19,49)_{\adfLo}$,

$(1,29,44,36,19,24,28,9,14,8,20,\adfsplit 34,15,17,18,37,4,45,6,52,31,13)_{\adfLo}$,

$(48,50,21,35,31,49,\infty,36,13,52,42,\adfsplit 20,17,28,47,51,41,40,16,8,44,0)_{\adfLo}$,

\adfLgap
$(\infty,29,39,8,3,26,53,30,35,23,18,\adfsplit 45,43,19,16,11,34,4,47,13,14,46)_{\adfLp}$,

$(20,42,28,29,23,10,19,27,21,38,48,\adfsplit 14,2,39,1,13,8,45,9,24,31,6)_{\adfLp}$,

$(35,14,41,53,9,16,2,19,43,15,5,\adfsplit 24,18,23,22,36,28,32,52,42,11,31)_{\adfLp}$,

$(11,53,12,45,36,22,27,42,33,1,16,\adfsplit 38,29,39,2,21,10,37,0,26,25,40)_{\adfLp}$,

$(47,49,52,36,46,1,10,25,9,\infty,50,\adfsplit 30,27,40,13,0,24,19,8,51,26,20)_{\adfLp}$,

\adfLgap
$(\infty,41,44,52,25,0,53,13,15,29,51,\adfsplit 40,14,39,35,37,23,22,18,26,19,48)_{\adfLq}$,

$(40,38,48,16,36,41,31,11,37,18,19,\adfsplit 51,22,32,30,4,8,45,49,25,33,29)_{\adfLq}$,

$(12,17,45,43,52,13,18,53,39,32,38,\adfsplit 6,42,46,10,0,25,3,36,27,21,35)_{\adfLq}$,

$(21,23,12,18,39,13,5,41,8,38,28,\adfsplit 22,15,44,6,32,1,7,4,14,25,11)_{\adfLq}$,

$(12,35,48,34,8,47,33,\infty,20,39,15,\adfsplit 29,31,46,44,50,38,22,19,45,10,16)_{\adfLq}$,

\adfLgap
$(\infty,46,45,0,38,4,11,9,31,43,41,\adfsplit 30,48,33,6,39,26,20,1,51,53,37)_{\adfLr}$,

$(6,38,36,42,5,41,20,53,40,29,15,\adfsplit 9,7,0,45,21,22,52,37,28,1,35)_{\adfLr}$,

$(0,10,4,29,43,46,38,42,47,52,2,\adfsplit 37,1,15,36,19,27,31,5,17,21,28)_{\adfLr}$,

$(21,50,46,0,32,26,6,41,51,1,39,\adfsplit 27,36,18,28,4,38,48,3,2,8,40)_{\adfLr}$,

$(21,38,49,17,7,36,19,22,20,16,30,\adfsplit 50,37,31,8,10,41,44,\infty,23,35,47)_{\adfLr}$,

\adfLgap
$(\infty,11,30,8,34,43,33,50,2,7,4,\adfsplit 1,37,42,49,5,44,26,48,14,10,13)_{\adfLs}$,

$(9,13,49,18,36,43,32,5,31,2,28,\adfsplit 17,44,3,25,42,14,6,10,24,21,27)_{\adfLs}$,

$(45,33,8,25,12,47,27,41,28,5,53,\adfsplit 10,49,48,17,9,40,29,4,32,34,52)_{\adfLs}$,

$(26,51,17,44,10,20,36,35,0,34,28,\adfsplit 29,6,14,41,39,9,3,30,11,23,24)_{\adfLs}$,

$(50,24,10,1,36,49,45,30,41,37,51,\adfsplit 22,13,47,34,27,31,15,46,8,9,\infty)_{\adfLs}$

\noindent and

$(\infty,46,0,7,5,26,32,37,45,29,51,\adfsplit 28,12,34,18,14,33,22,10,40,9,1)_{\adfLt}$,

$(33,27,21,35,14,10,49,34,18,2,46,\adfsplit 31,19,17,8,32,23,25,22,51,37,11)_{\adfLt}$,

$(45,26,43,0,9,16,34,19,5,4,11,\adfsplit 27,24,39,15,20,13,30,48,25,6,50)_{\adfLt}$,

$(32,1,26,4,11,42,8,37,30,9,14,\adfsplit 49,53,6,16,13,43,29,47,51,17,0)_{\adfLt}$,

$(49,8,48,23,50,40,27,0,4,24,47,\adfsplit 44,15,21,29,6,26,\infty,53,18,52,17)_{\adfLt}$

\noindent under the action of the mapping $\infty \mapsto \infty$, $x \mapsto x + 6$ (mod 54).

\adfVfy{55, \{\{54,9,6\},\{1,1,1\}\}, 54, -1, -1} 

Let the vertex set of $K_{67}$ be $Z_{67}$. The decompositions consist of

\adfLgap
$(40,6,35,50,31,10,39,16,8,2,49,\adfsplit 9,63,55,44,0,19,7,24,64,22,29)_{\adfLa}$,

\adfLgap
$(26,27,60,55,40,38,28,31,17,47,32,\adfsplit 48,23,29,43,3,7,30,13,0,24,5)_{\adfLb}$,

\adfLgap
$(21,61,36,9,13,31,12,22,44,25,54,\adfsplit 26,52,32,59,3,18,11,41,20,0,28)_{\adfLc}$,

\adfLgap
$(25,53,45,54,41,55,60,58,27,37,63,\adfsplit 23,16,59,64,2,0,18,7,6,36,43)_{\adfLd}$,

\adfLgap
$(56,58,43,0,40,63,44,33,25,24,47,\adfsplit 5,36,19,15,12,3,27,61,1,32,39)_{\adfLe}$,

\adfLgap
$(7,22,25,21,15,30,56,27,42,8,2,\adfsplit 45,43,51,19,23,31,33,32,28,66,57)_{\adfLf}$,

\adfLgap
$(54,43,35,10,27,47,21,2,30,41,26,\adfsplit 9,15,48,14,0,16,13,55,56,46,38)_{\adfLg}$,

\adfLgap
$(64,28,39,15,37,2,60,1,27,9,47,\adfsplit 33,16,8,41,0,21,4,58,43,54,7)_{\adfLh}$,

\adfLgap
$(59,22,20,61,11,7,12,16,37,54,39,\adfsplit 57,32,47,33,6,51,46,9,52,64,29)_{\adfLi}$,

\adfLgap
$(43,52,64,50,34,57,20,45,35,24,23,\adfsplit 49,51,5,59,6,4,0,1,13,2,28)_{\adfLj}$,

\adfLgap
$(14,7,59,30,34,57,3,2,55,62,60,\adfsplit 53,49,48,31,1,19,5,50,44,45,17)_{\adfLk}$,

\adfLgap
$(30,7,24,17,59,45,40,3,35,47,56,\adfsplit 20,61,50,11,10,21,66,46,23,43,38)_{\adfLl}$,

\adfLgap
$(32,28,18,61,52,2,45,31,27,44,49,\adfsplit 24,23,59,47,3,8,62,57,11,51,20)_{\adfLm}$,

\adfLgap
$(46,55,61,13,41,60,48,30,2,3,59,\adfsplit 28,22,43,45,1,11,40,20,21,57,17)_{\adfLn}$,

\adfLgap
$(11,61,47,0,40,20,59,49,37,52,55,\adfsplit 10,62,3,53,2,1,24,42,26,64,44)_{\adfLo}$,

\adfLgap
$(0,28,63,57,2,35,50,48,33,25,65,\adfsplit 55,45,31,37,3,21,22,10,38,30,32)_{\adfLp}$,

\adfLgap
$(27,56,35,0,60,21,63,20,53,18,11,\adfsplit 15,40,2,51,10,34,23,3,57,48,5)_{\adfLq}$,

\adfLgap
$(10,3,48,20,50,40,17,23,59,19,44,\adfsplit 55,39,31,45,0,34,1,13,28,24,2)_{\adfLr}$,

\adfLgap
$(45,28,35,41,52,64,20,1,57,60,50,\adfsplit 51,24,23,48,4,21,54,53,42,22,14)_{\adfLs}$

\noindent and

$(54,44,27,10,41,9,52,7,36,48,17,\adfsplit 42,38,23,45,5,11,40,26,30,12,39)_{\adfLt}$

\noindent under the action of the mapping $x \mapsto x + 1$ (mod 67).

\adfVfy{67, \{\{67,67,1\}\}, -1, -1, -1} 

Let the vertex set of $K_{88}$ be $Z_{87} \cup \{\infty\}$. The decompositions consist of

\adfLgap
$(9,84,83,29,63,44,22,46,14,42,86,\adfsplit 52,1,72,64,74,68,45,30,55,75,54)_{\adfLa}$,

$(27,74,62,24,79,40,22,50,33,7,17,\adfsplit 11,65,73,9,86,55,20,15,66,49,43)_{\adfLa}$,

$(58,59,55,30,40,31,1,71,32,46,47,\adfsplit 57,82,0,29,39,81,20,56,78,86,43)_{\adfLa}$,

$(73,2,24,6,12,52,10,4,13,83,38,\adfsplit \infty,78,27,69,17,56,42,49,18,22,20)_{\adfLa}$,

\adfLgap
$(20,47,84,50,57,16,22,26,37,9,46,\adfsplit 24,76,74,80,7,30,21,56,54,8,49)_{\adfLb}$,

$(86,9,17,52,40,73,21,62,84,24,35,\adfsplit 47,71,82,1,30,57,27,59,19,2,55)_{\adfLb}$,

$(33,53,71,47,62,18,0,21,82,45,55,\adfsplit 75,57,42,85,13,17,16,26,83,19,9)_{\adfLb}$,

$(41,15,1,83,19,72,10,4,0,42,27,\adfsplit 6,8,50,85,13,7,64,32,65,\infty,70)_{\adfLb}$,

\adfLgap
$(\infty,34,53,39,7,75,40,0,14,24,42,\adfsplit 5,82,74,59,36,76,63,23,21,84,20)_{\adfLc}$,

$(10,45,61,33,18,7,31,79,15,41,39,\adfsplit 57,86,19,80,1,63,47,17,4,48,64)_{\adfLc}$,

$(0,45,57,34,78,41,67,71,26,80,35,\adfsplit 50,25,54,79,85,13,56,40,75,82,11)_{\adfLc}$,

$(34,18,71,83,13,55,4,76,77,79,30,\adfsplit 3,23,21,36,22,14,86,8,29,39,26)_{\adfLc}$,

\adfLgap
$(\infty,86,55,54,25,0,37,39,14,10,31,\adfsplit 11,81,32,78,36,29,59,51,49,64,53)_{\adfLd}$,

$(25,62,39,46,44,32,15,19,80,5,0,\adfsplit 16,52,42,38,48,75,8,40,73,54,70)_{\adfLd}$,

$(60,65,9,56,81,17,78,10,2,85,39,\adfsplit 64,70,8,6,15,46,52,72,29,61,7)_{\adfLd}$,

$(63,54,28,44,37,64,68,3,0,2,60,\adfsplit 15,26,73,7,8,22,84,80,86,11,71)_{\adfLd}$,

\adfLgap
$(\infty,64,84,17,5,6,14,38,39,69,27,\adfsplit 49,79,33,32,53,37,80,47,21,60,66)_{\adfLe}$,

$(4,62,37,17,29,24,53,64,56,33,49,\adfsplit 84,9,7,47,80,72,0,61,26,73,79)_{\adfLe}$,

$(79,56,4,64,50,68,60,39,47,13,34,\adfsplit 63,41,11,73,10,42,16,30,8,75,53)_{\adfLe}$,

$(13,10,19,72,84,35,15,80,12,46,55,\adfsplit 77,25,42,69,27,44,61,26,31,36,37)_{\adfLe}$,

\adfLgap
$(\infty,14,43,66,73,3,59,61,33,11,30,\adfsplit 39,48,12,34,38,37,60,31,62,8,68)_{\adfLf}$,

$(8,61,49,39,64,12,79,13,21,78,84,\adfsplit 71,27,18,59,45,73,43,5,82,7,68)_{\adfLf}$,

$(85,40,7,75,71,62,66,20,35,82,2,\adfsplit 3,30,39,18,70,13,68,37,80,63,15)_{\adfLf}$,

$(38,44,43,37,59,5,16,26,70,6,63,\adfsplit 71,49,23,27,36,75,11,13,20,22,80)_{\adfLf}$,

\adfLgap
$(\infty,85,11,75,10,47,67,13,33,22,15,\adfsplit 62,72,58,0,80,66,6,38,86,57,64)_{\adfLg}$,

$(60,64,49,42,65,32,26,5,56,4,11,\adfsplit 29,6,70,3,9,33,57,30,52,16,14)_{\adfLg}$,

$(54,8,69,2,78,19,34,79,66,65,28,\adfsplit 36,64,82,23,4,67,25,5,51,59,45)_{\adfLg}$,

$(69,30,20,59,41,67,1,50,37,65,32,\adfsplit 7,35,79,53,61,39,66,5,58,72,42)_{\adfLg}$,

\adfLgap
$(\infty,80,34,66,26,54,4,85,6,37,79,\adfsplit 45,78,75,84,11,5,30,35,51,31,86)_{\adfLh}$,

$(73,84,44,86,48,82,34,66,22,40,67,\adfsplit 30,3,9,29,7,25,32,64,11,70,85)_{\adfLh}$,

$(79,76,65,15,82,80,63,71,11,67,35,\adfsplit 62,39,74,55,28,0,77,17,38,53,33)_{\adfLh}$,

$(22,81,39,14,74,35,24,30,27,42,85,\adfsplit 53,38,57,54,36,37,13,55,44,4,11)_{\adfLh}$,

\adfLgap
$(\infty,35,1,27,22,82,3,18,76,51,54,\adfsplit 64,32,79,46,26,61,10,4,20,81,31)_{\adfLi}$,

$(43,5,67,18,19,6,23,21,74,38,75,\adfsplit 78,25,20,64,62,7,41,14,84,13,80)_{\adfLi}$,

$(69,17,27,39,21,16,60,15,28,14,35,\adfsplit 65,46,62,72,45,71,54,47,20,80,1)_{\adfLi}$,

$(82,43,24,37,57,35,79,33,68,14,47,\adfsplit 64,39,65,49,58,77,2,84,48,12,18)_{\adfLi}$,

\adfLgap
$(\infty,37,21,44,47,23,46,52,33,76,1,\adfsplit 58,16,14,84,28,20,60,18,77,48,73)_{\adfLj}$,

$(12,52,85,60,72,2,29,74,54,55,19,\adfsplit 21,18,44,9,73,51,25,81,38,79,36)_{\adfLj}$,

$(75,8,48,53,28,60,22,47,6,41,83,\adfsplit 7,65,27,76,45,17,20,51,11,55,77)_{\adfLj}$,

$(60,72,61,69,0,1,64,3,62,23,86,\adfsplit 83,6,65,70,58,11,25,44,73,47,49)_{\adfLj}$,

\adfLgap
$(\infty,24,65,76,47,8,53,56,12,0,40,\adfsplit 29,41,72,16,32,74,18,73,44,61,11)_{\adfLk}$,

$(1,15,10,9,2,47,71,29,40,66,37,\adfsplit 72,80,13,49,83,7,42,31,54,81,48)_{\adfLk}$,

$(69,6,2,19,22,58,42,48,25,14,73,\adfsplit 7,67,45,50,41,47,63,43,72,51,16)_{\adfLk}$,

$(64,63,60,53,71,31,34,41,23,36,85,\adfsplit 32,83,81,37,18,79,80,15,30,69,58)_{\adfLk}$,

\adfLgap
$(\infty,57,86,85,34,23,7,54,47,4,8,\adfsplit 40,3,43,30,62,19,51,73,14,65,31)_{\adfLl}$,

$(14,78,67,41,60,19,74,49,80,24,18,\adfsplit 9,5,17,61,82,6,28,42,21,47,33)_{\adfLl}$,

$(64,72,81,15,23,42,80,27,17,19,40,\adfsplit 48,8,35,20,59,73,53,31,10,41,66)_{\adfLl}$,

$(25,75,76,78,35,30,67,0,66,68,59,\adfsplit 85,71,60,83,24,13,73,82,84,43,32)_{\adfLl}$,

\adfLgap
$(\infty,12,50,67,6,41,4,32,82,70,86,\adfsplit 85,15,20,69,39,25,51,33,34,24,30)_{\adfLm}$,

$(23,62,50,70,85,77,86,66,51,56,9,\adfsplit 15,79,27,65,41,47,55,7,73,30,16)_{\adfLm}$,

$(33,80,66,29,5,78,7,32,72,49,39,\adfsplit 34,58,17,82,40,3,51,76,86,2,45)_{\adfLm}$,

$(39,62,27,31,6,25,67,17,14,9,20,\adfsplit 37,52,15,75,64,74,71,29,47,13,40)_{\adfLm}$,

\adfLgap
$(\infty,7,56,60,20,0,29,19,14,45,55,\adfsplit 62,25,28,81,80,30,35,17,31,38,76)_{\adfLn}$,

$(63,40,60,80,27,29,53,81,37,84,64,\adfsplit 32,69,46,31,61,17,47,2,45,3,70)_{\adfLn}$,

$(37,59,20,56,55,78,57,60,8,44,13,\adfsplit 9,18,38,42,4,66,70,72,16,77,27)_{\adfLn}$,

$(27,80,50,36,19,13,14,5,40,16,61,\adfsplit 85,38,69,79,54,52,9,48,30,49,68)_{\adfLn}$,

\adfLgap
$(\infty,13,66,80,64,48,60,14,55,69,62,\adfsplit 76,51,50,79,58,84,39,12,41,17,72)_{\adfLo}$,

$(24,50,3,55,63,72,40,21,38,68,1,\adfsplit 2,81,70,36,43,26,49,7,6,76,9)_{\adfLo}$,

$(79,33,32,52,26,17,53,77,37,30,68,\adfsplit 74,23,31,73,16,13,20,18,25,42,4)_{\adfLo}$,

$(21,17,29,70,23,57,41,44,13,51,42,\adfsplit 11,1,86,76,3,40,20,27,78,71,73)_{\adfLo}$,

\adfLgap
$(\infty,81,23,67,8,52,12,13,57,6,16,\adfsplit 19,44,70,43,41,18,3,36,35,47,63)_{\adfLp}$,

$(4,13,73,76,68,69,51,71,60,55,37,\adfsplit 30,49,86,27,67,15,2,82,29,38,40)_{\adfLp}$,

$(74,6,11,62,55,83,20,68,49,17,33,\adfsplit 30,52,8,37,78,35,80,75,85,16,21)_{\adfLp}$,

$(10,85,57,11,15,6,78,16,4,42,86,\adfsplit 49,2,18,54,47,51,35,70,53,60,32)_{\adfLp}$,

\adfLgap
$(\infty,30,77,40,70,64,37,31,2,66,73,\adfsplit 29,74,24,23,14,9,41,19,3,54,16)_{\adfLq}$,

$(46,56,57,50,54,61,67,24,66,60,47,\adfsplit 83,75,20,86,35,30,3,13,73,27,28)_{\adfLq}$,

$(20,47,28,60,70,18,48,79,58,25,55,\adfsplit 68,17,71,83,16,37,6,56,26,1,84)_{\adfLq}$,

$(77,69,11,61,37,68,74,30,67,78,12,\adfsplit 54,76,39,29,47,51,3,55,56,81,73)_{\adfLq}$,

\adfLgap
$(\infty,44,57,31,14,13,5,33,74,60,51,\adfsplit 47,48,6,16,56,81,45,40,65,49,19)_{\adfLr}$,

$(61,62,79,63,85,33,4,29,57,27,76,\adfsplit 71,30,25,80,36,6,12,22,59,43,83)_{\adfLr}$,

$(76,4,62,8,60,24,64,17,43,27,16,\adfsplit 45,5,3,67,11,20,53,59,66,37,79)_{\adfLr}$,

$(42,24,53,12,57,76,9,35,73,85,5,\adfsplit 68,2,81,34,25,20,32,39,17,61,7)_{\adfLr}$,

\adfLgap
$(\infty,77,24,61,34,16,52,20,69,58,83,\adfsplit 86,81,29,62,71,54,57,82,35,84,7)_{\adfLs}$,

$(75,21,43,80,51,70,65,28,59,46,72,\adfsplit 74,52,19,3,85,15,67,45,48,54,27)_{\adfLs}$,

$(24,60,72,65,63,81,46,83,23,10,52,\adfsplit 38,53,48,34,47,85,18,68,86,76,43)_{\adfLs}$,

$(73,52,84,17,80,79,74,6,31,56,8,\adfsplit 38,78,59,61,72,62,39,34,15,49,26)_{\adfLs}$

\noindent and

$(\infty,62,76,45,15,13,72,25,60,66,23,\adfsplit 28,37,29,17,75,7,83,42,64,78,31)_{\adfLt}$,

$(46,45,20,59,83,12,81,0,11,13,71,\adfsplit 82,48,55,38,57,2,53,29,75,49,60)_{\adfLt}$,

$(27,14,17,5,80,76,74,70,61,29,18,\adfsplit 28,52,3,16,30,59,21,53,22,57,31)_{\adfLt}$,

$(28,5,61,35,30,16,58,12,38,79,66,\adfsplit 8,83,65,29,63,25,60,72,34,7,48)_{\adfLt}$

\noindent under the action of the mapping $\infty \mapsto \infty$, $x \mapsto x + 3$ (mod 87).\eproof

\adfVfy{88, \{\{87,29,3\},\{1,1,1\}\}, 87, -1, -1} 

{\lemma \label{lem:snark22 multipartite}
There exist decompositions of
$K_{33,33,33}$, $K_{11,11,11,11}$, $K_{22,22,22,55}$ and $K_{22,22,22,22,22,22,55}$
into each of the twenty non-trivial snarks on $22$ vertices.}\\


\noindent\textbf{Proof.}
Let the vertex set of $K_{33,33,33}$ be $Z_{99}$ partitioned according to residue classes modulo 3.
The decompositions consist of

\adfLgap
$(10,29,62,52,93,55,96,73,0,35,85,\adfsplit 77,24,19,59,3,57,88,98,31,61,18)_{\adfLa}$,

\adfLgap
$(91,6,15,46,26,16,41,63,38,84,58,\adfsplit 54,13,56,82,11,17,48,40,78,18,43)_{\adfLb}$,

\adfLgap
$(44,12,66,16,86,88,70,26,33,96,50,\adfsplit 31,20,94,1,0,72,36,65,71,10,79)_{\adfLc}$,

\adfLgap
$(8,24,21,75,31,83,68,97,80,37,81,\adfsplit 2,55,71,4,3,11,63,32,84,20,52)_{\adfLd}$,

\adfLgap
$(12,40,80,64,96,57,46,25,38,68,97,\adfsplit 69,11,16,77,5,13,8,3,70,51,91)_{\adfLe}$,

\adfLgap
$(61,56,11,50,12,63,10,88,91,79,23,\adfsplit 62,84,78,9,25,8,13,97,5,40,59)_{\adfLf}$,

\adfLgap
$(12,19,23,50,60,5,55,4,85,90,88,\adfsplit 32,38,41,28,6,10,18,87,70,20,74)_{\adfLg}$,

\adfLgap
$(94,6,74,89,76,20,15,22,57,55,78,\adfsplit 77,95,85,32,10,4,72,39,28,53,56)_{\adfLh}$,

\adfLgap
$(47,12,25,88,56,2,42,90,63,50,96,\adfsplit 10,92,20,81,1,39,64,40,93,31,17)_{\adfLi}$,

\adfLgap
$(84,71,17,2,16,61,90,9,39,1,63,\adfsplit 47,28,53,41,13,56,33,60,91,21,44)_{\adfLj}$,

\adfLgap
$(54,11,7,91,18,97,92,45,62,81,47,\adfsplit 74,61,46,13,6,15,29,32,33,73,37)_{\adfLk}$,

\adfLgap
$(78,70,65,94,89,66,40,24,84,74,23,\adfsplit 55,61,60,76,0,57,29,38,88,59,64)_{\adfLl}$,

\adfLgap
$(65,84,45,22,74,70,67,73,87,96,75,\adfsplit 69,18,64,29,7,68,71,38,13,3,21)_{\adfLm}$,

\adfLgap
$(16,42,62,59,86,58,52,12,45,27,80,\adfsplit 56,50,85,19,4,55,57,21,0,32,20)_{\adfLn}$,

\adfLgap
$(55,83,26,87,90,52,88,39,11,44,74,\adfsplit 18,79,93,70,16,28,29,40,59,24,54)_{\adfLo}$,

\adfLgap
$(19,89,84,24,63,36,25,7,79,35,86,\adfsplit 44,50,54,3,0,45,78,68,95,97,31)_{\adfLp}$,

\adfLgap
$(66,52,40,82,23,48,18,65,81,35,67,\adfsplit 94,96,32,3,47,2,58,68,91,69,34)_{\adfLq}$,

\adfLgap
$(89,1,3,4,0,27,83,73,32,63,58,\adfsplit 37,50,48,28,65,95,98,43,81,9,46)_{\adfLr}$,

\adfLgap
$(87,44,38,89,6,36,10,25,52,75,77,\adfsplit 65,72,92,69,55,81,34,8,59,0,1)_{\adfLs}$

\noindent and

$(75,20,13,97,15,96,47,9,44,54,70,\adfsplit 32,83,56,45,4,53,42,81,61,23,22)_{\adfLt}$

\noindent under the action of the mapping $x \mapsto x + 1$ (mod 99).

\adfVfy{99, \{\{99,99,1\}\}, -1, \{\{33,\{0,1,2\}\}\}, -1} 

Let the vertex set of $K_{11,11,11,11}$ be $Z_{44}$ partitioned according to residue classes modulo 4.
The decompositions consist of

\adfLgap
$(24,35,14,13,11,41,23,4,18,5,26,\adfsplit 17,19,0,33,1,2,7,8,27,29,6)_{\adfLa}$,

$(0,1,2,4,7,13,31,6,11,22,36,\adfsplit 25,12,35,20,34,3,30,37,16,28,42)_{\adfLa}$

\adfLgap
$(26,35,23,32,39,6,1,2,24,25,42,\adfsplit 19,28,15,33,0,5,3,21,10,31,4)_{\adfLb}$,

$(0,2,5,3,1,8,24,7,4,34,43,\adfsplit 30,13,40,33,28,17,42,27,6,22,29)_{\adfLb}$

\adfLgap
$(18,28,27,41,19,35,33,8,6,14,23,\adfsplit 9,12,25,22,0,1,3,11,7,10,16)_{\adfLc}$,

$(0,2,3,5,1,4,24,22,19,30,38,\adfsplit 29,14,8,25,12,9,21,27,35,18,28)_{\adfLc}$,

\adfLgap
$(16,31,11,38,34,13,5,22,39,8,36,\adfsplit 9,14,6,35,0,4,12,1,21,2,23)_{\adfLd}$,

$(0,3,5,6,1,9,2,35,32,23,42,\adfsplit 4,39,25,7,31,10,36,38,28,21,17)_{\adfLd}$,

\adfLgap
$(20,22,15,11,27,12,0,25,33,26,37,\adfsplit 5,18,31,8,2,1,6,7,39,28,21)_{\adfLe}$,

$(0,1,3,6,2,4,24,30,12,43,25,\adfsplit 22,16,5,39,17,35,15,38,13,26,8)_{\adfLe}$,

\adfLgap
$(22,21,17,11,32,7,10,36,30,29,0,\adfsplit 19,35,33,3,1,4,6,14,27,5,24)_{\adfLf}$,

$(0,2,5,6,3,7,12,14,33,20,38,\adfsplit 31,11,39,9,36,28,32,22,29,42,1)_{\adfLf}$,

\adfLgap
$(16,42,7,14,23,13,21,32,25,33,30,\adfsplit 15,10,3,12,0,1,4,2,6,11,9)_{\adfLg}$,

$(0,1,7,10,4,6,17,37,27,43,12,\adfsplit 32,20,18,25,2,3,38,41,31,36,23)_{\adfLg}$,

\adfLgap
$(27,14,13,22,31,28,40,26,12,4,23,\adfsplit 17,24,35,30,1,2,5,9,7,8,0)_{\adfLh}$,

$(0,3,6,7,2,5,1,15,13,33,16,\adfsplit 30,18,28,42,43,37,9,27,35,10,12)_{\adfLh}$,

\adfLgap
$(12,11,35,5,42,14,13,33,0,39,37,\adfsplit 2,34,9,24,6,4,7,15,20,1,32)_{\adfLi}$,

$(0,1,2,3,6,15,7,8,9,29,26,\adfsplit 12,42,24,39,4,23,19,37,34,41,10)_{\adfLi}$,

\adfLgap
$(38,28,12,11,15,18,10,17,40,4,25,\adfsplit 21,8,37,34,0,1,2,3,23,32,29)_{\adfLj}$,

$(0,1,3,6,7,12,8,34,9,36,30,\adfsplit 19,21,35,15,5,11,26,2,41,18,27)_{\adfLj}$,

\adfLgap
$(39,25,4,37,35,34,26,38,22,19,23,\adfsplit 36,28,1,8,3,2,0,5,7,13,27)_{\adfLk}$,

$(0,2,3,5,1,7,6,28,24,20,34,\adfsplit 12,17,43,9,25,10,13,22,30,31,40)_{\adfLk}$,

\adfLgap
$(5,24,42,35,27,19,37,13,16,34,22,\adfsplit 30,2,21,32,3,4,1,12,31,9,18)_{\adfLl}$,

$(0,1,2,5,3,4,21,7,26,19,33,\adfsplit 37,28,10,40,43,20,27,36,35,6,14)_{\adfLl}$,

\adfLgap
$(7,2,9,6,8,29,19,4,21,32,18,\adfsplit 36,13,20,38,11,14,3,1,12,22,27)_{\adfLm}$,

$(0,3,7,10,1,4,12,13,24,41,27,\adfsplit 31,43,21,22,9,6,26,14,17,36,35)_{\adfLm}$,

\adfLgap
$(11,34,20,37,39,24,41,7,36,3,42,\adfsplit 29,25,31,13,0,2,6,14,15,16,1)_{\adfLn}$,

$(0,2,6,7,3,5,28,13,10,22,40,\adfsplit 26,24,23,9,41,36,37,35,16,18,15)_{\adfLn}$,

\adfLgap
$(36,9,23,25,31,10,29,0,12,20,41,\adfsplit 13,14,27,42,2,3,4,1,16,8,7)_{\adfLo}$,

$(0,1,2,6,7,11,9,3,33,23,22,\adfsplit 14,26,19,40,32,28,5,10,21,39,30)_{\adfLo}$,

\adfLgap
$(18,31,33,15,13,17,12,36,14,30,42,\adfsplit 23,5,16,20,1,0,3,2,10,38,27)_{\adfLp}$,

$(0,5,7,9,10,12,21,2,31,11,37,\adfsplit 39,33,16,1,14,8,4,27,18,3,36)_{\adfLp}$,

\adfLgap
$(1,16,18,4,17,14,12,41,22,23,39,\adfsplit 29,11,35,30,0,3,2,5,32,8,7)_{\adfLq}$,

$(0,2,9,10,5,7,16,40,28,41,18,\adfsplit 17,31,15,6,3,11,1,14,33,38,24)_{\adfLq}$,

\adfLgap
$(17,15,31,39,8,29,2,28,42,33,22,\adfsplit 14,3,23,24,0,5,4,6,1,9,26)_{\adfLr}$,

$(0,1,5,6,2,3,32,12,20,17,14,\adfsplit 21,23,11,33,27,40,15,28,38,26,31)_{\adfLr}$,

\adfLgap
$(24,35,27,5,6,2,0,36,7,18,21,\adfsplit 23,14,38,8,4,9,1,11,31,13,16)_{\adfLs}$,

$(0,1,5,7,2,6,11,16,25,13,43,\adfsplit 8,26,37,34,28,27,42,9,3,22,40)_{\adfLs}$

\noindent and

$(27,37,2,16,31,40,24,15,14,6,41,\adfsplit 35,12,10,25,0,3,5,1,4,18,26)_{\adfLt}$,

$(0,2,5,9,1,3,14,12,28,7,37,\adfsplit 29,39,35,38,22,19,16,24,41,11,10)_{\adfLt}$

\noindent under the action of the mapping $x \mapsto x + 4$ (mod 44).

\adfVfy{44, \{\{44,11,4\}\}, -1, \{\{11,\{0,1,2,3\}\}\}, -1} 

Let the vertex set of $K_{22,22,22,55}$ be $\{0, 1, \dots, 120\}$ partitioned into
$\{3j + i: j = 0, 1, \dots, 21\}$, $i = 0, 1, 2$, and $\{66, 67, \dots, 120\}$.
The decompositions consist of

\adfLgap
$(98,14,32,103,53,66,46,112,63,52,83,\adfsplit 44,18,97,35,39,43,101,62,102,96,0)_{\adfLa}$,

$(32,0,72,16,59,61,50,18,14,109,49,\adfsplit 106,37,9,116,6,21,34,66,95,52,2)_{\adfLa}$,

$(107,26,37,42,86,17,85,36,81,7,3,\adfsplit 11,40,32,4,1,15,90,21,104,109,43)_{\adfLa}$,

$(60,29,19,42,91,41,38,64,95,59,112,\adfsplit 58,83,55,63,30,85,46,98,50,24,104)_{\adfLa}$,

$(30,83,104,29,119,42,19,23,116,37,100,\adfsplit 63,85,2,33,25,84,4,60,20,16,118)_{\adfLa}$,

$(105,10,4,92,45,104,9,111,67,19,93,\adfsplit 32,3,103,39,56,41,109,49,72,79,12)_{\adfLa}$,

$(25,38,92,27,85,22,5,32,41,80,102,\adfsplit 58,99,56,75,9,50,33,108,78,24,16)_{\adfLa}$,

\adfLgap
$(107,58,43,24,112,30,95,62,25,53,59,\adfsplit 116,37,115,109,46,72,13,3,48,74,23)_{\adfLb}$,

$(91,50,39,88,9,34,44,49,56,18,84,\adfsplit 1,63,104,71,59,0,20,13,112,32,40)_{\adfLb}$,

$(26,39,100,75,27,50,20,10,41,105,30,\adfsplit 93,21,13,94,92,115,0,46,57,78,35)_{\adfLb}$,

$(96,49,11,21,101,10,72,31,85,60,52,\adfsplit 23,9,26,65,98,84,30,61,20,111,24)_{\adfLb}$,

$(71,61,21,90,1,73,58,5,19,2,100,\adfsplit 12,98,42,84,6,110,22,113,25,44,97)_{\adfLb}$,

$(2,57,106,71,8,43,18,96,17,87,86,\adfsplit 24,32,33,22,53,110,6,77,10,59,114)_{\adfLb}$,

$(2,16,28,102,13,83,38,12,17,114,109,\adfsplit 36,104,25,73,27,87,6,78,56,22,68)_{\adfLb}$,

\adfLgap
$(30,58,7,78,48,9,112,113,44,17,39,\adfsplit 62,94,105,19,85,4,12,51,87,16,99)_{\adfLc}$,

$(57,92,2,38,41,45,86,43,46,78,51,\adfsplit 68,20,1,16,36,100,87,80,21,61,11)_{\adfLc}$,

$(96,43,39,29,74,91,44,17,31,16,93,\adfsplit 113,60,24,12,6,87,46,80,41,106,14)_{\adfLc}$,

$(1,84,32,101,43,40,107,93,14,29,58,\adfsplit 36,111,103,81,2,54,48,33,47,74,5)_{\adfLc}$,

$(100,63,8,64,72,49,46,115,9,86,85,\adfsplit 61,109,30,7,17,44,104,114,90,37,18)_{\adfLc}$,

$(84,11,48,59,63,115,96,38,10,43,58,\adfsplit 18,110,108,98,67,6,42,44,7,17,77)_{\adfLc}$,

$(53,3,1,114,56,69,72,5,37,15,26,\adfsplit 6,76,86,88,92,38,55,43,61,62,36)_{\adfLc}$,

\adfLgap
$(110,8,3,9,27,81,41,40,69,25,48,\adfsplit 71,62,89,32,44,61,31,105,78,51,18)_{\adfLd}$,

$(88,11,13,59,28,80,93,65,103,19,37,\adfsplit 116,7,62,86,27,102,83,31,6,17,63)_{\adfLd}$,

$(105,63,34,37,5,81,111,27,112,41,33,\adfsplit 113,18,70,59,42,71,4,80,49,89,12)_{\adfLd}$,

$(17,4,102,107,48,21,35,20,11,28,49,\adfsplit 40,82,80,66,0,43,8,55,5,114,104)_{\adfLd}$,

$(86,57,8,32,58,117,106,52,82,27,63,\adfsplit 87,4,75,67,18,35,28,50,15,119,99)_{\adfLd}$,

$(58,60,87,80,74,92,48,53,5,34,3,\adfsplit 12,73,83,98,47,32,4,28,54,104,9)_{\adfLd}$,

$(41,24,45,108,85,84,10,4,63,2,17,\adfsplit 81,65,105,89,9,1,61,19,35,88,56)_{\adfLd}$,

\adfLgap
$(110,23,21,52,95,16,2,113,15,53,81,\adfsplit 4,1,106,64,105,101,57,35,79,25,13)_{\adfLe}$,

$(57,37,80,50,51,75,17,53,119,110,18,\adfsplit 33,42,40,114,71,102,9,20,44,45,4)_{\adfLe}$,

$(34,115,93,29,48,0,65,53,97,71,20,\adfsplit 63,108,14,1,98,45,24,103,105,50,49)_{\adfLe}$,

$(89,48,13,27,55,72,17,56,112,76,11,\adfsplit 57,64,53,105,113,91,88,7,26,63,45)_{\adfLe}$,

$(75,1,48,5,79,26,62,87,33,19,106,\adfsplit 66,17,116,14,41,9,36,49,40,107,74)_{\adfLe}$,

$(101,52,53,33,97,20,12,7,25,11,89,\adfsplit 94,17,74,43,24,108,5,30,98,34,92)_{\adfLe}$,

$(30,13,99,46,75,112,25,29,27,84,58,\adfsplit 48,113,28,67,44,42,92,10,4,68,21)_{\adfLe}$,

\adfLgap
$(41,6,111,18,81,112,22,51,20,58,26,\adfsplit 108,71,45,25,24,9,61,13,107,53,97)_{\adfLf}$,

$(8,97,68,106,59,35,21,61,9,12,89,\adfsplit 83,78,80,52,2,50,0,19,29,100,105)_{\adfLf}$,

$(69,35,37,57,117,45,8,44,16,94,34,\adfsplit 54,78,21,36,118,52,26,110,40,71,1)_{\adfLf}$,

$(91,38,37,23,49,108,72,47,6,74,21,\adfsplit 0,90,64,95,8,52,30,86,62,15,13)_{\adfLf}$,

$(91,8,20,13,55,88,114,64,27,118,4,\adfsplit 32,90,66,53,77,6,41,2,36,103,54)_{\adfLf}$,

$(5,60,103,12,10,69,63,64,84,58,50,\adfsplit 112,41,8,116,29,114,28,48,102,6,110)_{\adfLf}$,

$(18,107,40,25,16,0,109,99,70,60,100,\adfsplit 4,26,8,6,111,22,80,89,28,77,23)_{\adfLf}$,

\adfLgap
$(22,62,48,44,89,91,2,35,45,71,12,\adfsplit 19,92,98,47,77,31,51,17,32,39,49)_{\adfLg}$,

$(40,51,113,95,84,61,57,9,0,13,78,\adfsplit 23,92,107,70,45,26,65,22,64,15,75)_{\adfLg}$,

$(41,70,91,106,40,20,5,25,59,61,89,\adfsplit 79,78,52,60,58,11,64,6,109,13,103)_{\adfLg}$,

$(24,72,26,104,5,15,115,78,13,56,43,\adfsplit 53,62,112,107,83,99,88,49,8,19,45)_{\adfLg}$,

$(69,62,30,53,6,73,89,52,114,117,5,\adfsplit 55,110,0,48,42,12,24,50,88,7,101)_{\adfLg}$,

$(79,10,44,20,96,86,19,34,3,92,6,\adfsplit 62,101,65,16,106,37,28,115,51,95,57)_{\adfLg}$,

$(100,1,9,35,17,73,96,44,60,105,27,\adfsplit 24,85,10,58,112,3,21,106,47,102,52)_{\adfLg}$,

\adfLgap
$(11,63,37,79,85,87,39,3,14,43,8,\adfsplit 80,74,19,110,31,26,96,40,38,18,48)_{\adfLh}$,

$(36,104,38,67,15,9,109,113,52,0,119,\adfsplit 59,24,92,74,48,11,14,61,16,75,76)_{\adfLh}$,

$(100,8,19,57,33,24,78,76,112,95,94,\adfsplit 25,18,34,59,26,96,54,55,13,9,118)_{\adfLh}$,

$(119,0,25,50,47,1,108,115,37,109,96,\adfsplit 8,49,112,31,17,114,6,59,2,61,107)_{\adfLh}$,

$(17,48,70,60,118,22,36,15,32,101,33,\adfsplit 64,41,76,7,28,102,20,62,77,11,103)_{\adfLh}$,

$(94,1,46,21,15,35,93,3,26,115,107,\adfsplit 52,66,27,12,40,54,108,68,34,88,5)_{\adfLh}$,

$(17,25,86,64,116,120,37,29,45,67,56,\adfsplit 27,73,23,50,31,19,75,101,60,77,32)_{\adfLh}$,

\adfLgap
$(0,91,58,112,53,18,114,84,7,8,113,\adfsplit 12,55,74,117,2,110,27,38,56,54,70)_{\adfLi}$,

$(115,10,19,61,82,116,9,3,63,95,7,\adfsplit 35,8,101,41,24,78,106,25,17,33,65)_{\adfLi}$,

$(6,69,83,13,24,1,35,48,77,75,118,\adfsplit 43,106,3,45,11,36,5,49,84,14,103)_{\adfLi}$,

$(108,7,0,46,86,118,34,55,44,114,45,\adfsplit 38,65,72,49,83,81,42,6,29,9,43)_{\adfLi}$,

$(93,11,53,47,28,63,107,85,80,73,90,\adfsplit 36,32,9,48,29,61,116,104,13,99,57)_{\adfLi}$,

$(98,6,5,16,87,70,40,15,11,105,19,\adfsplit 54,58,55,50,117,66,111,44,14,7,28)_{\adfLi}$,

$(5,4,109,58,112,97,3,8,11,63,1,\adfsplit 91,69,12,119,45,23,32,37,117,55,70)_{\adfLi}$,

\adfLgap
$(68,54,1,39,64,78,8,5,113,88,17,\adfsplit 30,102,50,43,75,55,40,114,84,15,3)_{\adfLj}$,

$(95,56,16,51,112,109,47,38,17,58,54,\adfsplit 83,91,60,82,35,50,53,97,7,39,40)_{\adfLj}$,

$(7,81,82,106,39,27,56,35,25,48,96,\adfsplit 93,37,109,40,61,107,66,15,41,42,50)_{\adfLj}$,

$(62,10,9,116,70,105,25,40,16,30,8,\adfsplit 92,99,74,94,15,22,12,21,32,3,81)_{\adfLj}$,

$(28,5,103,107,108,46,49,57,26,35,76,\adfsplit 106,90,105,33,52,25,30,17,44,100,68)_{\adfLj}$,

$(3,92,113,109,9,40,14,56,59,10,87,\adfsplit 74,105,66,15,52,51,12,47,53,85,84)_{\adfLj}$,

$(29,57,120,61,71,47,18,28,32,36,114,\adfsplit 40,56,31,19,60,7,92,75,54,118,20)_{\adfLj}$,

\adfLgap
$(64,105,90,59,44,43,51,52,16,108,93,\adfsplit 92,20,53,1,42,85,99,112,45,56,40)_{\adfLk}$,

$(9,74,118,109,17,22,12,29,31,44,39,\adfsplit 27,115,76,35,86,42,65,19,36,62,67)_{\adfLk}$,

$(45,10,113,108,67,114,16,32,35,24,19,\adfsplit 38,48,97,0,105,71,41,17,69,54,58)_{\adfLk}$,

$(12,115,22,86,28,39,71,97,8,35,43,\adfsplit 16,7,101,65,27,114,30,74,70,25,31)_{\adfLk}$,

$(11,80,42,46,52,30,88,4,73,53,108,\adfsplit 43,45,87,50,6,101,9,72,91,58,32)_{\adfLk}$,

$(61,53,45,112,82,74,26,103,27,48,46,\adfsplit 41,2,99,93,94,4,60,51,32,30,101)_{\adfLk}$,

$(51,85,95,31,17,29,43,64,118,116,119,\adfsplit 19,53,38,9,56,115,96,117,16,48,21)_{\adfLk}$,

\adfLgap
$(65,42,99,12,111,116,35,22,29,78,56,\adfsplit 50,48,54,6,114,43,73,14,110,7,0)_{\adfLl}$,

$(106,61,63,46,107,118,34,17,62,105,25,\adfsplit 38,109,39,91,102,78,33,14,43,65,31)_{\adfLl}$,

$(109,25,12,43,69,112,7,110,8,63,36,\adfsplit 17,22,114,66,48,78,62,59,60,56,58)_{\adfLl}$,

$(67,53,7,65,98,99,35,101,1,88,16,\adfsplit 40,42,55,72,107,20,64,6,76,36,87)_{\adfLl}$,

$(40,109,113,81,19,48,8,16,0,3,91,\adfsplit 104,93,79,55,64,32,39,51,54,117,90)_{\adfLl}$,

$(111,32,52,59,101,0,63,29,34,105,108,\adfsplit 94,90,56,17,11,54,45,70,61,100,47)_{\adfLl}$,

$(88,19,12,6,27,33,82,92,75,5,62,\adfsplit 1,37,72,70,14,115,8,24,42,4,55)_{\adfLl}$,

\adfLgap
$(119,35,8,64,93,77,1,45,51,27,59,\adfsplit 42,81,85,104,103,2,44,25,46,21,82)_{\adfLm}$,

$(43,98,80,66,20,49,12,44,4,13,103,\adfsplit 73,101,60,53,46,30,19,75,110,2,10)_{\adfLm}$,

$(106,9,5,3,102,112,30,75,22,99,11,\adfsplit 38,21,89,58,118,41,17,2,72,13,49)_{\adfLm}$,

$(56,30,84,92,104,79,24,3,43,17,52,\adfsplit 36,41,14,4,77,16,94,98,65,44,9)_{\adfLm}$,

$(18,112,97,95,54,25,41,30,5,40,103,\adfsplit 11,68,6,109,86,1,32,108,22,51,36)_{\adfLm}$,

$(72,28,42,54,84,30,37,81,13,110,14,\adfsplit 119,50,111,38,32,80,39,55,100,91,33)_{\adfLm}$,

$(92,14,34,22,37,12,45,110,93,70,55,\adfsplit 96,58,20,33,1,66,30,106,65,35,21)_{\adfLm}$,

\adfLgap
$(107,23,6,12,55,118,98,95,96,61,26,\adfsplit 33,8,60,78,25,10,84,113,45,39,65)_{\adfLn}$,

$(22,73,83,91,8,59,24,4,0,49,110,\adfsplit 90,70,19,107,35,28,33,21,58,32,2)_{\adfLn}$,

$(65,118,112,68,38,57,39,31,63,7,80,\adfsplit 84,72,61,106,49,46,59,21,2,8,99)_{\adfLn}$,

$(103,32,39,5,30,111,34,79,102,45,19,\adfsplit 40,64,18,109,93,80,89,24,2,46,9)_{\adfLn}$,

$(21,34,90,84,41,26,10,64,45,11,69,\adfsplit 105,81,62,55,0,53,37,38,67,82,6)_{\adfLn}$,

$(12,22,82,99,87,103,28,35,19,23,45,\adfsplit 4,67,81,40,77,56,24,53,33,101,66)_{\adfLn}$,

$(55,26,89,79,66,96,54,3,20,59,32,\adfsplit 16,76,110,43,80,0,18,45,56,23,75)_{\adfLn}$,

\adfLgap
$(8,101,36,93,56,17,13,111,4,23,39,\adfsplit 85,11,115,49,89,46,97,5,110,30,40)_{\adfLo}$,

$(12,100,102,103,8,52,15,44,46,35,86,\adfsplit 82,37,26,43,0,14,95,78,54,19,2)_{\adfLo}$,

$(13,83,94,38,36,42,1,37,72,116,77,\adfsplit 81,12,35,23,47,7,29,30,103,67,0)_{\adfLo}$,

$(36,28,105,102,68,72,17,54,39,27,22,\adfsplit 65,87,49,92,85,73,33,106,35,47,4)_{\adfLo}$,

$(71,13,27,58,78,86,55,83,5,59,40,\adfsplit 30,19,68,70,42,113,15,1,54,50,69)_{\adfLo}$,

$(58,82,81,36,20,55,22,28,79,110,99,\adfsplit 5,45,19,60,50,84,32,89,46,15,119)_{\adfLo}$,

$(51,17,91,35,84,57,9,45,1,99,52,\adfsplit 80,76,21,63,5,12,11,22,59,75,69)_{\adfLo}$,

\adfLgap
$(38,33,103,37,102,52,23,35,27,67,79,\adfsplit 70,42,25,41,11,58,26,54,10,111,117)_{\adfLp}$,

$(90,54,15,40,66,119,58,102,38,94,61,\adfsplit 11,51,111,64,100,57,105,53,2,25,42)_{\adfLp}$,

$(6,23,113,116,77,21,1,51,42,12,75,\adfsplit 115,101,50,53,34,26,2,64,110,93,8)_{\adfLp}$,

$(61,48,66,2,55,26,60,35,39,111,58,\adfsplit 110,95,64,45,47,24,49,101,113,89,22)_{\adfLp}$,

$(72,14,12,62,61,43,79,108,82,13,21,\adfsplit 23,29,63,66,22,28,110,50,4,36,102)_{\adfLp}$,

$(61,80,76,72,45,62,26,22,25,54,104,\adfsplit 84,114,63,53,5,41,27,94,119,88,57)_{\adfLp}$,

$(27,1,17,103,98,73,60,113,21,25,2,\adfsplit 4,37,77,88,72,116,99,16,39,3,8)_{\adfLp}$,

\adfLgap
$(36,100,10,74,2,33,113,95,22,41,8,\adfsplit 117,60,106,64,92,59,58,110,34,50,21)_{\adfLq}$,

$(6,101,96,81,61,56,14,50,48,21,36,\adfsplit 71,63,66,49,79,26,102,53,9,38,57)_{\adfLq}$,

$(53,43,76,58,99,18,47,17,57,119,34,\adfsplit 91,117,66,52,65,15,41,54,25,55,79)_{\adfLq}$,

$(118,56,51,34,19,115,17,113,107,119,5,\adfsplit 49,16,32,42,73,3,108,98,58,33,35)_{\adfLq}$,

$(76,58,40,22,97,78,48,95,39,117,35,\adfsplit 11,64,4,42,70,119,20,69,8,77,54)_{\adfLq}$,

$(106,19,41,18,78,90,57,103,29,52,49,\adfsplit 10,15,105,117,30,119,37,58,54,44,80)_{\adfLq}$,

$(43,68,5,85,0,34,79,16,8,15,50,\adfsplit 117,75,51,20,97,49,63,17,94,93,45)_{\adfLq}$,

\adfLgap
$(72,43,7,3,73,53,57,87,41,55,75,\adfsplit 84,62,61,45,60,8,11,92,119,47,6)_{\adfLr}$,

$(33,106,89,102,65,27,4,38,51,15,64,\adfsplit 37,94,101,25,119,56,103,7,60,24,68)_{\adfLr}$,

$(40,8,75,41,115,85,28,25,37,6,50,\adfsplit 57,96,93,33,78,19,83,18,26,38,104)_{\adfLr}$,

$(6,115,52,80,48,34,8,71,55,25,75,\adfsplit 73,41,14,67,33,64,49,36,23,81,39)_{\adfLr}$,

$(28,113,117,82,47,19,18,12,15,29,0,\adfsplit 112,84,79,38,92,36,51,52,86,50,4)_{\adfLr}$,

$(111,49,28,8,77,104,17,112,9,93,7,\adfsplit 33,65,20,75,44,69,23,66,64,13,41)_{\adfLr}$,

$(9,70,81,52,11,27,4,63,91,94,85,\adfsplit 83,44,115,58,22,17,6,14,111,68,12)_{\adfLr}$,

\adfLgap
$(33,113,31,91,42,48,86,114,19,22,14,\adfsplit 119,37,46,67,110,11,25,12,23,15,84)_{\adfLs}$,

$(83,42,59,2,34,35,69,25,86,40,101,\adfsplit 75,24,110,28,65,10,36,7,108,13,23)_{\adfLs}$,

$(37,79,95,33,28,47,8,59,73,82,83,\adfsplit 96,52,106,12,23,19,26,108,104,34,56)_{\adfLs}$,

$(49,96,83,2,56,11,27,40,109,117,103,\adfsplit 112,24,21,7,65,43,47,72,115,39,44)_{\adfLs}$,

$(111,3,7,51,26,109,117,82,31,58,6,\adfsplit 18,1,24,14,85,86,107,33,27,37,11)_{\adfLs}$,

$(98,47,19,26,33,55,92,73,78,22,27,\adfsplit 102,16,60,25,90,14,0,110,66,53,42)_{\adfLs}$,

$(44,3,9,115,85,74,11,75,23,30,53,\adfsplit 0,57,8,24,71,89,40,94,97,18,35)_{\adfLs}$

\noindent and

$(65,15,19,51,77,98,12,118,37,110,57,\adfsplit 5,28,63,79,49,82,81,56,3,109,1)_{\adfLt}$,

$(56,49,82,88,47,111,0,59,6,46,34,\adfsplit 14,69,67,23,68,95,21,38,4,51,22)_{\adfLt}$,

$(0,72,117,95,7,10,50,25,38,45,59,\adfsplit 90,92,42,85,13,101,62,61,18,84,36)_{\adfLt}$,

$(3,19,93,70,68,83,5,60,61,4,1,\adfsplit 18,79,89,51,106,67,62,57,55,65,49)_{\adfLt}$,

$(77,41,58,51,89,99,53,98,28,81,32,\adfsplit 54,33,17,90,11,110,71,5,26,83,36)_{\adfLt}$,

$(64,0,113,42,31,81,35,63,71,75,61,\adfsplit 40,111,53,19,14,96,69,39,44,77,2)_{\adfLt}$,

$(44,58,87,42,95,90,63,41,8,37,36,\adfsplit 32,49,30,119,66,84,78,60,19,56,22)_{\adfLt}$

\noindent under the action of the mapping $x \mapsto x + 3$ (mod 66) for $x < 66$,
$x \mapsto 66 + (x - 66 + 5 \mathrm{~(mod~55)})$ for $x \ge 66$.

\adfVfy{121, \{\{66,22,3\},\{55,11,5\}\}, -1, \{\{22,\{0,1,2\}\},\{55,\{3\}\}\}, -1} 

Let the vertex set of $K_{22,22,22,22,22,22,55}$ be $\{0, 1, \dots, 186\}$ partitioned into
$\{6j + i: j = 0, 1, \dots, 21\}$, $i = 0, 1, \dots, 5$, and $\{132, 133, \dots, 186\}$.
The decompositions consist of

\adfLgap
$(181,64,78,15,138,94,150,126,132,16,28,\adfsplit 5,106,79,130,128,84,51,48,179,116,118)_{\adfLa}$,

$(178,27,80,79,82,73,87,48,10,125,135,\adfsplit 101,140,91,85,43,146,120,15,47,182,8)_{\adfLa}$,

$(92,173,165,91,20,111,80,73,14,33,22,\adfsplit 150,35,178,102,87,4,26,176,32,118,89)_{\adfLa}$,

$(142,59,111,104,113,173,134,123,14,117,105,\adfsplit 77,166,72,73,146,67,110,29,114,108,115)_{\adfLa}$,

$(180,75,41,155,29,156,28,72,21,83,167,\adfsplit 76,138,88,32,50,127,9,149,63,38,154)_{\adfLa}$,

$(73,26,155,48,145,46,65,6,140,27,121,\adfsplit 113,99,98,67,152,86,167,7,54,179,10)_{\adfLa}$,

$(18,31,38,50,143,65,35,10,109,169,162,\adfsplit 51,25,30,33,8,156,103,3,6,12,81)_{\adfLa}$,

$(87,151,109,127,113,139,149,117,157,108,106,\adfsplit 7,181,67,23,123,51,13,147,148,12,27)_{\adfLa}$,

$(1,184,52,20,161,26,171,91,57,46,2,\adfsplit 88,148,113,107,35,72,157,102,164,77,112)_{\adfLa}$,

$(60,150,2,67,92,22,132,153,20,134,47,\adfsplit 97,76,3,105,21,24,8,173,122,79,11)_{\adfLa}$,

\adfLgap
$(96,16,57,23,114,176,80,110,63,37,4,\adfsplit 185,21,66,150,135,84,71,28,74,20,123)_{\adfLb}$,

$(181,42,57,113,60,158,22,39,111,12,63,\adfsplit 16,159,52,162,68,122,103,126,54,73,119)_{\adfLb}$,

$(131,14,147,125,33,115,86,121,63,151,88,\adfsplit 17,94,79,60,181,133,5,32,20,126,185)_{\adfLb}$,

$(139,77,62,163,14,140,84,10,122,114,99,\adfsplit 164,7,8,156,183,43,81,80,21,63,76)_{\adfLb}$,

$(12,16,40,44,145,125,42,118,56,70,174,\adfsplit 76,107,63,144,88,102,97,78,71,27,152)_{\adfLb}$,

$(178,15,4,47,27,166,91,86,54,101,57,\adfsplit 143,41,1,70,126,164,21,150,50,123,132)_{\adfLb}$,

$(121,154,6,32,162,119,184,104,118,133,63,\adfsplit 56,179,34,158,138,66,109,33,67,24,8)_{\adfLb}$,

$(184,57,91,49,114,52,16,152,111,79,29,\adfsplit 132,58,116,83,128,30,166,59,167,180,37)_{\adfLb}$,

$(18,160,43,103,165,52,171,126,30,80,8,\adfsplit 5,170,85,182,17,62,100,114,163,147,64)_{\adfLb}$,

$(148,81,124,149,65,14,89,166,146,115,79,\adfsplit 56,85,35,31,16,92,176,117,67,97,137)_{\adfLb}$,

\adfLgap
$(11,87,13,148,61,48,143,71,32,88,68,\adfsplit 180,66,53,136,78,144,86,57,83,175,7)_{\adfLc}$,

$(32,111,118,18,79,149,169,57,158,64,21,\adfsplit 10,98,155,185,53,72,75,76,106,73,85)_{\adfLc}$,

$(152,75,78,128,22,112,109,137,5,157,116,\adfsplit 72,56,47,120,172,11,150,20,147,131,79)_{\adfLc}$,

$(48,13,184,123,156,74,65,106,50,182,163,\adfsplit 44,34,69,4,90,27,46,178,97,71,35)_{\adfLc}$,

$(46,146,142,87,15,76,30,11,104,7,134,\adfsplit 6,69,45,23,128,118,73,48,182,84,158)_{\adfLc}$,

$(162,100,43,79,81,95,16,98,101,163,68,\adfsplit 178,18,34,75,59,76,125,67,120,106,96)_{\adfLc}$,

$(126,3,89,49,162,5,88,151,64,134,137,\adfsplit 15,42,4,23,168,20,176,63,114,55,7)_{\adfLc}$,

$(159,109,72,36,14,66,180,171,146,8,52,\adfsplit 90,94,101,127,35,80,74,102,145,105,168)_{\adfLc}$,

$(159,24,29,107,75,163,42,30,76,39,180,\adfsplit 25,136,22,124,164,32,149,68,165,116,28)_{\adfLc}$,

$(109,88,58,161,50,145,146,42,3,8,87,\adfsplit 53,74,90,180,174,124,186,77,127,63,20)_{\adfLc}$,

\adfLgap
$(130,8,39,175,119,177,68,44,91,76,98,\adfsplit 75,11,102,58,47,168,5,144,9,70,1)_{\adfLd}$,

$(164,35,21,126,91,45,151,2,99,44,32,\adfsplit 85,90,183,94,163,128,14,51,66,17,142)_{\adfLd}$,

$(23,145,34,86,125,84,91,184,93,33,49,\adfsplit 48,124,155,2,150,138,9,111,31,130,72)_{\adfLd}$,

$(168,29,30,81,133,73,105,157,16,61,146,\adfsplit 46,135,90,10,148,33,137,126,163,95,58)_{\adfLd}$,

$(13,138,52,158,81,99,162,113,34,20,150,\adfsplit 39,116,19,23,154,36,62,3,161,94,37)_{\adfLd}$,

$(185,62,94,124,166,87,115,140,56,111,120,\adfsplit 131,136,159,145,170,91,63,113,128,89,82)_{\adfLd}$,

$(163,42,124,11,161,129,51,7,98,88,176,\adfsplit 34,85,90,6,156,80,10,121,117,120,134)_{\adfLd}$,

$(125,156,118,146,107,40,9,59,98,91,159,\adfsplit 74,184,5,94,131,17,43,152,169,26,105)_{\adfLd}$,

$(174,57,124,102,59,35,79,142,51,16,149,\adfsplit 45,110,107,7,85,145,49,152,0,71,81)_{\adfLd}$,

$(88,102,62,162,49,119,46,108,128,84,22,\adfsplit 31,29,147,161,167,120,60,45,105,20,44)_{\adfLd}$,

\adfLgap
$(138,105,50,67,169,8,99,45,65,35,82,\adfsplit 153,120,179,71,130,77,101,68,145,66,109)_{\adfLe}$,

$(124,138,37,127,33,111,84,24,168,157,126,\adfsplit 57,97,14,153,142,66,19,28,23,3,12)_{\adfLe}$,

$(166,121,31,129,140,70,117,47,46,97,185,\adfsplit 178,164,116,119,4,11,88,81,96,179,174)_{\adfLe}$,

$(26,183,47,125,128,96,141,135,73,117,173,\adfsplit 83,99,74,98,51,174,1,105,176,54,38)_{\adfLe}$,

$(44,1,131,148,65,166,80,25,128,59,109,\adfsplit 172,42,61,158,72,112,41,33,16,26,28)_{\adfLe}$,

$(183,35,121,109,142,51,60,65,80,82,155,\adfsplit 132,117,154,100,7,56,85,66,114,152,175)_{\adfLe}$,

$(137,112,55,29,2,156,162,78,150,120,81,\adfsplit 98,99,14,85,58,160,52,116,126,23,15)_{\adfLe}$,

$(70,134,164,127,67,121,74,112,34,105,96,\adfsplit 55,135,111,113,90,124,5,177,166,57,32)_{\adfLe}$,

$(11,84,180,28,29,112,104,114,48,43,167,\adfsplit 54,124,172,87,55,80,56,156,136,81,3)_{\adfLe}$,

$(81,175,83,53,126,55,122,184,144,117,41,\adfsplit 61,13,105,176,8,166,156,56,43,106,96)_{\adfLe}$,

\adfLgap
$(107,74,91,86,18,166,77,88,58,146,27,\adfsplit 73,120,110,29,104,171,184,111,93,37,12)_{\adfLf}$,

$(149,22,103,131,138,172,92,38,70,99,4,\adfsplit 28,112,184,86,114,95,181,117,49,164,118)_{\adfLf}$,

$(99,167,24,55,3,40,122,17,26,181,14,\adfsplit 73,185,90,19,106,0,38,172,160,69,85)_{\adfLf}$,

$(27,54,47,41,176,33,62,142,100,96,179,\adfsplit 165,61,185,122,67,76,66,114,97,119,138)_{\adfLf}$,

$(142,109,116,95,183,126,11,167,12,61,74,\adfsplit 152,97,119,65,19,28,135,111,118,75,46)_{\adfLf}$,

$(103,113,52,105,63,102,178,93,127,8,174,\adfsplit 115,16,14,73,45,157,145,7,65,0,64)_{\adfLf}$,

$(50,173,10,113,110,85,19,47,160,182,83,\adfsplit 33,171,22,102,17,94,126,15,176,59,148)_{\adfLf}$,

$(115,180,81,76,110,108,23,168,144,24,34,\adfsplit 156,17,125,145,96,48,54,99,172,45,70)_{\adfLf}$,

$(62,57,139,118,42,175,83,87,129,86,125,\adfsplit 158,168,41,54,31,94,63,173,140,50,101)_{\adfLf}$,

$(158,113,25,105,167,1,45,174,26,12,154,\adfsplit 37,33,166,110,24,144,95,14,66,31,186)_{\adfLf}$,

\adfLgap
$(90,74,185,119,54,3,62,53,182,149,37,\adfsplit 51,20,118,56,168,184,108,102,164,15,44)_{\adfLg}$,

$(80,140,138,67,108,28,7,50,14,44,141,\adfsplit 112,157,51,121,123,95,16,13,131,96,77)_{\adfLg}$,

$(44,84,25,100,183,112,28,14,42,155,63,\adfsplit 133,21,106,20,85,163,171,77,55,69,81)_{\adfLg}$,

$(140,72,73,37,164,19,106,125,84,176,59,\adfsplit 98,142,139,118,177,42,58,44,62,167,87)_{\adfLg}$,

$(119,174,18,46,100,77,99,144,168,111,91,\adfsplit 181,48,38,70,39,67,133,23,93,24,109)_{\adfLg}$,

$(13,82,130,8,126,174,134,23,1,152,27,\adfsplit 79,80,57,124,89,184,22,67,48,103,180)_{\adfLg}$,

$(105,92,48,100,136,41,65,80,38,147,68,\adfsplit 178,13,170,82,102,180,110,108,23,7,63)_{\adfLg}$,

$(52,69,153,136,156,46,107,78,81,22,117,\adfsplit 126,137,54,168,165,66,53,51,94,56,30)_{\adfLg}$,

$(57,110,70,55,145,165,129,98,92,152,49,\adfsplit 94,184,61,12,47,115,112,156,29,26,120)_{\adfLg}$,

$(158,105,23,47,77,132,14,45,176,186,44,\adfsplit 35,181,64,60,111,25,104,160,103,147,78)_{\adfLg}$,

\adfLgap
$(31,176,158,183,52,41,72,79,30,18,76,\adfsplit 116,51,97,169,138,47,135,45,113,62,59)_{\adfLh}$,

$(153,65,77,92,22,37,147,114,75,139,94,\adfsplit 107,134,76,43,113,41,69,6,171,30,10)_{\adfLh}$,

$(84,103,148,158,157,137,1,112,99,122,50,\adfsplit 151,119,120,140,73,169,33,19,104,65,2)_{\adfLh}$,

$(4,43,12,167,107,28,104,154,118,49,78,\adfsplit 130,59,32,141,56,165,178,126,112,69,99)_{\adfLh}$,

$(82,61,55,48,172,87,78,14,116,129,49,\adfsplit 2,109,123,139,93,184,20,75,39,64,157)_{\adfLh}$,

$(91,165,183,155,105,0,118,8,23,84,107,\adfsplit 166,66,14,115,76,174,71,52,30,79,160)_{\adfLh}$,

$(119,40,121,181,56,20,141,164,30,31,134,\adfsplit 87,103,154,113,83,125,107,150,123,7,34)_{\adfLh}$,

$(122,175,168,97,60,3,100,50,185,106,72,\adfsplit 110,129,54,41,26,17,93,128,108,142,58)_{\adfLh}$,

$(152,41,102,54,74,18,146,176,118,157,16,\adfsplit 71,20,19,48,6,75,177,160,110,178,1)_{\adfLh}$,

$(165,94,106,80,24,20,183,48,141,85,152,\adfsplit 126,57,118,5,75,40,31,74,136,72,42)_{\adfLh}$,

\adfLgap
$(183,47,65,11,105,139,78,120,32,52,59,\adfsplit 163,64,39,3,50,41,109,90,96,76,153)_{\adfLi}$,

$(80,48,124,95,147,85,29,125,79,161,46,\adfsplit 31,27,123,68,149,24,77,53,139,45,9)_{\adfLi}$,

$(45,144,148,98,36,61,118,84,71,126,60,\adfsplit 50,73,30,13,159,138,87,123,1,35,62)_{\adfLi}$,

$(48,13,170,185,51,95,43,34,47,99,44,\adfsplit 163,18,156,62,175,181,171,86,24,82,40)_{\adfLi}$,

$(18,109,17,174,178,20,84,180,93,85,39,\adfsplit 182,36,124,118,146,160,151,120,1,59,43)_{\adfLi}$,

$(134,131,79,25,181,182,12,4,110,88,29,\adfsplit 45,115,158,54,37,159,50,55,103,2,167)_{\adfLi}$,

$(143,57,69,95,120,161,80,54,137,106,10,\adfsplit 145,175,61,174,86,121,59,89,88,8,147)_{\adfLi}$,

$(8,91,11,126,183,86,50,73,3,58,159,\adfsplit 30,167,55,111,17,82,118,87,33,170,104)_{\adfLi}$,

$(138,113,55,16,111,141,65,59,185,150,67,\adfsplit 122,180,30,52,42,19,23,100,157,12,181)_{\adfLi}$,

$(125,120,31,24,147,167,14,34,69,132,97,\adfsplit 36,8,43,23,122,1,71,174,45,176,114)_{\adfLi}$,

\adfLgap
$(179,71,32,46,31,114,177,100,111,105,119,\adfsplit 49,21,121,76,88,4,126,142,118,2,50)_{\adfLj}$,

$(71,54,70,141,55,149,144,98,33,68,105,\adfsplit 130,10,0,142,164,35,23,159,49,128,51)_{\adfLj}$,

$(4,133,85,90,59,11,56,66,67,145,21,\adfsplit 91,161,102,29,36,51,95,142,160,7,45)_{\adfLj}$,

$(33,120,30,11,116,113,160,165,1,155,40,\adfsplit 21,79,53,114,59,106,176,24,102,133,52)_{\adfLj}$,

$(132,83,60,84,49,172,124,125,98,177,77,\adfsplit 103,39,163,30,175,73,158,118,69,167,123)_{\adfLj}$,

$(55,92,164,122,24,174,26,112,63,113,71,\adfsplit 87,98,76,148,16,111,35,166,72,32,31)_{\adfLj}$,

$(31,42,57,123,85,141,183,10,54,180,60,\adfsplit 49,101,75,71,26,166,20,110,119,14,152)_{\adfLj}$,

$(4,37,178,136,145,182,43,124,49,54,76,\adfsplit 156,174,58,13,65,1,83,161,60,168,115)_{\adfLj}$,

$(3,10,73,179,140,139,93,68,11,84,42,\adfsplit 109,176,145,104,23,113,131,33,55,170,70)_{\adfLj}$,

$(25,5,183,132,105,50,14,71,67,87,156,\adfsplit 29,96,99,12,60,80,133,164,98,153,111)_{\adfLj}$,

\adfLgap
$(58,152,168,79,69,100,80,90,142,72,49,\adfsplit 146,29,8,134,115,85,60,102,7,101,112)_{\adfLk}$,

$(122,184,0,37,8,64,47,106,185,134,9,\adfsplit 65,102,25,104,44,57,139,158,71,46,115)_{\adfLk}$,

$(179,117,15,99,17,73,142,141,157,32,0,\adfsplit 81,84,67,2,103,104,128,57,185,113,160)_{\adfLk}$,

$(57,34,19,131,77,101,81,38,157,114,120,\adfsplit 97,113,133,149,2,118,16,93,98,177,148)_{\adfLk}$,

$(145,40,84,28,85,96,174,172,175,99,111,\adfsplit 176,43,113,89,131,147,7,166,36,19,75)_{\adfLk}$,

$(60,87,129,136,145,110,126,25,99,131,47,\adfsplit 5,177,93,183,21,97,130,104,105,80,53)_{\adfLk}$,

$(144,105,44,102,54,48,133,51,93,121,77,\adfsplit 181,17,82,90,56,91,89,88,163,96,146)_{\adfLk}$,

$(78,182,141,57,9,112,43,7,58,10,168,\adfsplit 156,128,38,84,23,178,169,64,41,121,73)_{\adfLk}$,

$(112,171,53,127,90,6,133,168,130,27,68,\adfsplit 88,87,154,144,185,122,119,104,70,13,18)_{\adfLk}$,

$(141,82,89,67,29,14,147,2,140,4,131,\adfsplit 100,3,69,185,71,135,145,40,92,79,24)_{\adfLk}$,

\adfLgap
$(131,150,176,130,77,6,68,81,21,153,151,\adfsplit 49,185,100,122,169,161,99,0,91,120,118)_{\adfLl}$,

$(107,36,111,135,29,101,42,152,39,130,62,\adfsplit 87,124,53,142,94,116,110,86,72,14,166)_{\adfLl}$,

$(21,100,137,171,30,77,103,14,33,11,31,\adfsplit 134,114,124,184,72,40,93,89,139,87,173)_{\adfLl}$,

$(17,167,64,88,23,118,44,164,51,33,50,\adfsplit 25,71,163,95,83,113,133,109,12,140,34)_{\adfLl}$,

$(139,88,123,20,105,1,49,24,163,82,32,\adfsplit 36,77,155,61,116,147,72,73,97,65,100)_{\adfLl}$,

$(32,34,22,82,138,135,89,161,66,13,112,\adfsplit 83,44,51,0,99,67,176,8,178,31,72)_{\adfLl}$,

$(68,137,79,154,52,89,86,64,103,71,54,\adfsplit 51,12,174,67,159,85,164,119,39,116,60)_{\adfLl}$,

$(53,154,153,78,10,121,120,63,104,38,141,\adfsplit 118,173,178,54,132,14,4,6,60,15,180)_{\adfLl}$,

$(178,88,37,77,181,185,63,74,82,64,101,\adfsplit 128,141,175,54,71,65,20,112,68,46,142)_{\adfLl}$,

$(5,150,6,64,63,3,147,137,160,85,7,\adfsplit 52,15,24,176,54,12,82,129,107,45,167)_{\adfLl}$,

\adfLgap
$(127,98,77,81,137,31,78,68,90,29,16,\adfsplit 17,61,155,51,172,84,66,8,119,5,100)_{\adfLm}$,

$(54,175,163,178,45,84,52,123,100,113,124,\adfsplit 179,170,120,115,31,3,106,148,144,47,104)_{\adfLm}$,

$(100,102,144,162,172,163,14,4,94,57,33,\adfsplit 68,136,56,160,149,60,124,64,49,77,23)_{\adfLm}$,

$(8,182,136,167,4,10,109,113,111,21,158,\adfsplit 91,36,61,156,108,83,81,140,62,130,34)_{\adfLm}$,

$(30,99,38,175,64,71,160,180,90,37,70,\adfsplit 67,35,14,141,125,20,151,153,21,108,102)_{\adfLm}$,

$(142,103,74,5,178,58,81,4,112,169,45,\adfsplit 132,127,77,41,51,159,27,40,164,71,128)_{\adfLm}$,

$(184,51,63,2,56,41,30,47,150,161,156,\adfsplit 122,104,99,53,153,12,93,86,18,103,76)_{\adfLm}$,

$(144,121,68,102,48,27,79,140,87,46,86,\adfsplit 150,76,106,14,44,131,122,142,108,11,3)_{\adfLm}$,

$(169,126,40,61,59,154,163,87,147,27,90,\adfsplit 70,142,111,113,1,94,34,31,161,96,26)_{\adfLm}$,

$(94,51,15,71,181,40,76,17,30,110,173,\adfsplit 166,127,22,48,13,124,171,101,77,27,168)_{\adfLm}$,

\adfLgap
$(8,46,4,168,184,113,126,157,3,68,58,\adfsplit 48,130,61,134,89,64,149,104,84,123,5)_{\adfLn}$,

$(151,76,16,79,154,103,57,125,99,114,116,\adfsplit 11,96,142,143,87,3,34,135,47,67,146)_{\adfLn}$,

$(169,129,20,15,155,124,71,181,91,24,69,\adfsplit 49,21,30,110,159,152,146,92,85,109,108)_{\adfLn}$,

$(100,11,71,20,79,150,98,137,160,159,16,\adfsplit 81,60,62,101,114,103,39,161,140,21,84)_{\adfLn}$,

$(108,147,7,182,104,116,5,178,69,16,4,\adfsplit 78,23,84,31,177,110,124,20,96,3,41)_{\adfLn}$,

$(72,185,140,50,16,96,7,25,11,173,69,\adfsplit 109,101,43,97,120,105,129,152,113,131,56)_{\adfLn}$,

$(144,47,105,130,81,48,148,10,83,77,49,\adfsplit 73,101,153,57,157,94,114,80,119,43,61)_{\adfLn}$,

$(170,50,37,46,95,66,161,179,181,68,100,\adfsplit 150,40,56,180,62,19,87,163,55,116,5)_{\adfLn}$,

$(71,161,43,174,20,13,137,147,6,17,30,\adfsplit 96,90,143,166,53,92,93,183,61,31,148)_{\adfLn}$,

$(149,4,131,111,6,42,87,171,101,79,105,\adfsplit 110,46,166,158,135,99,1,116,27,106,70)_{\adfLn}$,

\adfLgap
$(87,185,136,146,57,104,92,121,102,12,147,\adfsplit 178,114,19,125,47,55,123,180,48,90,61)_{\adfLo}$,

$(80,105,21,156,119,2,72,88,97,39,184,\adfsplit 169,86,179,28,84,36,134,106,71,89,55)_{\adfLo}$,

$(97,160,166,146,48,10,74,65,35,11,81,\adfsplit 55,157,79,57,106,153,26,76,154,101,102)_{\adfLo}$,

$(9,24,157,91,93,23,7,88,44,14,132,\adfsplit 68,154,127,161,2,174,42,15,115,125,39)_{\adfLo}$,

$(37,110,51,60,31,163,173,52,50,32,56,\adfsplit 27,69,29,150,41,166,28,184,92,102,160)_{\adfLo}$,

$(122,17,109,147,110,84,99,157,90,25,57,\adfsplit 162,58,143,9,130,96,128,93,142,183,3)_{\adfLo}$,

$(110,49,72,139,34,107,132,161,103,20,56,\adfsplit 164,65,8,112,157,12,97,96,62,118,66)_{\adfLo}$,

$(62,52,126,29,26,132,155,51,49,85,125,\adfsplit 65,130,136,135,75,180,96,44,10,36,178)_{\adfLo}$,

$(183,10,21,116,102,114,64,169,119,133,106,\adfsplit 138,63,173,25,98,40,115,26,175,67,96)_{\adfLo}$,

$(136,40,103,36,91,92,112,131,124,57,30,\adfsplit 69,85,170,176,165,110,46,8,19,7,113)_{\adfLo}$,

\adfLgap
$(102,141,119,73,37,92,99,51,172,163,154,\adfsplit 160,161,97,80,7,55,4,82,59,120,56)_{\adfLp}$,

$(42,182,5,136,56,125,70,172,91,53,8,\adfsplit 102,68,112,169,146,167,7,128,76,82,24)_{\adfLp}$,

$(185,66,114,117,52,146,128,88,157,36,83,\adfsplit 11,1,33,40,135,165,99,107,131,42,126)_{\adfLp}$,

$(79,178,140,18,40,19,104,42,148,162,21,\adfsplit 75,11,33,83,157,120,74,133,173,119,69)_{\adfLp}$,

$(110,91,166,141,5,154,37,9,48,119,74,\adfsplit 55,156,85,175,41,94,124,107,25,135,30)_{\adfLp}$,

$(23,164,61,74,84,86,15,101,45,16,146,\adfsplit 30,60,159,33,70,92,56,131,18,19,165)_{\adfLp}$,

$(172,47,34,22,74,180,24,38,72,14,10,\adfsplit 7,158,119,60,69,174,25,76,132,35,96)_{\adfLp}$,

$(165,96,73,55,47,63,133,24,12,44,151,\adfsplit 153,180,131,53,48,127,9,130,22,169,107)_{\adfLp}$,

$(151,70,107,112,133,93,100,37,105,152,162,\adfsplit 115,84,159,113,108,24,7,13,144,16,92)_{\adfLp}$,

$(42,4,67,113,134,108,9,92,56,158,117,\adfsplit 64,99,138,116,65,0,11,97,90,31,144)_{\adfLp}$,

\adfLgap
$(59,171,180,36,10,94,19,129,31,86,110,\adfsplit 43,82,18,120,57,8,152,83,138,4,6)_{\adfLq}$,

$(70,174,14,73,69,74,152,136,66,33,114,\adfsplit 84,85,23,149,128,143,125,6,64,176,101)_{\adfLq}$,

$(135,86,123,126,52,11,48,139,33,149,58,\adfsplit 50,71,103,76,156,94,121,142,57,145,49)_{\adfLq}$,

$(109,93,64,47,152,115,101,21,114,143,0,\adfsplit 41,54,141,12,11,121,46,133,98,66,44)_{\adfLq}$,

$(185,43,110,52,4,175,154,84,23,42,118,\adfsplit 74,155,26,91,129,164,115,27,54,77,176)_{\adfLq}$,

$(46,92,158,133,159,83,9,58,55,40,172,\adfsplit 48,80,105,167,5,107,155,100,143,79,162)_{\adfLq}$,

$(68,149,27,148,89,47,2,65,105,51,22,\adfsplit 100,52,5,95,90,158,0,80,82,1,176)_{\adfLq}$,

$(143,55,29,76,88,126,60,153,17,166,81,\adfsplit 146,102,10,67,35,73,58,70,59,137,99)_{\adfLq}$,

$(180,12,105,98,88,103,144,130,139,34,179,\adfsplit 40,162,84,99,101,147,117,125,14,104,108)_{\adfLq}$,

$(98,51,83,161,16,125,150,12,39,108,140,\adfsplit 1,92,54,32,24,11,135,29,141,63,162)_{\adfLq}$,

\adfLgap
$(72,39,110,151,31,19,118,160,35,90,58,\adfsplit 113,5,66,109,149,82,180,52,102,103,142)_{\adfLr}$,

$(177,112,1,82,138,125,114,143,0,171,25,\adfsplit 8,22,31,127,96,77,60,63,106,14,179)_{\adfLr}$,

$(165,23,119,25,52,22,63,130,114,180,24,\adfsplit 110,45,5,154,127,36,78,123,10,32,83)_{\adfLr}$,

$(0,144,40,55,86,130,131,152,65,104,68,\adfsplit 13,60,125,43,163,182,76,139,117,22,47)_{\adfLr}$,

$(159,11,90,93,133,164,55,168,65,42,19,\adfsplit 121,15,28,167,162,54,138,71,44,27,79)_{\adfLr}$,

$(174,101,48,54,118,147,160,105,171,1,81,\adfsplit 120,22,10,60,182,65,97,162,87,83,78)_{\adfLr}$,

$(7,166,71,164,47,88,108,38,61,92,173,\adfsplit 8,29,3,27,24,175,5,176,144,89,78)_{\adfLr}$,

$(54,20,136,71,45,159,2,116,145,48,50,\adfsplit 99,100,6,69,157,59,161,94,91,89,166)_{\adfLr}$,

$(176,44,80,106,70,109,23,54,175,155,52,\adfsplit 118,185,133,104,71,47,95,114,75,110,163)_{\adfLr}$,

$(152,66,50,81,59,178,60,146,25,135,83,\adfsplit 102,7,31,5,43,156,63,64,78,183,15)_{\adfLr}$,

\adfLgap
$(122,94,85,139,51,26,14,63,37,27,52,\adfsplit 95,145,7,56,178,167,30,65,148,68,61)_{\adfLs}$,

$(91,183,178,159,73,119,24,14,123,88,20,\adfsplit 148,163,90,44,42,92,21,58,12,32,157)_{\adfLs}$,

$(19,133,48,126,123,114,10,149,107,44,182,\adfsplit 82,9,131,30,53,162,171,83,37,140,7)_{\adfLs}$,

$(124,181,156,96,4,86,51,63,137,68,81,\adfsplit 78,36,154,56,1,17,14,134,18,161,37)_{\adfLs}$,

$(48,109,40,52,168,157,1,66,129,74,118,\adfsplit 81,139,77,141,97,170,107,86,2,88,57)_{\adfLs}$,

$(15,172,127,18,6,29,50,165,98,97,102,\adfsplit 110,10,25,83,157,69,153,115,117,74,34)_{\adfLs}$,

$(133,100,104,115,47,38,140,181,16,117,147,\adfsplit 39,44,85,174,8,161,96,175,26,43,29)_{\adfLs}$,

$(107,171,118,68,84,75,131,156,100,43,76,\adfsplit 158,79,2,9,35,51,66,80,136,38,145)_{\adfLs}$,

$(139,55,101,67,10,53,26,19,62,155,98,\adfsplit 43,144,170,87,127,108,175,66,54,69,169)_{\adfLs}$,

$(8,134,146,155,55,29,117,54,40,96,118,\adfsplit 67,49,157,57,177,78,33,106,126,7,129)_{\adfLs}$

\noindent and

$(83,157,97,60,20,37,179,78,47,105,4,\adfsplit 164,26,175,49,14,87,99,150,156,9,82)_{\adfLt}$,

$(114,105,89,130,65,160,26,106,75,47,49,\adfsplit 112,81,68,25,182,86,148,128,18,70,67)_{\adfLt}$,

$(174,48,80,51,172,28,36,175,61,163,113,\adfsplit 8,67,122,146,23,165,136,138,64,6,32)_{\adfLt}$,

$(160,100,85,84,6,49,70,64,92,154,22,\adfsplit 166,111,167,81,66,89,36,80,79,134,141)_{\adfLt}$,

$(79,23,124,40,160,73,125,184,17,164,129,\adfsplit 77,68,72,76,118,115,86,138,151,89,83)_{\adfLt}$,

$(78,133,56,157,68,100,36,118,33,73,139,\adfsplit 129,105,185,83,167,47,122,175,91,21,96)_{\adfLt}$,

$(41,118,57,34,17,176,167,16,87,131,20,\adfsplit 50,183,137,59,54,93,83,36,7,10,58)_{\adfLt}$,

$(67,168,0,80,15,31,173,2,177,131,88,\adfsplit 169,60,78,99,56,180,51,135,98,82,12)_{\adfLt}$,

$(154,94,11,128,132,177,64,16,41,33,92,\adfsplit 119,27,108,109,30,176,85,84,48,15,76)_{\adfLt}$,

$(154,131,96,85,168,93,37,181,51,183,136,\adfsplit 14,76,99,80,25,21,48,116,151,88,133)_{\adfLt}$

\noindent under the action of the mapping $x \mapsto x + 3$ (mod 132) for $x < 132$,
$x \mapsto 132 + (x - 132 + 5 \mathrm{~(mod~55)})$ for $x \ge 132$.\eproof

\adfVfy{187, \{\{132,44,3\},\{55,11,5\}\}, -1, \{\{22,\{0,1,2,3,4,5\}\},\{55,\{6\}\}\}, -1} 

\vskip 2mm
Theorem~\ref{thm:snark22 22} follows from
Lemmas~\ref{lem:snark22 designs} and \ref{lem:snark22 multipartite}, and Proposition~\ref{prop:d=3, v=22}.


\section{Goldberg's snark \#3}
\label{sec:Goldberg3}
\begin{figure}[h]
\begin{center}
\includegraphics[width=\adfGlarge]{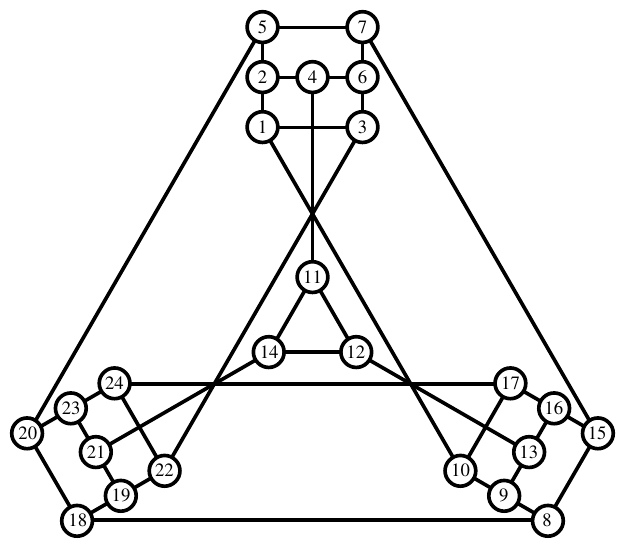}
\caption{Goldberg's snark \#3}
\end{center}
\end{figure}
\newcommand{\adfGc}{\mathrm{GS3}}
Goldberg's snark \#3 is represented by the ordered 24-tuple of its vertices
(1, 2, \dots, 24)${}_{\adfGc}$. Its edge set, as supplied with the {\sc Mathematica} system, \cite{Math-Snarks}, is
 \{$\{9,10\}$, $\{8,9\}$, $\{9,13\}$, $\{13,16\}$, $\{15,16\}$, $\{16,17\}$, $\{10,17\}$, $\{8,15\}$, $\{5,7\}$, $\{1,3\}$, $\{1,2\}$,
   $\{2,5\}$, $\{2,4\}$, $\{4,6\}$, $\{6,7\}$, $\{3,6\}$, $\{18,20\}$, $\{22,24\}$, $\{19,22\}$, $\{18,19\}$, $\{19,21\}$, $\{21,23\}$,
   $\{20,23\}$, $\{23,24\}$, $\{7,15\}$, $\{12,13\}$, $\{1,10\}$, $\{3,22\}$, $\{4,11\}$, $\{17,24\}$, $\{14,21\}$, $\{5,20\}$,
   $\{8,18\}$, $\{11,14\}$, $\{12,14\}$, $\{11,12\}$\}.

{\lemma \label{lem:Goldberg3 designs}
There exist Goldberg's snark \#3 designs of order $64$, $73$, $136$ and $145$.}\\


\noindent\textbf{Proof.}
Let the vertex set of $K_{64}$ be $Z_{63} \cup \{\infty\}$. The decomposition consists of

$(\infty,11,58,48,37,22,13,33,26,21,6,54,\adfsplit 60,45,10,0,8,52,19,2,5,14,59,57)_{\adfGc}$,

$(62,24,28,58,51,3,49,8,23,1,57,45,\adfsplit 41,26,54,2,22,30,60,37,42,13,31,0)_{\adfGc}$,

\vskip 1mm
$(58,15,35,14,5,56,26,60,25,50,2,30,\adfsplit 29,34,1,46,62,32,16,22,23,28,53,0)_{\adfGc}$,

$(40,36,52,8,16,44,34,10,2,12,19,50,\adfsplit 49,23,11,4,55,20,31,17,56,35,29,32)_{\adfGc}$

\noindent under the action of the mapping $\infty \mapsto \infty$, $x \mapsto x + 3$ (mod 63)
for the first two graphs, and $x \mapsto x + 9$ (mod 63) for the last two.

\adfVfy{64, \{\{63,7,9\},\{1,1,1\}\}, 63, -1, \{2,\{\{63,3,3\},\{1,1,1\}\}\}} 

Let the vertex set of $K_{73}$ be $Z_{73}$. The decomposition consists of

$(42,3,33,22,53,28,7,11,14,25,0,1,\adfsplit 16,8,21,34,5,36,71,12,38,45,54,9)_{\adfGc}$

\noindent under the action of the mapping $x \mapsto x + 1$ (mod 73).

\adfVfy{73, \{\{73,73,1\}\}, -1, -1, -1} 

Let the vertex set of $K_{136}$ be $Z_{135} \cup \{\infty\}$. The decomposition consists of

$(\infty,76,105,41,117,33,51,62,85,17,19,110,\adfsplit 116,13,104,81,39,25,57,18,79,111,119,58)_{\adfGc}$,

$(114,46,53,101,70,84,41,45,16,35,26,117,\adfsplit 2,62,97,86,98,12,67,80,65,118,7,55)_{\adfGc}$,

$(31,134,119,20,75,72,86,29,16,32,60,131,\adfsplit 23,3,62,14,40,68,85,125,108,67,43,28)_{\adfGc}$,

$(37,67,94,117,33,61,20,3,103,17,9,19,\adfsplit 53,96,18,49,76,112,5,108,122,30,2,75)_{\adfGc}$,

$(14,119,124,29,78,99,83,53,7,15,126,68,\adfsplit 96,105,102,92,28,107,38,75,43,84,63,25)_{\adfGc}$,

\vskip 1mm
$(97,99,108,60,115,103,54,4,13,88,16,114,\adfsplit 67,30,78,69,52,123,84,72,121,10,100,57)_{\adfGc}$,

$(88,0,67,121,9,31,48,55,97,7,79,51,\adfsplit 37,100,99,82,93,13,73,120,64,78,118,102)_{\adfGc}$

\noindent under the action of the mapping $\infty \mapsto \infty$, $x \mapsto x + 3$ (mod 135)
for the first five graphs, and $x \mapsto x + 9$ (mod 135) for the last two.

\adfVfy{136, \{\{135,15,9\},\{1,1,1\}\}, 135, -1, \{5,\{\{135,3,3\},\{1,1,1\}\}\}} 

Let the vertex set of $K_{145}$ be $Z_{145}$. The decomposition consists of

$(2,106,141,113,58,22,38,62,85,74,0,1,\adfsplit 3,4,5,13,8,6,15,18,28,30,45,59)_{\adfGc}$,

$(0,18,19,39,40,64,2,1,29,59,3,34,\adfsplit 71,73,51,4,103,82,24,5,126,84,74,13)_{\adfGc}$

\noindent under the action of the mapping $x \mapsto x + 1$ (mod 145).\eproof

\adfVfy{145, \{\{145,145,1\}\}, -1, -1, -1} 

{\lemma \label{lem:Goldberg3 multipartite}
There exist decompositions of the complete multipartite graphs
$K_{12,12,12}$, $K_{24,24,15}$, $K_{72,72,63}$, $K_{24,24,24,24}$ and $K_{24,24,24,21}$
into Goldberg's snark \#3.}\\


\noindent\textbf{Proof.}
Let the vertex set of $K_{12,12,12}$ be $Z_{36}$ partitioned according to residue classes modulo 3.
The decomposition consists of

$(11,30,12,13,28,14,24,18,19,6,2,7,\adfsplit 35,27,32,22,29,10,15,3,1,8,23,25)_{\adfGc}$

\noindent under the action of the mapping $x \mapsto x + 3$ (mod 36).

\adfVfy{36, \{\{36,12,3\}\}, -1, \{\{12,\{0,1,2\}\}\}, -1} 

Let the vertex set of $K_{24,24,15}$ be $\{0, 1, \dots, 62\}$ partitioned into
$\{2j + i: j = 0, 1, \dots, 23\}$, $i = 0, 1$, and $\{48, 49, \dots, 62\}$.
The decompositions consist of

$(10,27,13,12,61,56,32,18,17,55,39,50,\adfsplit 40,8,29,54,47,58,19,25,60,22,42,52)_{\adfGc}$,

$(16,52,37,1,12,20,23,48,27,58,54,14,\adfsplit 49,9,30,35,36,0,43,19,61,10,46,11)_{\adfGc}$,

$(37,32,2,23,17,61,42,24,57,30,14,50,\adfsplit 33,29,13,20,7,54,10,8,58,39,9,52)_{\adfGc}$

\noindent under the action of the mapping $x \mapsto x + 4$ (mod 48) for $x < 48$,
$x \mapsto 48 + (x - 48 + 5 \mathrm{~(mod~15)})$ for $x \ge 48$.

\adfVfy{63, \{\{48,12,4\},\{15,3,5\}\}, -1, \{\{24,\{0,1\}\},\{15,\{2\}\}\}, -1} 

Let the vertex set of $K_{72,72,63}$ be $\{0, 1, \dots, 206\}$ partitioned into
$\{2j + i: j = 0, 1, \dots, 71\}$, $i = 0, 1$, and $\{144, 145, \dots, 206\}$.
The decompositions consist of

$(186,47,13,140,171,169,62,82,188,43,99,120,\adfsplit 95,158,174,78,147,39,130,8,65,164,165,139)_{\adfGc}$,

$(189,48,14,141,172,170,63,83,191,44,100,121,\adfsplit 96,159,175,79,148,40,131,9,66,165,166,140)_{\adfGc}$,

$(192,49,15,142,173,171,64,84,194,45,101,122,\adfsplit 97,160,176,80,149,41,132,10,67,166,167,141)_{\adfGc}$,

$(195,50,16,143,174,172,65,85,197,46,102,123,\adfsplit 98,161,177,81,150,42,133,11,68,167,168,142)_{\adfGc}$,

$(150,83,32,84,164,47,86,191,13,6,174,52,\adfsplit 203,55,1,90,153,115,18,142,202,113,3,136)_{\adfGc}$,

$(151,84,33,85,165,48,87,194,14,7,175,53,\adfsplit 206,56,2,91,154,116,19,143,205,114,4,137)_{\adfGc}$,

$(152,85,34,86,166,49,88,197,15,8,176,54,\adfsplit 182,57,3,92,155,117,20,0,181,115,5,138)_{\adfGc}$,

$(153,86,35,87,167,50,89,200,16,9,177,55,\adfsplit 185,58,4,93,156,118,21,1,184,116,6,139)_{\adfGc}$,

$(32,141,3,171,84,22,13,143,148,61,12,89,\adfsplit 118,180,82,47,28,158,128,9,109,190,184,63)_{\adfGc}$,

$(102,153,89,87,112,186,39,64,199,121,78,3,\adfsplit 60,196,192,127,50,193,22,125,105,183,38,5)_{\adfGc}$,

$(192,46,48,180,79,57,140,116,201,38,3,187,\adfsplit 6,12,183,35,193,47,2,92,196,15,137,102)_{\adfGc}$

\noindent under the action of the mapping $x \mapsto x + 4$ (mod 144) for $x < 144$,
$x \mapsto 144 + (x + 4 \mathrm{~(mod~36)})$ for $144 \le x < 180$,
$x \mapsto 180 + (x - 180 + 12 \mathrm{~(mod~27)})$ for $x \ge 180$.

\adfVfy{207, \{\{144,36,4\},\{36,9,4\},\{27,9,12\}\}, -1, \{\{72,\{0,1\}\},\{63,\{2\}\}\}, -1} 

Let the vertex set of $K_{24,24,24,24}$ be $Z_{96}$ partitioned according to residue classes modulo 4.
The decomposition consists of

$(33,62,24,92,51,85,36,44,47,84,91,78,\adfsplit 72,37,18,41,2,6,80,49,3,34,82,7)_{\adfGc}$

\noindent under the action of the mapping $x \mapsto x + 1$ (mod 96).

\adfVfy{96, \{\{96,96,1\}\}, -1, \{\{24,\{0,1,2,3\}\}\}, -1} 

Let the vertex set of $K_{24,24,24,21}$ be $\{0, 1, \dots, 92\}$ partitioned into
$\{2j + i: j = 0, 1, \dots, 23\}$, $i = 0, 1, 2$, and $\{72, 73, \dots, 92\}$.
The decompositions consist of

$(13,72,44,57,65,81,48,34,8,45,83,41,\adfsplit 10,46,53,18,80,59,42,88,47,19,15,20)_{\adfGc}$,

$(88,68,45,60,30,58,20,63,2,54,53,61,\adfsplit 91,6,85,5,84,29,48,1,75,19,51,23)_{\adfGc}$,

$(38,75,37,41,61,57,85,28,59,48,67,3,\adfsplit 90,5,47,54,35,17,31,60,15,73,32,6)_{\adfGc}$,

$(49,42,35,28,11,91,48,27,78,24,90,5,\adfsplit 66,7,61,38,82,81,30,46,87,72,68,3)_{\adfGc}$,

$(44,40,4,17,63,87,11,2,76,24,21,75,\adfsplit 28,43,70,14,85,52,36,84,72,86,26,66)_{\adfGc}$

\noindent under the action of the mapping $x \mapsto x + 4$ (mod 72) for $x < 72$,
$x \mapsto 72 + (x - 72 + 7 \mathrm{~(mod~21)})$ for $x \ge 72$.\eproof

\adfVfy{93, \{\{72,18,4\},\{21,3,7\}\}, -1, \{\{24,\{0,1,2\}\},\{21,\{3\}\}\}, -1} 

\vskip 2mm
Theorem~\ref{thm:Goldberg3 24} follows from Lemmas~\ref{lem:Goldberg3 designs} and \ref{lem:Goldberg3 multipartite},
and Proposition~\ref{prop:d=3, v=24}.


\section{The two Celmins--Swart snarks and the two 26-vertex Blanu\v{s}a snarks}
\label{sec:CelminsSwart}
\label{sec:Blanusa26 26}
\begin{figure}[h]
\begin{center}
\includegraphics[width=\adfGlarge]{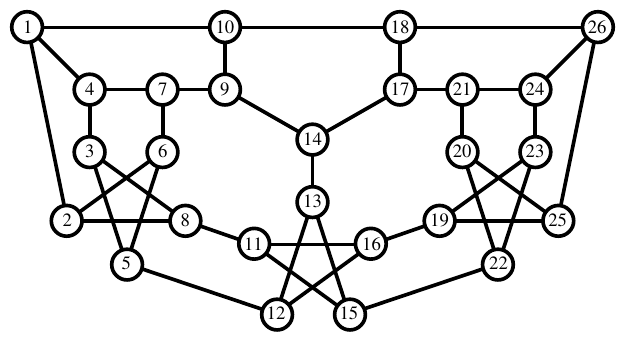}

\includegraphics[width=\adfGlarge]{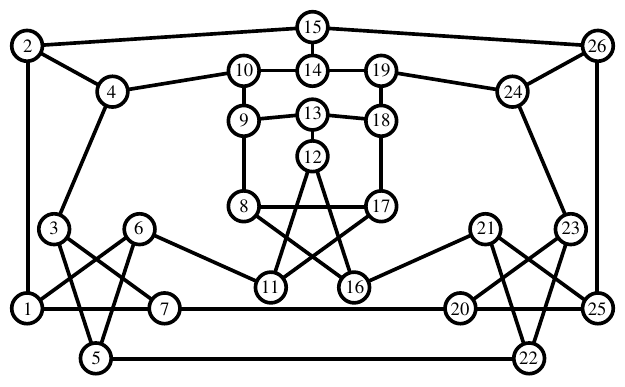}
\caption{The Celmins--Swart snarks \#1 and \#2}
\end{center}
\end{figure}
\begin{figure}[h]
\begin{center}
\includegraphics[width=\adfGlarge]{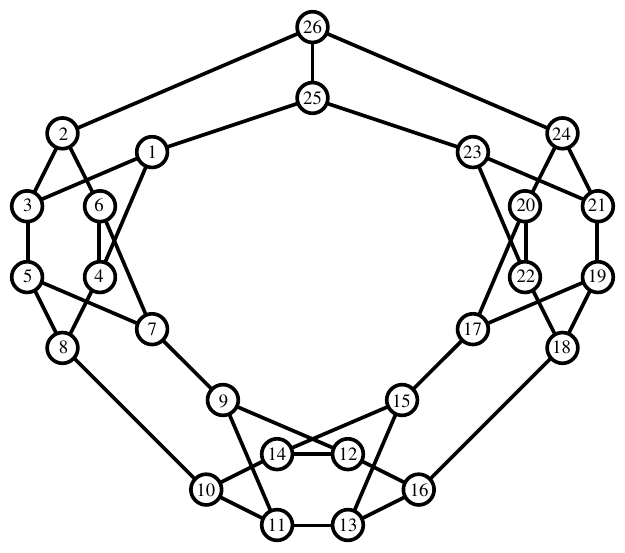}
\caption{The Blanu\v{s}a snark (2,1)}
\end{center}
\end{figure}
\begin{figure}[h]
\begin{center}
\includegraphics[width=\adfGlarge]{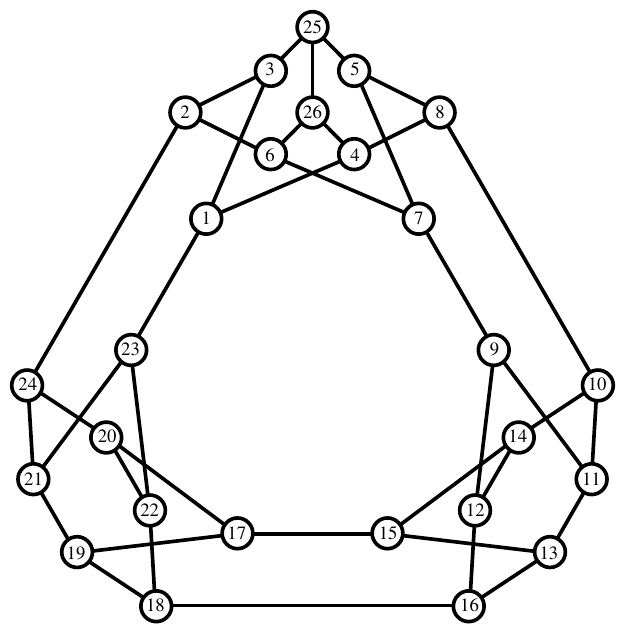}
\caption{The Blanu\v{s}a snark (2,2)}
\end{center}
\end{figure}
\newcommand{\adfCWa}{\mathrm{CS1}}
\newcommand{\adfCWb}{\mathrm{CS2}}
\newcommand{\adfBba}{\mathrm{B21}}
\newcommand{\adfBbb}{\mathrm{B22}}
The Celmins--Swart snarks and the 26-vertex Blanu\v{s}a snarks are each represented by
the ordered 26-tuple of its vertices:
(1, 2, \dots, 26)${}_{\adfCWa}$ for Celmins--Swart snark \#1,
(1, 2, \dots, 26)${}_{\adfCWb}$ for Celmins--Swart snark \#2,
(1, 2, \dots, 26)${}_{\adfBba}$ for the $(2,1)$-Blanu\v{s}a snark and
(1, 2, \dots, 26)${}_{\adfBbb}$ for the $(2,2)$-Blanu\v{s}a snark.
The edge sets, as supplied with the {\sc Mathematica} system, \cite{Math-Snarks}, are respectively

$\adfCWa$:
  \{$\{1,2\}$, $\{1,4\}$, $\{1,10\}$, $\{2,6\}$, $\{2,8\}$, $\{3,4\}$, $\{3,5\}$, $\{3,8\}$, $\{4,7\}$, $\{5,6\}$, $\{5,12\}$, $\{6,7\}$, $\{7,9\}$,
   $\{8,11\}$, $\{9,10\}$, $\{9,14\}$, $\{10,18\}$, $\{11,15\}$, $\{11,16\}$, $\{12,13\}$, $\{12,16\}$, $\{13,14\}$, $\{13,15\}$, $\{14,17\}$,
   $\{15,22\}$, $\{16,19\}$, $\{17,18\}$, $\{17,21\}$, $\{18,26\}$, $\{19,23\}$, $\{19,25\}$, $\{20,21\}$, $\{20,22\}$, $\{20,25\}$,
   $\{21,24\}$, $\{22,23\}$, $\{23,24\}$, $\{24,26\}$, $\{25,26\}$\},

$\adfCWb$:
   \{$\{1,6\}$, $\{1,2\}$, $\{3,4\}$, $\{3,7\}$, $\{2,4\}$, $\{3,5\}$, $\{5,6\}$, $\{1,7\}$, $\{7,20\}$, $\{8,9\}$, $\{8,17\}$, $\{9,10\}$, $\{4,10\}$,
   $\{10,14\}$, $\{6,11\}$, $\{11,12\}$, $\{12,13\}$, $\{12,16\}$, $\{9,13\}$, $\{13,18\}$, $\{14,15\}$, $\{14,19\}$, $\{2,15\}$,
   $\{15,26\}$, $\{8,16\}$, $\{16,21\}$, $\{11,17\}$, $\{17,18\}$, $\{18,19\}$, $\{19,24\}$, $\{20,25\}$, $\{21,22\}$, $\{5,22\}$,
   $\{22,23\}$, $\{20,23\}$, $\{23,24\}$, $\{24,26\}$, $\{21,25\}$, $\{25,26\}$\},

$\adfBba$:
  \{$\{1,3\}$, $\{1,4\}$, $\{2,3\}$, $\{2,6\}$, $\{3,5\}$, $\{4,6\}$, $\{4,8\}$, $\{5,7\}$, $\{5,8\}$, $\{6,7\}$, $\{7,9\}$, $\{8,10\}$, $\{9,11\}$, $\{9,12\}$, $\{10,11\}$, $\{10,14\}$, $\{11,13\}$, $\{12,14\}$, $\{12,16\}$, $\{13,15\}$, $\{13,16\}$, $\{14,15\}$, $\{15,17\}$, $\{16,18\}$, $\{17,19\}$, $\{17,20\}$, $\{18,19\}$, $\{18,22\}$, $\{19,21\}$, $\{20,22\}$, $\{20,24\}$, $\{21,23\}$, $\{21,24\}$, $\{22,23\}$, $\{23,25\}$, $\{24,26\}$, $\{1,25\}$, $\{2,26\}$, $\{25,26\}$\}
and

$\adfBbb$:
  \{$\{1,3\}$, $\{1,4\}$, $\{2,3\}$, $\{2,6\}$, $\{4,8\}$, $\{5,7\}$, $\{5,8\}$, $\{6,7\}$, $\{7,9\}$, $\{8,10\}$, $\{9,11\}$, $\{9,12\}$, $\{10,11\}$, $\{10,14\}$, $\{11,13\}$, $\{12,14\}$, $\{12,16\}$, $\{13,15\}$, $\{13,16\}$, $\{14,15\}$, $\{15,17\}$, $\{16,18\}$, $\{17,19\}$, $\{17,20\}$, $\{18,19\}$, $\{18,22\}$, $\{19,21\}$, $\{20,22\}$, $\{20,24\}$, $\{21,23\}$, $\{21,24\}$, $\{22,23\}$, $\{23,1\}$, $\{24,2\}$, $\{3,25\}$, $\{4,26\}$, $\{5,25\}$, $\{6,26\}$, $\{25,26\}$\}.


{\lemma \label{lem:CelminsSwart designs}
There exist designs of order $40$ and $79$ for each of the two Celmins--Swart snarks
and the two 26-vertex Blanu\v{s}a snarks.}\\


\noindent\textbf{Proof.}
Let the vertex set of $K_{40}$ be $Z_{40}$. The decompositions consist of

$(15,20,1,11,30,37,10,2,7,16,22,6,0,\adfsplit 27,23,31,36,4,17,19,38,5,35,12,33,32)_{\adfCWa}$,

$(36,30,29,15,9,7,23,25,10,32,21,12,1,\adfsplit 20,16,28,5,34,35,18,26,27,19,2,33,6)_{\adfCWa}$,

$(30,35,37,22,7,13,21,0,27,33,19,5,17,\adfsplit 23,4,16,20,25,9,11,38,14,26,29,18,8)_{\adfCWa}$,

$(11,26,29,13,6,32,18,28,23,12,15,24,14,16,\adfsplit 3,0,8,21,36,38,17,31,10,39,2,9)_{\adfCWa}$,

\adfLgap
$(1,34,24,28,9,31,19,8,10,32,18,4,3,\adfsplit 12,36,11,6,2,25,22,30,15,26,13,33,27)_{\adfCWb}$,

$(20,37,36,4,15,23,8,30,12,14,26,18,33,\adfsplit 7,24,3,38,25,27,5,32,0,10,31,13,17)_{\adfCWb}$,

$(32,38,9,33,0,31,18,2,19,36,37,22,11,\adfsplit 26,29,1,13,35,5,23,12,17,21,20,6,15)_{\adfCWb}$,

$(27,0,36,7,23,39,22,24,29,9,11,2,20,\adfsplit 21,19,18,25,35,31,32,38,5,16,30,14,13)_{\adfCWb}$,

\adfLgap
$(12,6,32,36,25,26,18,31,5,22,34,16,37,\adfsplit 14,33,11,17,24,28,7,30,29,21,27,3,10)_{\adfBba}$,

$(22,19,35,16,24,15,21,6,36,31,0,2,18,\adfsplit 28,3,20,38,17,5,9,4,14,13,23,37,7)_{\adfBba}$,

$(0,2,16,19,28,18,33,20,6,7,35,23,12,\adfsplit 26,38,25,22,37,3,4,1,14,10,11,9,8)_{\adfBba}$,

$(22,10,15,23,33,8,29,16,34,24,7,17,18,\adfsplit 9,27,37,1,39,35,20,3,31,5,13,28,0)_{\adfBba}$

\noindent and

$(1,5,0,19,17,4,26,32,37,3,29,31,13,\adfsplit 9,33,10,20,18,28,35,15,30,36,27,34,21)_{\adfBbb}$,

$(31,1,6,29,0,33,19,32,23,13,24,8,34,\adfsplit 30,9,14,7,5,11,4,2,35,18,37,20,22)_{\adfBbb}$,

$(36,18,14,17,34,12,35,27,10,28,32,15,22,\adfsplit 37,38,3,30,23,31,19,20,6,8,16,1,7)_{\adfBbb}$,

$(29,28,4,16,0,35,33,7,5,31,17,15,8,\adfsplit 22,34,32,38,9,13,1,39,12,10,2,14,11)_{\adfBbb}$

\noindent under the action of the mapping $x \mapsto x + 8$ (mod 40).

\adfVfy{40, \{\{40,5,8\}\}, -1, -1, -1} 

Let the vertex set of $K_{79}$ be $Z_{79}$. The decompositions consist of

$(62,53,77,59,5,19,74,6,51,24,17,3,56,\adfsplit 7,57,13,1,29,0,9,55,21,42,71,31,12)_{\adfCWa}$,

\adfLgap
$(17,24,76,77,12,22,68,57,20,74,0,40,38,\adfsplit 7,53,2,9,3,23,6,45,59,25,46,72,42)_{\adfCWb}$,

\adfLgap
$(69,10,42,12,14,46,52,0,67,71,38,36,3,\adfsplit 15,40,6,1,2,11,8,66,68,7,61,70,72)_{\adfBba}$

\noindent and

$(42,73,38,32,69,72,2,48,43,37,0,41,34,\adfsplit 18,66,1,4,6,12,29,25,58,28,47,10,17)_{\adfBbb}$

\noindent under the action of the mapping $x \mapsto x + 1$ (mod 79).\eproof

\adfVfy{79, \{\{79,79,1\}\}, -1, -1, -1} 

{\lemma \label{lem:CelminsSwart multipartite}
There exist decompositions of $K_{39,39,39}$ and $K_{13,13,13,13}$ into each of the two Celmins--Swart snarks
and the two 26-vertex Blanu\v{s}a snarks.}\\


\noindent\textbf{Proof.}
Let the vertex set of $K_{39,39,39}$ be $Z_{117}$ partitioned according to residue classes modulo 3.
The decompositions consist of

$(64,44,53,72,111,85,5,10,78,92,32,26,51,\adfsplit 68,28,19,0,31,3,8,46,90,91,86,73,84)_{\adfCWa}$,

\adfLgap
$(27,94,49,54,66,101,98,50,84,107,102,5,19,\adfsplit 79,87,3,40,8,21,0,47,10,41,37,25,116)_{\adfCWb}$,

\adfLgap
$(75,53,40,5,8,51,85,22,78,17,79,115,11,\adfsplit 90,112,21,3,1,29,34,48,77,19,98,71,109)_{\adfBba}$

\noindent and

$(47,25,72,43,59,110,10,45,62,67,75,31,11,\adfsplit 87,4,12,9,1,44,82,90,27,106,2,100,33)_{\adfBbb}$

\noindent under the action of the mapping $x \mapsto x + 1$ (mod 117).

\adfVfy{117, \{\{117,117,1\}\}, -1, \{\{39,\{0,1,2\}\}\}, -1} 

Let the vertex set of $K_{13,13,13,13}$ be $Z_{52}$ partitioned according to residue classes modulo 4.
The decompositions consist of

$(15,34,11,0,26,21,39,4,1,12,35,3,32,\adfsplit 46,9,41,47,25,31,14,49,28,18,24,48,6)_{\adfCWa}$,

$(4,31,49,26,28,21,23,32,6,5,34,13,16,\adfsplit 11,7,43,18,39,2,50,9,1,48,46,29,20)_{\adfCWa}$,

\adfLgap
$(26,32,27,25,48,7,13,5,38,15,50,11,36,\adfsplit 18,41,20,16,3,8,39,35,34,4,9,46,19)_{\adfCWb}$,

$(10,20,5,46,42,27,40,23,17,15,24,22,4,\adfsplit 16,41,29,1,18,49,11,30,33,6,0,25,7)_{\adfCWb}$,

\adfLgap
$(24,5,6,3,7,34,17,48,15,18,13,2,26,\adfsplit 8,35,40,49,14,16,46,33,19,4,44,43,50)_{\adfBba}$,

$(18,7,37,45,19,30,41,24,50,49,15,35,12,\adfsplit 36,21,9,31,16,17,25,11,46,28,20,51,10)_{\adfBba}$

\noindent and

$(18,21,24,16,11,3,48,9,15,34,37,6,42,\adfsplit 27,33,23,28,13,26,19,12,32,49,22,10,29)_{\adfBbb}$,

$(23,28,17,12,39,31,2,34,15,16,45,25,43,\adfsplit 42,8,18,35,40,49,4,38,46,9,27,44,13)_{\adfBbb}$

\noindent under the action of the mapping $x \mapsto x + 4$ (mod 52).\eproof

\adfVfy{52, \{\{52,13,4\}\}, -1, \{\{13,\{0,1,2,3\}\}\}, -1} 

\vskip 2mm
Theorem~\ref{thm:Snark partial} (i) follows from
Lemmas~\ref{lem:CelminsSwart designs} and \ref{lem:CelminsSwart multipartite}, and
Proposition~\ref{prop:e + 1, 2e + 1, e^3, (e/3)^4 => et + 1}.


\section{The flower snark J7}
\label{sec:FlowerSnarkJ7}
\begin{figure}[h]
\begin{center}
\includegraphics[width=\adfGlarge]{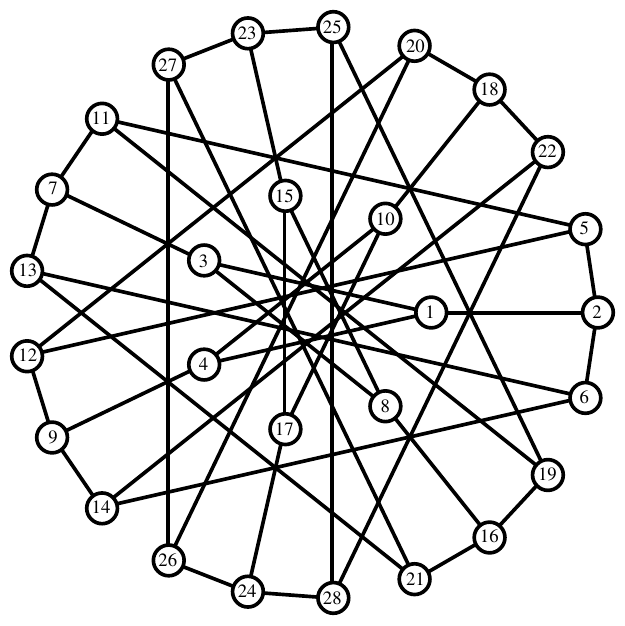}
\caption{The Flower snark J7}
\end{center}
\end{figure}
\newcommand{\adfFJg}{\mathrm{FJ7}}
The flower snark J7 is represented by the ordered 28-tuple of its vertices
(1, 2, \dots, 28)${}_{\adfFJg}$ and edge set
  \{$\{1,2\}$, $\{1,3\}$, $\{1,4\}$, $\{2,5\}$, $\{2,6\}$, $\{3,7\}$, $\{3,8\}$, $\{4,9\}$, $\{4,10\}$, $\{5,11\}$, $\{5,12\}$,
   $\{6,13\}$, $\{6,14\}$, $\{7,11\}$, $\{7,13\}$, $\{8,15\}$, $\{8,16\}$, $\{9,12\}$, $\{9,14\}$, $\{10,17\}$, $\{10,18\}$,
   $\{11,19\}$, $\{12,20\}$, $\{13,21\}$, $\{14,22\}$, $\{15,17\}$, $\{15,23\}$, $\{16,19\}$, $\{16,21\}$, $\{17,24\}$,
   $\{18,20\}$, $\{18,22\}$, $\{19,25\}$, $\{20,26\}$, $\{21,27\}$, $\{22,28\}$, $\{23,25\}$, $\{23,27\}$, $\{24,26\}$, $\{24,28\}$,
   $\{25,28\}$, $\{26,27\}$\}
as supplied with the {\sc Mathematica} system, \cite{Math-Snarks}.


{\lemma \label{lem:FlowerSnarkJ7 designs}
There exist flower snark J7 designs of order $28$, $85$ and $169$.}\\


\noindent\textbf{Proof.}
Let the vertex set of $K_{28}$ be $Z_{27} \cup \{\infty\}$. The decomposition consists of

$(\infty,15,26,1,25,5,16,8,7,4,18,12,11,3,10,\adfsplit 23,19,17,24,13,22,14,2,0,21,20,6,9)_{\adfFJg}$

\noindent under the action of the mapping $\infty \mapsto \infty$, $x \mapsto x + 3$ (mod 27).

\adfVfy{28, \{\{27,9,3\},\{1,1,1\}\}, 27, -1, -1} 

Let the vertex set of $K_{85}$ be $Z_{85}$. The decomposition consists of

$(21,37,77,11,72,22,52,32,4,83,28,25,64,24,\adfsplit 58,0,1,2,5,3,6,10,27,56,38,39,73,19)_{\adfFJg}$

\noindent under the action of the mapping $x \mapsto x + 1$ (mod 85).

\adfVfy{85, \{\{85,85,1\}\}, -1, -1, -1} 

Let the vertex set of $K_{169}$ be $Z_{169}$. The decomposition consists of

$(117,136,116,167,165,145,149,121,130,92,78,4,90,142,\adfsplit 134,40,152,7,110,30,156,111,3,52,138,148,27,91)_{\adfFJg}$,

$(128,126,164,112,25,80,149,5,116,137,86,19,1,3,\adfsplit
12,16,26,9,121,36,68,83,69,75,99,151,24,161)_{\adfFJg}$

\noindent under the action of the mapping $x \mapsto x + 1$ (mod 169).\eproof

\adfVfy{169, \{\{169,169,1\}\}, -1, -1, -1} 

{\lemma \label{lem:FlowerSnarkJ7 multipartite}
There exist decompositions of $K_{42,42,42}$ and $K_{7,7,7,7}$ into the flower snark J7.}\\


\noindent\textbf{Proof.}
Let the vertex set of $K_{42,42,42}$ be $Z_{126}$ partitioned according to residue classes modulo 3.
The decomposition consists of

$(42,113,74,16,22,63,93,99,47,80,83,0,4,7,\adfsplit 1,2,3,6,10,5,15,89,123,37,17,117,100,105)_{\adfFJg}$

\noindent under the action of the mapping $x \mapsto x + 1$ (mod 126).

\adfVfy{126, \{\{126,126,1\}\}, -1, \{\{42,\{0,1,2\}\}\}, -1} 

Let the vertex set of $K_{7,7,7,7}$ be $Z_{28}$ partitioned according to residue classes modulo 4.
The decomposition consists of

$(25,16,7,11,5,22,20,17,9,2,3,8,13,4,\adfsplit 10,23,27,1,24,18,0,26,12,6,15,19,14,21)_{\adfFJg}$

\noindent under the action of the mapping $x \mapsto x + 4$ (mod 28).\eproof

\adfVfy{28, \{\{28,7,4\}\}, -1, \{\{7,\{0,1,2,3\}\}\}, -1} 

{\lemma \label{lem:FlowerSnarkJ7 28 mod 84}
There exist a flower snark J7 design of order $n$ if $n \equiv 28$ \textup{(mod 84)}.}\\

\noindent\textbf{Proof.} The case $n=28$ follows from Lemma~\ref{lem:FlowerSnarkJ7 designs}.
There exists a 4-GDD of type $4^{3t+1}$ for $t\ge 1$, \cite{BSH} (see also \cite{Ge}).
Inflate by a factor of $7$, lay a complete graph $K_{28}$ (from Lemma~\ref{lem:FlowerSnarkJ7 designs})
on each of the inflated groups and replace each block by a complete 4-partite graph $K_{7,7,7,7}$
(from Lemma~\ref{lem:FlowerSnarkJ7 multipartite}) to yield
a flower snark J7 design of order $84t + 28$ for $t \ge 1$.\eproof


\vskip 2mm
Theorem~\ref{thm:Snark partial} (ii) follows from
Lemmas~\ref{lem:FlowerSnarkJ7 designs} and \ref{lem:FlowerSnarkJ7 multipartite}, and
Proposition~\ref{prop:2e + 1, 4e + 1, e^3 => 2et + 1} for $n \equiv 1$ (mod 84),
and Lemma~\ref{lem:FlowerSnarkJ7 28 mod 84} for $n \equiv 28$ (mod 84).


\section{The double star snark}
\label{sec:DoubleStar}
\begin{figure}[h]
\begin{center}
\includegraphics[width=\adfGlarge]{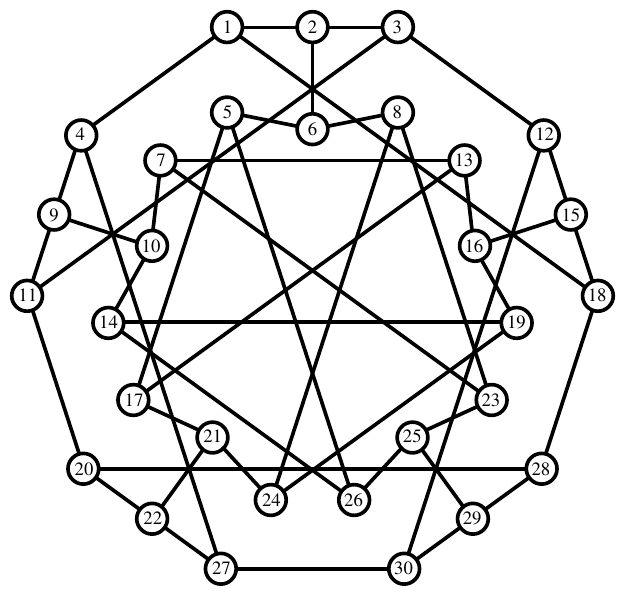}
\caption{The double star snark}
\end{center}
\end{figure}
\newcommand{\adfDS}{\mathrm{DS}}
The double star snark is represented by the ordered 30-tuple of its vertices
(1, 2, \dots, 30)${}_{\adfDS}$ and edge set
  \{$\{13,16\}$, $\{15,16\}$, $\{12,15\}$, $\{16,19\}$, $\{15,18\}$, $\{23,25\}$, $\{25,29\}$, $\{28,29\}$, $\{25,26\}$, $\{29,30\}$,
   $\{21,24\}$, $\{21,22\}$, $\{22,27\}$, $\{17,21\}$, $\{20,22\}$, $\{10,14\}$, $\{9,10\}$, $\{9,11\}$, $\{7,10\}$, $\{4,9\}$, $\{5,6\}$,
   $\{2,6\}$, $\{1,2\}$, $\{6,8\}$, $\{2,3\}$, $\{14,19\}$, $\{7,13\}$, $\{5,26\}$, $\{8,23\}$, $\{13,17\}$, $\{19,24\}$, $\{7,23\}$,
   $\{14,26\}$, $\{8,24\}$, $\{5,17\}$, $\{1,18\}$, $\{12,30\}$, $\{20,28\}$, $\{4,27\}$, $\{3,11\}$, $\{18,28\}$, $\{27,30\}$, $\{11,20\}$,
   $\{1,4\}$, $\{3,12\}$\}
as supplied with the {\sc Mathematica} system, \cite{Math-Snarks}.


{\lemma \label{lem:DoubleStar designs}
There exist double star snark designs of order $46$, $91$ and $181$.}\\


\noindent\textbf{Proof.}
Let the vertex set of $K_{46}$ be $Z_{46}$. The decomposition consists of

$(41,10,13,24,6,16,4,21,5,43,39,12,32,0,23,\adfsplit 28,19,27,8,2,40,17,38,15,36,22,33,37,9,11)_{\adfDS}$

\noindent under the action of the mapping $x \mapsto x + 2$ (mod 46).

\adfVfy{46, \{\{46,23,2\}\}, -1, -1, -1} 

Let the vertex set of $K_{91}$ be $Z_{91}$. The decomposition consists of

$(10,1,48,50,59,11,49,47,27,69,16,34,4,77,31,88,\adfsplit 6,85,12,20,0,19,8,25,70,42,83,15,82,52)_{\adfDS}$

\noindent under the action of the mapping $x \mapsto x + 1$ (mod 91).

\adfVfy{91, \{\{91,91,1\}\}, -1, -1, -1} 

Let the vertex set of $K_{181}$ be $Z_{181}$. The decomposition consists of

$(158,77,118,129,105,50,98,83,70,25,133,143,14,162,117,\adfsplit 79,1,26,88,21,0,2,3,5,6,12,9,4,16,20)_{\adfDS}$,

$(0,8,22,16,1,29,2,5,34,54,57,52,36,3,4,\adfsplit 68,72,43,105,7,134,64,111,58,15,128,131,142,102,10)_{\adfDS}$

\noindent under the action of the mapping $x \mapsto x + 1$ (mod 181).\eproof

\adfVfy{181, \{\{181,181,1\}\}, -1, -1, -1} 

{\lemma \label{lem:DoubleStar multipartite}
There exists a decomposition of $K_{15,15,15}$ into the double star snark.}\\


\noindent\textbf{Proof.}
Let the vertex set of $K_{15,15,15}$ be $Z_{45}$ partitioned according to residue classes modulo 3.
The decomposition consists of

$(26,15,14,7,25,38,34,31,17,21,33,6,18,22,28,\adfsplit 32,5,36,27,37,19,39,24,2,43,20,41,29,12,40)_{\adfDS}$

\noindent under the action of the mapping $x \mapsto x + 3$ (mod 45).\eproof

\adfVfy{45, \{\{45,15,3\}\}, -1, \{\{15,\{0,1,2\}\}\}, -1} 

\vskip 2mm
Theorem~\ref{thm:Snark partial} (iii) follows from
Lemmas~\ref{lem:DoubleStar designs} and \ref{lem:DoubleStar multipartite}, and
Proposition~\ref{prop:e + 1, 2e + 1, 4e + 1, (e/3)^3 => et + 1}.


\section{The two 34-vertex Blanu\v{s}a snarks}
\label{sec:Blanusa34 34}
\begin{figure}[h]
\begin{center}
\includegraphics[width=\adfGlarge]{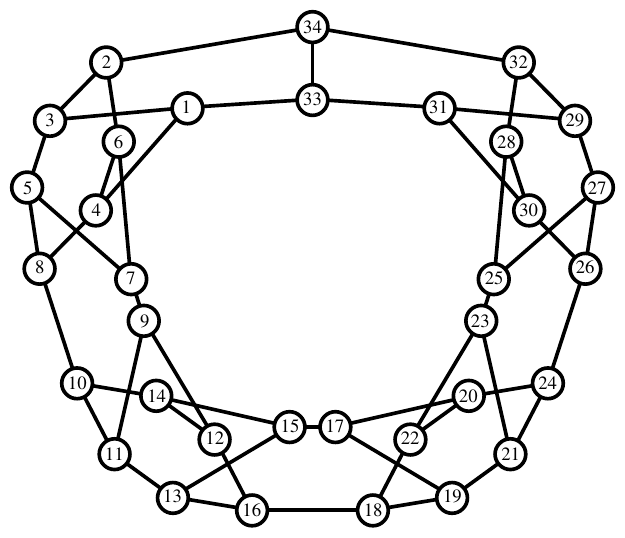}
\caption{The Blanu\v{s}a snark (3,1)}
\end{center}
\end{figure}
\begin{figure}[h]
\begin{center}
\includegraphics[width=\adfGlarge]{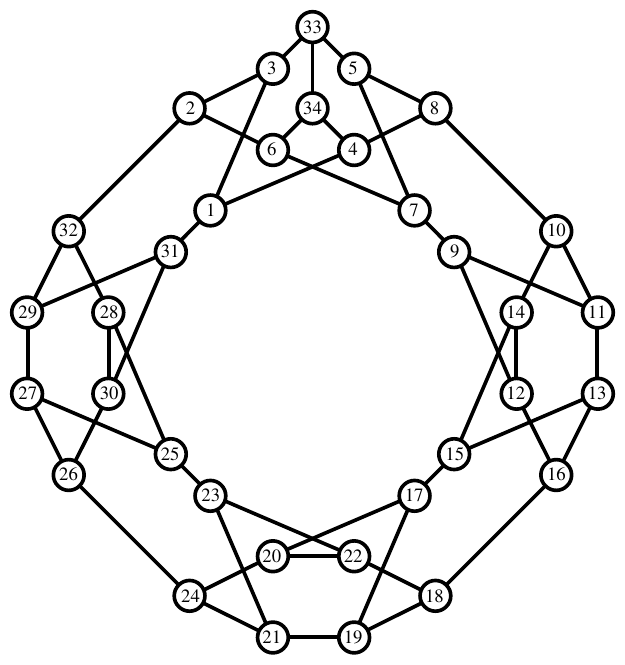}
\caption{The Blanu\v{s}a snark (3,2)}
\end{center}
\end{figure}
\newcommand{\adfBca}{\mathrm{B31}}
\newcommand{\adfBcb}{\mathrm{B32}}
The two 34-vertex Blanu\v{s}a snarks are each represented by the ordered 36-tuple of its vertices,
(1, 2, \dots, 34)${}_{\adfBca}$ for the $(3,1)$-Blanu\v{s}a snark and
(1, 2, \dots, 34)${}_{\adfBcb}$ for the $(3,2)$-Blanu\v{s}a snark.
The edge sets, as supplied with the {\sc Mathematica} system, \cite{Math-Snarks}, are respectively

$\adfBca$:
\{$\{1,3\}$, $\{1,4\}$, $\{2,3\}$, $\{2,6\}$, $\{3,5\}$, $\{4,6\}$, $\{4,8\}$, $\{5,7\}$, $\{5,8\}$, $\{6,7\}$, $\{7,9\}$, $\{8,10\}$, $\{9,11\}$, $\{9,12\}$, $\{10,11\}$, $\{10,14\}$, $\{11,13\}$, $\{12,14\}$, $\{12,16\}$, $\{13,15\}$, $\{13,16\}$, $\{14,15\}$, $\{15,17\}$, $\{16,18\}$, $\{17,19\}$, $\{17,20\}$, $\{18,19\}$, $\{18,22\}$, $\{19,21\}$, $\{20,22\}$, $\{20,24\}$, $\{21,23\}$, $\{21,24\}$, $\{22,23\}$, $\{23,25\}$, $\{24,26\}$, $\{25,27\}$, $\{25,28\}$, $\{26,27\}$, $\{26,30\}$, $\{27,29\}$, $\{28,30\}$, $\{28,32\}$, $\{29,31\}$, $\{29,32\}$, $\{30,31\}$, $\{31,33\}$, $\{32,34\}$, $\{33,1\}$, $\{34,2\}$, $\{33,34\}$\}
and

$\adfBcb$:
\{$\{1,3\}$, $\{1,4\}$, $\{2,3\}$, $\{2,6\}$, $\{4,8\}$, $\{5,7\}$, $\{5,8\}$, $\{6,7\}$, $\{7,9\}$, $\{8,10\}$, $\{9,11\}$, $\{9,12\}$, $\{10,11\}$, $\{10,14\}$, $\{11,13\}$, $\{12,14\}$, $\{12,16\}$, $\{13,15\}$, $\{13,16\}$, $\{14,15\}$, $\{15,17\}$, $\{16,18\}$, $\{17,19\}$, $\{17,20\}$, $\{18,19\}$, $\{18,22\}$, $\{19,21\}$, $\{20,22\}$, $\{20,24\}$, $\{21,23\}$, $\{21,24\}$, $\{22,23\}$, $\{23,25\}$, $\{24,26\}$, $\{25,27\}$, $\{25,28\}$, $\{26,27\}$, $\{26,30\}$, $\{27,29\}$, $\{28,30\}$, $\{28,32\}$, $\{29,31\}$, $\{29,32\}$, $\{30,31\}$, $\{31,1\}$, $\{32,2\}$, $\{3,33\}$, $\{4,34\}$, $\{5,33\}$, $\{6,34\}$, $\{33,34\}$\}.


{\lemma \label{lem:Blanusa34 designs}
There exist designs of order $52$ and $103$ for each of the two 34-vertex Blanu\v{s}a snarks.}\\


\noindent\textbf{Proof.}
Let the vertex set of $K_{52}$ be $Z_{52}$. The decomposition consists of

$(43,6,27,36,38,41,17,12,39,30,22,10,7,28,25,50,32,\adfsplit 26,21,5,8,48,47,45,3,44,24,29,14,13,46,15,40,11)_{\adfBca}$,

$(48,29,23,18,45,9,7,11,12,38,28,24,5,35,47,15,32,\adfsplit 39,26,36,33,19,0,34,14,50,25,17,27,13,21,46,4,42)_{\adfBca}$

\noindent and

$(6,20,37,23,45,44,26,11,29,28,27,35,48,50,10,39,46,\adfsplit 34,47,3,36,31,4,25,33,1,15,40,38,18,12,30,24,49)_{\adfBcb}$,

$(11,40,3,47,44,38,45,28,8,14,41,0,32,12,35,39,13,\adfsplit 6,26,7,21,34,30,17,51,5,49,31,33,42,1,36,10,9)_{\adfBcb}$

\noindent under the action of the mapping $x \mapsto x + 4$ (mod 52).

\adfVfy{52, \{\{52,13,4\}\}, -1, -1, -1} 

Let the vertex set of $K_{103}$ be $Z_{103}$. The decomposition consists of

$(0,1,2,3,6,8,14,15,4,26,39,18,5,41,21,35,40,\adfsplit 7,58,9,10,31,52,36,11,80,43,67,92,17,63,22,25,61)_{\adfBca}$

\noindent and

$(0,1,2,3,4,5,10,11,19,21,6,7,17,35,51,34,8,\adfsplit 9,27,28,48,50,12,72,38,14,69,95,13,43,66,55,37,76)_{\adfBcb}$

\noindent under the action of the mapping $x \mapsto x + 1$ (mod 103).\eproof

\adfVfy{103, \{\{103,103,1\}\}, -1, -1, -1} 

{\lemma \label{lem:Blanusa34 multipartite}
There exists decompositions of $K_{51,51,51}$ and $K_{17,17,17,17}$
into each of the two 34-vertex Blanu\v{s}a snarks.}\\


\noindent\textbf{Proof.}
Let the vertex set of $K_{51,51,51}$ be $Z_{153}$ partitioned according to residue classes modulo 3.
The decomposition consists of

$(0,1,2,4,7,11,24,15,5,29,45,25,8,3,31,54,6,\adfsplit 10,41,38,75,81,13,124,60,47,148,19,69,108,14,119,97,39)_{\adfBca}$

\noindent and

$(0,1,2,4,3,6,13,17,5,28,45,21,7,47,69,41,8,\adfsplit 9,37,31,12,62,97,68,11,114,55,109,125,32,73,51,67,110)_{\adfBcb}$

\noindent under the action of the mapping $x \mapsto x + 1$ (mod 153).

\adfVfy{153, \{\{153,153,1\}\}, -1, \{\{51,\{0,1,2\}\}\}, -1} 


Let the vertex set of $K_{17,17,17,17}$ be $Z_{68}$ partitioned according to residue classes modulo 4.
The decomposition consists of

$(17,53,40,42,13,43,60,39,29,0,35,51,46,10,27,33,41,\adfsplit 16,15,8,38,62,23,37,61,19,54,14,49,64,26,48,36,50)_{\adfBca}$,

$(46,35,8,39,55,1,30,64,13,15,2,62,29,20,27,31,49,\adfsplit 37,52,40,34,51,36,43,61,22,0,3,57,16,42,66,9,4)_{\adfBca}$

\noindent and

$(64,34,59,29,35,56,42,26,7,13,40,10,1,12,22,44,61,\adfsplit 37,62,52,32,46,19,31,21,48,2,55,15,17,54,41,28,47)_{\adfBcb}$,

$(2,50,44,13,4,17,7,18,25,39,28,20,30,33,12,35,1,\adfsplit 8,46,47,63,61,10,24,11,67,22,42,49,57,51,0,21,43)_{\adfBcb}$

\noindent under the action of the mapping $x \mapsto x + 4$ (mod 68).\eproof

\adfVfy{68, \{\{68,17,4\}\}, -1, \{\{17,\{0,1,2,3\}\}\}, -1} 

\vskip 2mm
Theorem~\ref{thm:Snark partial} (iv) follows from
Lemmas~\ref{lem:Blanusa34 designs} and \ref{lem:Blanusa34 multipartite}, and
Proposition~\ref{prop:e + 1, 2e + 1, e^3, (e/3)^4 => et + 1}.


\section{Zamfirescu's graph}
\label{sec:Zamfirescu}
\begin{figure}[h]
\begin{center}
\includegraphics[width=\adfGlarge]{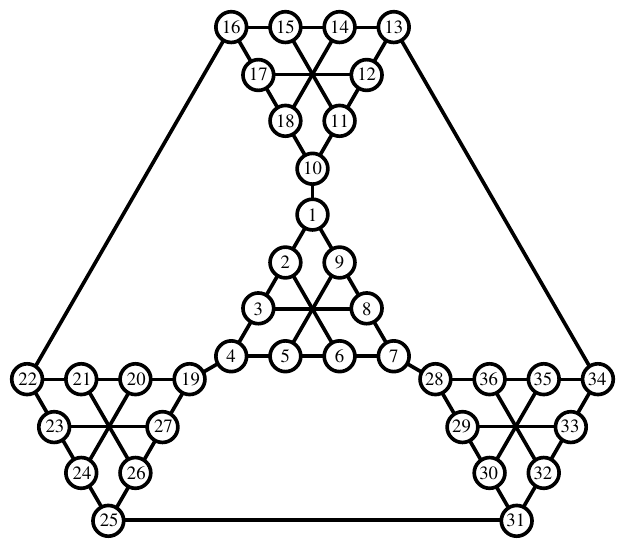}
\caption{Zamfirescu'g graph}
\end{center}
\end{figure}
\newcommand{\adfZ}{\mathrm{Z}}
Zamfirescu's graph is represented by the ordered 36-tuple of its vertices
(1, 2, \dots, 36)${}_{\adfZ}$ and edge set
  \{$\{1,2\}$, $\{2,3\}$, $\{3,4\}$, $\{4,5\}$, $\{5,6\}$, $\{6,7\}$, $\{7,8\}$, $\{8,9\}$, $\{9,1\}$, $\{3,8\}$, $\{6,2\}$, $\{9,5\}$,
   $\{10,11\}$, $\{11,12\}$, $\{12,13\}$, $\{13,14\}$, $\{14,15\}$, $\{15,16\}$, $\{16,17\}$, $\{17,18\}$, $\{18,10\}$, $\{12,17\}$, $\{15,11\}$, $\{18,14\}$,
   $\{19,20\}$, $\{20,21\}$, $\{21,22\}$, $\{22,23\}$, $\{23,24\}$, $\{24,25\}$, $\{25,26\}$, $\{26,27\}$, $\{27,19\}$, $\{21,26\}$, $\{24,20\}$, $\{27,23\}$,
   $\{28,29\}$, $\{29,30\}$, $\{30,31\}$, $\{31,32\}$, $\{32,33\}$, $\{33,34\}$, $\{34,35\}$, $\{35,36\}$, $\{36,28\}$, $\{30,35\}$, $\{33,29\}$, $\{36,32\}$,
   $\{1,10\}$, $\{4,19\}$, $\{7,28\}$, $\{16,22\}$, $\{25,31\}$, $\{34,13\}$\}
as supplied with the {\sc Mathematica} system, \cite{Math-Snarks}.


{\lemma \label{lem:Zamfirescu designs}
There exist Zamfirescu's graph designs of order $109$ and $217$.}\\


\noindent\textbf{Proof.}
Let the vertex set of $K_{109}$ be $Z_{109}$. The decomposition consists of

$(0,1,3,6,2,7,14,22,13,10,24,4,16,31,8,25,43,64,\adfsplit 28,52,5,50,9,78,11,39,66,45,80,17,83,51,15,76,46,104)_{\adfZ}$

\noindent under the action of the mapping $x \mapsto x + 1$ (mod 109).

\adfVfy{109, \{\{109,109,1\}\}, -1, -1, -1} 

Let the vertex set of $K_{217}$ be $Z_{217}$. The decomposition consists of

$(193,8,55,59,182,146,209,22,98,184,20,149,178,74,122,79,40,102,\adfsplit 134,99,85,68,26,3,119,156,153,141,126,164,169,86,172,73,133,16)_{\adfZ}$,

$(153,23,75,65,0,1,3,9,16,5,13,25,4,17,37,10,50,94,\adfsplit 6,24,58,84,11,69,2,112,167,54,95,14,128,35,180,61,125,174)_{\adfZ}$

\noindent under the action of the mapping $x \mapsto x + 1$ (mod 217).\eproof

\adfVfy{217, \{\{217,217,1\}\}, -1, -1, -1} 

{\lemma \label{lem:Zamfirescu multipartite}
There exists a decomposition of $K_{18,18,18}$ into Zamfirescu's graph.}\\


\noindent\textbf{Proof.}
Let the vertex set of $K_{18,18,18}$ be $Z_{54}$ partitioned according to residue classes modulo 3.
The decomposition consists of

$(21,25,35,40,14,6,7,45,31,28,30,29,15,47,13,2,43,0,\adfsplit 9,4,5,37,12,8,16,48,32,23,46,53,24,50,33,20,18,52)_{\adfZ}$

\noindent under the action of the mapping $x \mapsto x + 3$ (mod 54).\eproof

\adfVfy{54, \{\{54,18,3\}\}, -1, \{\{18,\{0,1,2\}\}\}, -1} 

\vskip 2mm
Theorem~\ref{thm:Snark partial} (v) follows from
Lemmas~\ref{lem:Zamfirescu designs} and \ref{lem:Zamfirescu multipartite}, and
Proposition~\ref{prop:2e + 1, 4e + 1, (e/3)^3 => 2et + 1}.


\section{Goldberg's snark \#5}
\label{sec:Goldberg5}
\begin{figure}[h]
\begin{center}
\includegraphics[width=\adfGlarge]{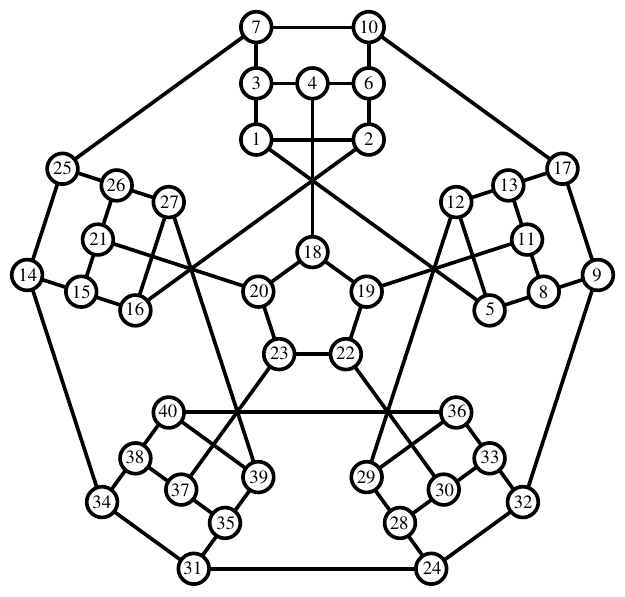}
\caption{Goldberg's snark \#5}
\end{center}
\end{figure}
\newcommand{\adfGe}{\mathrm{GS5}}
Goldberg's snark \#5 is represented by the ordered 40-tuple of its vertices
(1, 2, \dots, 40)${}_{\adfGe}$. The edge set is
  \{$\{15,16\}$, $\{14,15\}$, $\{15,21\}$, $\{21,26\}$, $\{25,26\}$, $\{26,27\}$, $\{16,27\}$, $\{14,25\}$, $\{7,10\}$, $\{2,6\}$, $\{6,10\}$, $\{4,6\}$,
   $\{3,4\}$, $\{3,7\}$, $\{1,3\}$, $\{1,2\}$, $\{9,17\}$, $\{5,8\}$, $\{8,9\}$, $\{8,11\}$, $\{11,13\}$, $\{13,17\}$, $\{12,13\}$, $\{5,12\}$, $\{24,32\}$,
   $\{28,29\}$, $\{24,28\}$, $\{28,30\}$, $\{30,33\}$, $\{32,33\}$, $\{33,36\}$, $\{29,36\}$, $\{31,34\}$, $\{38,40\}$, $\{34,38\}$, $\{37,38\}$,
   $\{35,37\}$, $\{31,35\}$, $\{35,39\}$, $\{39,40\}$, $\{2,16\}$, $\{7,25\}$, $\{10,17\}$, $\{1,5\}$, $\{9,32\}$, $\{12,29\}$, $\{24,31\}$, $\{36,40\}$,
   $\{14,34\}$, $\{27,39\}$, $\{20,21\}$, $\{4,18\}$, $\{11,19\}$, $\{22,30\}$, $\{23,37\}$, $\{19,22\}$, $\{22,23\}$, $\{20,23\}$, $\{18,20\}$, $\{18,19\}$\},
as supplied with the {\sc Mathematica} system, \cite{Math-Snarks}.

{\lemma \label{lem:Goldberg5 designs}
There exist Goldberg's snark \#5 designs of order $40$, $121$ and $241$.}\\


\noindent\textbf{Proof.}
Let the vertex set of $K_{40}$ be $Z_{39} \cup \{\infty\}$. The decomposition consists of the graphs

$(\infty,23,18,26,4,35,2,25,21,7,19,36,34,20,\adfsplit 17,3,24,30,0,31,1,33,6,38,22,8,29,\adfsplit 16,13,28,14,32,15,10,9,5,37,11,27,12)_{\adfGe}$

\noindent under the action of the mapping $\infty \mapsto \infty$, $x \mapsto x + 3$ (mod 39).

\adfVfy{40, \{\{39,13,3\},\{1,1,1\}\}, 39, -1, -1} 

Let the vertex set of $K_{121}$ be $Z_{121}$. The decomposition consists of the graphs

$(0,1,2,5,4,10,8,11,3,20,22,17,36,\adfsplit 6,21,37,53,23,43,45,68,9,69,7,33,96,\adfsplit 66,38,70,64,44,55,109,86,83,32,118,40,26,91)_{\adfGe}$

\noindent under the action of the mapping $x \mapsto x + 1$ (mod 121).

\adfVfy{121, \{\{121,121,1\}\}, -1, -1, -1} 

Let the vertex set of $K_{241}$ be $Z_{241}$. The decomposition consists of the graphs

$(169,233,1,108,6,57,27,190,239,224,144,192,130,\adfsplit 39,187,216,202,74,134,142,84,158,178,32,136,127,\adfsplit 56,212,184,79,5,133,83,154,3,64,55,66,69,104)_{\adfGe}$,

$(195,23,31,224,164,219,38,231,190,216,118,45,22,\adfsplit 96,176,170,92,13,200,12,24,1,33,0,5,80,\adfsplit 41,15,98,91,87,86,232,198,172,214,147,88,125,230)_{\adfGe}$

\noindent under the action of the mapping $x \mapsto x + 1$ (mod 241).\eproof

\adfVfy{241, \{\{241,241,1\}\}, -1, -1, -1} 

{\lemma \label{lem:Goldberg5 multipartite}
There exist decompositions of the complete multipartite graphs
$K_{60,60,60}$, $K_{20,20,20,20}$ and $K_{40,40,40,40}$ into Goldberg's snark \#5.}\\


\noindent\textbf{Proof.}
Let the vertex set of $K_{60,60,60}$ be $Z_{180}$ partitioned according to residue classes modulo 3.
The decomposition consists of the graph

$(0,1,2,6,5,14,9,15,4,28,31,22,51,\adfsplit 3,25,48,77,34,65,59,96,7,99,8,41,139,\adfsplit 83,49,101,93,58,78,173,122,113,55,10,63,16,131)_{\adfGe}$

\noindent under the action of the mapping $x \mapsto x + 1$ (mod 180).

\adfVfy{180, \{\{180,180,1\}\}, -1, \{\{60,\{0,1,2\}\}\}, -1} 

Let the vertex set of $K_{20,20,20,20}$ be $Z_{80}$ partitioned according to residue classes modulo 4.
The decomposition consists of the graph

$(31,48,56,18,1,5,29,50,41,71,43,47,38,\adfsplit 14,35,58,72,4,62,51,57,63,40,0,3,6,\adfsplit 8,19,12,60,74,59,42,9,27,55,25,64,37,10)_{\adfGe}$

\noindent under the action of the mapping $x \mapsto x + 2$ (mod 80).

\adfVfy{80, \{\{80,40,2\}\}, -1, \{\{20,\{0,1,2,3\}\}\}, -1} 

Let the vertex set of $K_{40,40,40,40}$ be $Z_{160}$ partitioned according to residue classes modulo 4.
The decomposition consists of the graph

$(37,98,22,49,88,40,148,41,36,93,2,94,129,\adfsplit 78,15,8,6,59,21,109,150,139,24,0,1,4,\adfsplit 25,11,12,33,18,62,128,149,104,61,47,106,23,152)_{\adfGe}$

\noindent under the action of the mapping $x \mapsto x + 1$ (mod 160).\eproof

\adfVfy{160, \{\{160,160,1\}\}, -1, \{\{40,\{0,1,2,3\}\}\}, -1} 

{\lemma \label{lem:Goldberg5 40 mod 120}
There exist a Goldberg's snark \#5 design of order $n$ if $n \equiv 40$ \textup{(mod 120)}.}\\

\noindent\textbf{Proof.}
There exist a 4-GDD of type $2^{3t+1}$ for $t \ge 2$, \cite{BSH} (see also \cite{Ge}).
Inflate by a factor of 20, lay a complete graph $K_{40}$ from Lemma~\ref{lem:Goldberg5 designs}
on each of the inflated groups and replace each block by a complete 4-partite graph
$K_{20,20,20,20}$ from Lemma~\ref{lem:Goldberg5 multipartite} to yield a
Goldberg's snark \#5 design of order $120t + 40$ for $t \ge 2$.
The design for the remaining case, $n = 160$ is constructed from $K_{40}$ (Lemma~\ref{lem:Goldberg5 designs})
and the 4-partite graph $K_{40,40,40,40}$ from Lemma~\ref{lem:Goldberg5 multipartite}.\eproof


\vskip 2mm
Theorem~\ref{thm:Snark partial} (vi) follows from
Lemmas~\ref{lem:Goldberg5 designs} and \ref{lem:Goldberg5 multipartite},
and Proposition~\ref{prop:2e + 1, 4e + 1, e^3 => 2et + 1} for $n \equiv 1$ (mod 120),
and Lemma~\ref{lem:Goldberg5 40 mod 120} for $n \equiv 40$ (mod 120).


\section{The Szekeres and Watkins snarks}
\label{sec:SzekeresWatkins}
\label{sec:last}
\begin{figure}[h]
\begin{center}
\includegraphics[width=\adfGlarge]{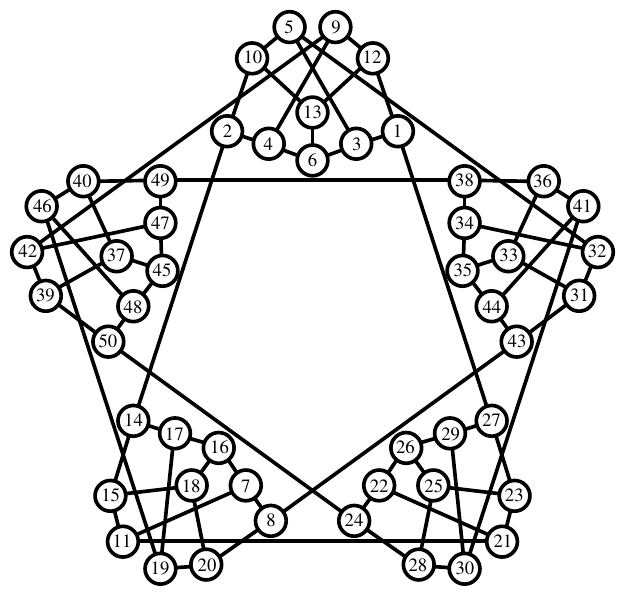}
\caption{The Szekeres snark}
\end{center}
\end{figure}
\begin{figure}[h]
\begin{center}
\includegraphics[width=\adfGlarge]{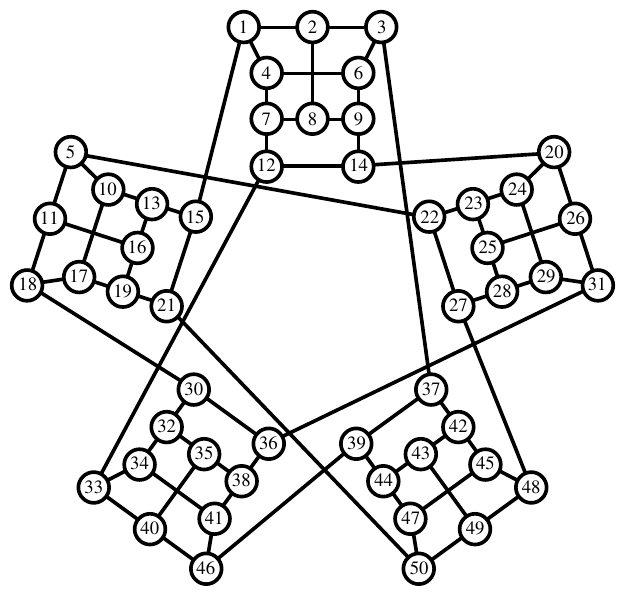}
\caption{The Watkins snark}
\end{center}
\end{figure}
\newcommand{\adfSze}{\mathrm{Sze}}
\newcommand{\adfWat}{\mathrm{Wat}}
The Szekeres snark and the Watkins snark are each represented by the ordered 50-tuple of its vertices,
(1, 2, \dots, 50)${}_{\adfSze}$ for the Szekeres snark and
(1, 2, \dots, 50)${}_{\adfWat}$ for the Watkins snark.
The edge sets, as supplied with the {\sc Mathematica} system, \cite{Math-Snarks}, are respectively

$\adfSze$:
  \{$\{29,30\}$, $\{27,29\}$, $\{22,24\}$, $\{21,22\}$, $\{21,23\}$, $\{23,25\}$, $\{25,28\}$, $\{28,30\}$, $\{22,26\}$, $\{25,26\}$,
   $\{26,29\}$, $\{23,27\}$, $\{24,28\}$, $\{3,6\}$, $\{5,10\}$, $\{1,3\}$, $\{3,5\}$, $\{1,12\}$, $\{2,10\}$, $\{4,6\}$, $\{4,9\}$, $\{9,12\}$,
   $\{12,13\}$, $\{6,13\}$, $\{10,13\}$, $\{2,4\}$, $\{16,17\}$, $\{19,20\}$, $\{14,17\}$, $\{17,19\}$, $\{14,15\}$, $\{8,20\}$, $\{7,16\}$, $\{7,11\}$,
   $\{11,15\}$, $\{15,18\}$, $\{16,18\}$, $\{18,20\}$, $\{7,8\}$, $\{35,44\}$, $\{36,41\}$, $\{43,44\}$, $\{41,44\}$, $\{31,43\}$, $\{36,38\}$,
   $\{34,35\}$, $\{32,34\}$, $\{31,32\}$, $\{31,33\}$, $\{33,35\}$, $\{33,36\}$, $\{34,38\}$, $\{45,47\}$, $\{39,42\}$, $\{47,49\}$, $\{42,47\}$,
   $\{40,49\}$, $\{39,50\}$, $\{45,48\}$, $\{46,48\}$, $\{40,46\}$, $\{37,40\}$, $\{37,45\}$, $\{37,39\}$, $\{48,50\}$, $\{30,41\}$, $\{9,42\}$,
   $\{11,21\}$, $\{5,32\}$, $\{19,46\}$, $\{1,27\}$, $\{2,14\}$, $\{8,43\}$, $\{38,49\}$, $\{24,50\}$\}
and

$\adfWat$:
  \{$\{1,2\}$, $\{1,4\}$, $\{1,15\}$, $\{2,3\}$, $\{2,8\}$, $\{3,6\}$, $\{3,37\}$, $\{4,6\}$, $\{4,7\}$, $\{5,10\}$, $\{5,11\}$, $\{5,22\}$,
   $\{6,9\}$, $\{7,8\}$, $\{7,12\}$, $\{8,9\}$, $\{9,14\}$, $\{10,13\}$, $\{10,17\}$, $\{11,16\}$, $\{11,18\}$, $\{12,14\}$, $\{12,33\}$,
   $\{13,15\}$, $\{13,16\}$, $\{14,20\}$, $\{15,21\}$, $\{16,19\}$, $\{17,18\}$, $\{17,19\}$, $\{18,30\}$, $\{19,21\}$, $\{20,24\}$,
   $\{20,26\}$, $\{21,50\}$, $\{22,23\}$, $\{22,27\}$, $\{23,24\}$, $\{23,25\}$, $\{24,29\}$, $\{25,26\}$, $\{25,28\}$, $\{26,31\}$,
   $\{27,28\}$, $\{27,48\}$, $\{28,29\}$, $\{29,31\}$, $\{30,32\}$, $\{30,36\}$, $\{31,36\}$, $\{32,34\}$, $\{32,35\}$, $\{33,34\}$,
   $\{33,40\}$, $\{34,41\}$, $\{35,38\}$, $\{35,40\}$, $\{36,38\}$, $\{37,39\}$, $\{37,42\}$, $\{38,41\}$, $\{39,44\}$, $\{39,46\}$,
   $\{40,46\}$, $\{41,46\}$, $\{42,43\}$, $\{42,45\}$, $\{43,44\}$, $\{43,49\}$, $\{44,47\}$, $\{45,47\}$, $\{45,48\}$, $\{47,50\}$,
   $\{48,49\}$, $\{49,50\}$\}.


{\lemma \label{lem:SzekeresWatkins designs}
There exist Szekeres snark and Watkins snark designs of orders $76$ and $151$.}\\


\noindent\textbf{Proof.}
Let the vertex set of $K_{76}$ be $Z_{76}$. The decompositions consist of

$(54,47,25,16,39,41,38,22,48,32,57,58,28,14,69,59,31,\adfsplit 27,45,46,36,37,65,26,20,44,42,4,19,23,70,12,1,67,\adfsplit 18,10,62,11,0,29,40,35,7,60,21,64,6,66,51,74)_{\adfSze}$,

$(22,4,66,47,59,41,44,50,53,70,8,56,65,74,3,58,61,\adfsplit 67,28,1,20,12,38,23,0,17,2,26,43,51,9,49,63,64,\adfsplit 16,5,29,42,25,57,35,33,11,69,73,27,71,36,68,75)_{\adfSze}$

\noindent and

$(35,55,62,65,49,22,17,33,21,24,26,72,60,40,6,71,38,\adfsplit 47,48,36,45,4,7,8,67,1,23,15,34,16,64,10,29,11,\adfsplit 59,0,18,27,66,9,70,32,39,68,42,14,46,69,3,44)_{\adfWat}$,

$(48,22,57,28,15,19,71,10,4,40,43,74,27,63,41,49,64,\adfsplit 30,18,68,1,52,29,26,24,9,61,32,5,37,11,51,56,17,\adfsplit 39,13,53,21,6,35,2,3,42,47,54,69,70,72,34,45)_{\adfWat}$

\noindent under the action of the mapping $x \mapsto x + 4$ (mod 76).

\adfVfy{76, \{\{76,19,4\}\}, -1, -1, -1} 

Let the vertex set of $K_{151}$ be $Z_{151}$. The decompositions consist of

$(47,54,128,27,56,15,131,51,71,99,137,21,5,23,3,29,82,\adfsplit 92,24,93,129,18,132,102,104,123,33,19,110,80,105,103,100,43,\adfsplit 25,121,63,97,88,72,7,49,10,14,118,57,84,89,120,2)_{\adfSze}$

\noindent and

$(93,6,13,146,107,110,42,144,149,46,26,61,79,140,94,100,20,\adfsplit 3,75,67,97,1,134,33,27,11,72,112,21,130,120,116,53,49,\adfsplit 127,92,23,129,71,121,8,50,85,128,138,39,63,118,56,5)_{\adfWat}$

\noindent under the action of the mapping $x \mapsto x + 1$ (mod 151).\eproof

\adfVfy{151, \{\{151,151,1\}\}, -1, -1, -1} 

{\lemma \label{lem:SzekeresWatkins multipartite}
There exist decompositions of $K_{75,75,75}$ and $K_{25,25,25,25}$ into each of the
Szekeres snark and the Watkins snark.}\\


\noindent\textbf{Proof.}
Let the vertex set of $K_{75,75,75}$ be $Z_{225}$ partitioned according to residue classes modulo 3.
The decompositions consist of

$(186,132,25,92,162,42,89,205,12,214,45,28,161,100,119,34,96,\adfsplit 216,175,26,152,142,157,18,36,8,101,29,30,88,217,209,192,144,\adfsplit 143,53,61,22,75,134,165,49,176,82,81,39,47,196,9,173)_{\adfSze}$

\noindent and

$(204,53,99,109,0,161,186,139,48,4,121,185,197,70,16,141,59,\adfsplit 201,49,32,30,152,223,66,165,1,205,182,7,203,89,45,87,94,\adfsplit 169,208,217,143,132,209,102,17,93,23,148,145,12,177,193,158)_{\adfWat}$

\noindent under the action of the mapping $x \mapsto x + 1$ (mod 225).

\adfVfy{225, \{\{225,225,1\}\}, -1, \{\{75,\{0,1,2\}\}\}, -1} 

Let the vertex set of $K_{25,25,25,25}$ be $Z_{100}$ partitioned according to residue classes modulo 4.
The decompositions consist of

$(32,79,37,88,71,75,81,72,53,0,14,67,50,34,48,68,45,\adfsplit 87,66,57,85,43,4,77,15,78,10,30,20,23,8,25,1,7,\adfsplit 56,63,12,54,94,65,80,52,35,13,11,96,62,21,83,27)_{\adfSze}$,

$(12,6,87,11,40,56,72,1,98,37,79,61,86,33,18,90,96,\adfsplit 92,55,14,81,67,42,52,43,73,85,54,27,49,26,7,95,97,\adfsplit 38,0,62,74,28,35,23,17,80,21,65,58,10,4,9,66)_{\adfSze}$

\noindent and

$(10,53,36,1,3,59,46,67,96,17,14,75,23,45,25,64,43,\adfsplit 89,74,32,84,66,0,6,51,85,48,9,91,18,80,27,44,82,\adfsplit 37,49,79,30,86,83,63,12,54,11,87,13,93,26,92,71)_{\adfWat}$,

$(75,26,5,40,92,67,31,9,14,13,81,34,47,1,70,0,23,\adfsplit 82,58,42,93,87,74,44,91,72,4,18,29,21,66,84,88,85,\adfsplit 41,12,7,15,8,90,60,6,20,53,24,37,50,97,35,28)_{\adfWat}$

\noindent under the action of the mapping $x \mapsto x + 4$ (mod 100).\eproof

\adfVfy{100, \{\{100,25,4\}\}, -1, \{\{25,\{0,1,2,3\}\}\}, -1} 

\vskip 2mm
Theorem~\ref{thm:Snark partial} (vii) follows from
Lemmas~\ref{lem:SzekeresWatkins designs} and \ref{lem:SzekeresWatkins multipartite}, and
Proposition~\ref{prop:e + 1, 2e + 1, e^3, (e/3)^4 => et + 1}.
The proof of Theorem~\ref{thm:Snark partial} is complete.



\end{document}